\documentclass[reqno,11pt]{amsart}
\usepackage{geometry}
\usepackage[numbers,sort&compress]{natbib}
\usepackage{mathtools,amssymb,amsthm,mathrsfs,color,lineno,paralist,graphicx,float}
\usepackage[colorlinks,
linkcolor=red,
anchorcolor=green,
citecolor=blue,
]{hyperref}
\usepackage{xcolor}

\usepackage[T1]{fontenc}
\usepackage[utf8]{inputenc}
\usepackage{rotating}
\usepackage{pxfonts}

\usepackage{rotating}
\usepackage{tikz}
\usetikzlibrary{positioning}
\usepackage{calc}
\definecolor{bleu1}{RGB}{0,57,128}
\def\bleu1{\color{bleu1}}

\usepackage{etoolbox}
\makeatletter
\patchcmd{\@maketitle}{\topskip42\p@}{\topskip12\p@}{}{}
\makeatother
\makeatletter
\patchcmd{\@setauthors}{\@topsep30\p@}{\@topsep18\p@}{}{}
\patchcmd{\@setabstracta}{\skip@20\p@}{\skip@12\p@}{}{}
\makeatother
\patchcmd{\section}{\normalfont}{\normalfont \bleu1}{}{}
\patchcmd{\subsection}{\normalfont}{\normalfont \bleu1}{}{}
\patchcmd{\subsubsection}{\normalfont}{\normalfont \bleu1}{}{}
\renewcommand{\proofname}{\it \bleu1 Proof}

\newcommand{\SL}{\mathrm{SL}(2,\mathbb{R})}

\def\a{\alpha}

\def\e{\varepsilon}

\let\newpf\proof \let\proof\relax 
\newenvironment{pf}{\newpf[\proofname]}{\qed\endtrivlist}

\renewcommand{\emph}[1]{\textit{#1}}

\newcommand{\ba}{\overline{A}}

\def\be{\begin{equation}}
\def\ee{\end{equation}}

\def\ba{{\begin{align}}}
\def\ea{{\end{align}}}

\def\bm{\begin{matrix}}
\def\em{\end{matrix}}

\def\u{{\mathbb U}}

\def\a{{\alpha}}

\def\SL{{\mathrm{SL}}}

\def\0{{\mathbf 0}}

\newtheorem{Theorem}{Theorem}[section]
\newtheorem{Lemma}{Lemma}[section]
\newtheorem{Proposition}{Proposition}[section]
\newtheorem{Corollary}{Corollary}[section]
\newtheorem{Remark}{Remark}[section]

\newtheorem{Definition}{Definition}[section]

\numberwithin{equation}{section}

\theoremstyle{definition}

\newtheorem{definition}{Definition}[section]

\def\ssm{\smallsetminus}

\renewcommand{\setminus}{\ssm}

\newcommand{\dist}{\operatorname{dist}}

\renewcommand{\mod}{\operatorname{mod}}

\newcommand{\C}{{\mathbb C}}

\newcommand{\N}{{\mathbb N}}
\newcommand{\Q}{{\mathbb Q}}
\newcommand{\R}{{\mathbb R}}
\newcommand{\T}{{\mathbb T}}

\newcommand{\Z}{{\mathbb Z}}

\def\B0{{\bold{0}}}

\catcode`\@=12

\def\Empty{}
\newcommand\oplabel[1]{
  \def\OpArg{#1} \ifx \OpArg\Empty {} \else
    \label{#1}
  \fi}

\newcommand{\comm}[1]{}
\newcommand{\comment}[1]{}

\begin{document}

\title{Intrinsic Symplectic Structure and Sharp Arithmetic Universality}
\thanks{This paper is a substantially revised version of the preprint previously posted under the title “Hidden subcriticality, symplectic structure, and universality of sharp arithmetic spectral results for type I operators” which is superseded by the present version and is not intended for publication. The present version substantially revises the exposition, introduces new concepts clarifying the main results, rewrites parts of the proofs, reformulates a central part of the argument in arbitrary dimension 
$k$ rather than only dimension 
$2$ and presents the structural framework in a form encompassing both subcritical and supercritical regimes; the proofs of the main universality results are essentially unchanged.}
\author{Lingrui Ge}
\address{Beijing International Center for Mathematical Research, Peking University, Beijing, China}
\email{gelingrui@bicmr.pku.edu.cn}

\author{Svetlana Jitomirskaya}
\address{Department of Mathematics, University of California, Berkeley, CA 94720, USA}
\email{sjitomi@berkeley.edu}

\begin{abstract}

We show that  formal eigenvalue equations of the Aubry duals of analytic one-frequency Schr\"odinger operators admit
intrinsic analytic
$Sp(2k,\C)$  structures, where  $k=k(E)$ is the T-acceleration in global theory.

 For
trigonometric potentials those structures govern the center dynamics of partially
hyperbolic dual cocycles; for general analytic potentials they persist, without loss of
analyticity, as an intrinsic object even when the dual operator has infinite range and no finite-dimensional  cocycle exists.

For $k=1$, we also introduce the concept of projectively real cocycles: Hermitian-symplectic
systems whose projective action is algebraically conjugate, up to a scalar phase, to that
of a real $\SL(2,\R)$ cocycle. This allows us to define a rotation pair and establish a
rotation--IDS correspondence in the general analytic setting, where standard dynamical methods fail.

Using these tools, we solve two 
spectral
        arithmetic conjectures: universality of the sharp arithmetic transition in frequency
(AAJ) and of the absolute continuity of the integrated density of states for all
frequencies, throughout the class of non-critical Type I operators, an
open and conjecturally dense set. We also prove
universality of sharp $1/2$-H\"older continuity of the integrated density of states for
non-critical
Type I operators with Diophantine frequencies, establishing part of You's
conjecture.

The new route to sharp spectral results for general
analytic potentials through intrinsic symplectic centers, introduced in this paper, underlies
the robust and supercritical dry Ten Martini theorems for Type~I
operators \cite{gjy1,lxz} and, in higher center dimensions, the general
supercritical Ten Martini theorem \cite{gjy2}. The higher-dimensional construction may also lead to further spectral
applications beyond the Type~I regime.
       
\end{abstract}

\maketitle
\vspace{-2em}
\tableofcontents

\section{Introduction}
\subsection{Overview 
}
\label{subsec:intro-overview}

The primary goal of this paper is to establish \emph{universality of sharp arithmetic spectral phenomena}
for analytic one-frequency Schrödinger operators
\begin{equation}\label{sch} (H_{v, \alpha, \theta}u)_n = u_{n+1} +
  u_{n-1} + v(\theta + n\alpha)u_n, \quad n \in \mathbb{Z},
  \theta\in\T, v\in C^\omega(\T,\R)\end{equation}
extending celebrated almost Mathieu ($v(\theta)=2\lambda\cos 2\pi \theta$)
results to a broad and robust class of analytic potentials.

The almost Mathieu operator
has played a central role in the modern theory of operators
\eqref{sch}, both as a foundational solid-state physics model
\cite{Peierls} 
and as a
mathematical prototype. From early on, its study has inspired new methods in many fields, including dynamical
systems, semiclassical analysis,
and multiscale analysis (e.g. \cite{H, HelSj,fsw}). It has served as the primary testing
ground for the interaction between quasiperiodicity, spectral theory,
and arithmetic properties of the frequency.

Historically, the development of the
spectral theory of operators \eqref{sch} has been shaped by the almost Mathieu
      problems on Simon's lists \cite{sim15,simXXI}, all by now fully
      solved \cite{J,aj,ayz,JLiu, ak, jkras, ad,bjk}.

      At the same time, the almost Mathieu operator is  as central and iconic, as it is highly
   special, with all celebrated almost Mathieu \emph{proofs} utilizing its
   specific features---most notably its reflection symmetry
 (evenness) and self-duality --- that do not withstand small analytic
 perturbations. Yet its physical origin and
      relevance certainly suggest that the
      properties of non-critical
      almost Mathieu operators should be robust with respect to small
   analytic perturbations of $v$. 
   
A central feature of operators \eqref{sch} is that
for irrational $\alpha$ there is an interplay between their
spectral properties and approximation arithmetics of $\alpha. $ For
the almost Mathieu operator such arithmetic (non)transitions are
\emph{sharp}. The corresponding
results---such as the AAJ  transition in frequency, the Ten Martini
Problem, and absolute continuity
of the integrated density of states (IDS) for all irrational frequencies---  have long been regarded as among the most rigid and delicate
features of the almost Mathieu operator and 
resisted extension.

The dependence of the existing proofs of these sharp results on the special properties of the almost Mathieu is structural.
Reflection symmetry leads, in particular, to the real $SL(2,\mathbb R)$
dual dynamics required by the classical rotation-number
and reducibility arguments, while finite range supplies
finite-dimensional dynamics.
General analytic perturbations do not preserve either,
 so these methods cannot be applied directly.

In this work we show that these sharp arithmetic phenomena are not
exceptional, nor symmetry based, but rather manifestations of
a deeper geometric structure present in general analytic
settings.

For trigonometric polynomials $v,$ 
duals of operators \eqref{sch} are finite
range, and their eigenvalue equations naturally give rise to symplectic
dynamics. However, this structure 
is inherently tied to the dual cocycle that changes wildly with small
perturbations. A central discovery of \cite{gjyz} was that for
trigonometric polynomials $v$ the dual cocycle is partially hyperbolic with center
of dimension $2k$, where $k$ is the
 \emph{T-acceleration} in global theory. 
 
Here we show that the symplectic structure on this center is intrinsic to the dual
eigenvalue equation. We prove convergence of the center dynamics under trigonomet-
ric polynomial approximation, \emph{without any loss of analyticity}. This defines an intrinsic
$Sp(2k,\C)$ structure for general analytic potentials, even when the full dual equation has
 no finite-dimensional transfer-matrix formulation, thus resolving the fundamental obstruction
that previously prevented duality-based arguments from extending
beyond finite-range settings.

For $k=1$, the intrinsic symplectic center has an additional
structure. We introduce the notion of \emph{projectively real
cocycles} and show that, in a suitable analytic frame, the center
cocycle is a scalar phase times a real $SL(2,\mathbb R)$ cocycle.
This uncovers the real dynamics needed for rotation theory
without reflection symmetry, resolving a longstanding obstacle
to extending sharp arithmetic results to general analytic
potentials.

Utilizing projective reality, we introduce a notion of fibered rotation pair for the two-dimensional
center and establish its relation to the IDS. In this
framework, the integrated density of states is expressed in terms of rotation data associated with the
intrinsic $Sp(2,\mathbb C)$ structure, without reflection symmetry and even in the absence of  finite-range dual operator. This correspondence
transfers absolute continuity and sharp H\"older regularity from
rotation numbers to the IDS and is also used in the
reducibility-to-localization argument.

Our spectral results apply to the class of \emph{Type~I} energies,
introduced in \cite{gjy}\footnote{The earlier preprint \cite{gjy} is superseded by \cite{gjy1}.} and characterized by a simple acceleration pattern in
Avila’s global theory. Within this class we prove universality of the major sharp arithmetic spectral
features previously established in the almost Mathieu setting. From this perspective, Type~I emerges as the
natural \emph{universality class} for sharp arithmetic phenomena in the supercritical regime.

The Type~I condition is open in each $C^\omega_h$ whose strip extends past the relevant first turning points \cite{gjy1}, so our
results are analytically robust. Moreover, it goes far beyond just  a
neighborhood of the existing symmetric models. In fact, a natural
conjecture is that \emph{Type I is generic}, i.e. that
\emph{Type~I energies are (open and) dense in the spectrum} for generic (i.e. open and
dense) analytic one-frequency Schrödinger operators. We discuss
supporting evidence in Appendix \ref{app}, see also Conjecture 3.2 in \cite{you}.

Since the first version of this paper, our two-dimensional center construction was also
used for several other outstanding results for the entire class of {\it Type I} operators:  the  ten martini theorem
\cite{gjy1}, the dry ten martini theorem for supercritical case
\cite{lxz}, and  the phase counterpart of our sharp arithmetic phase transitions result (i.e. the extension of the almost Mathieu result of \cite{JLiu1})
and exact exponential gap asymptotics for supercritical 
case, without assuming evenness (where the projectively real center and its subcriticality were also crucial) \cite{flz}. The higher dimensional center construction is a key ingredient in the general supercritical
Ten Martini theorem \cite{gjy2}.
\subsection{Main results}\label{subsec:intro-results}
We now state the main results of the paper; see precise definitions and complete formal
statements, when different, in Section 1.4-1.6.

The three main universality results apply to Type~I, a large and robust class of
energies/operators. 

\medskip
Let $L(E)$ be the Lyapunov exponent and $\beta(\alpha)$ denote the \emph{arithmetic exponent} associated with the frequency $\alpha$.
Their interplay establishes  the precise threshold between localization and singular continuous spectrum in the sharp
arithmetic transition, as described in what is known as the Aubry--André--Jitomirskaya
(AAJ) conjecture, originally formulated for the almost Mathieu operator \cite{jcong}: the spectrum is pure point for a.e. $\theta$ on $\{E : L(E)>\beta(\alpha)\}$ and purely
singular continuous for all $\theta$ on 
$
\{E : 0 < L(E) < \beta(\alpha)\}.
$
\begin{Theorem}\label{thm:intro-aaj}
For any real analytic $v$ and irrational $\alpha$, the sharp arithmetic transition in frequency (i.e. the AAJ) holds universally for Type
  I energies $E$ of operators \(H_{v,\alpha,\theta}\).\footnote{Our sharp arithmetic transition result was originally announced  at QMath15 (UC Davis, September 2022).}\end{Theorem}

The next result concerns the integrated density of states (IDS), $N(E)$.
Its absolute continuity  \emph{for all $\alpha$}
was previously known only for the almost Mathieu operator and long
conjectured to hold universally for \emph{non-critical} operators. We extend
this result to all non-critical Type~I operators.

\begin{Theorem}
 
  The absolute continuity of $N(E)$ for all $\alpha$ is universal for every non-critical Type
  I operator \(H_{v,\alpha,\theta}\).
\end{Theorem}

\medskip
The next result concerns sharp Hölder regularity. H\"older exponent
deteriorates as frequencies get more Liouville \cite{alsz}, but is
expected to be the same for a.e. (in fact, all Diophantine) $\alpha.$ Sharp \(1/2\)-Hölder
regularity for Diophantine $\alpha$ is indeed universal throughout
the subcritical regime \cite{aj1,avila2}. In the supercritical regime,
the modulus of continuity becomes also dependent on  the acceleration. According to You’s conjecture \cite{You,you}, for Diophantine
\(\alpha\)  and supercritical operators, \(1/2\)-Hölder regularity is expected \emph{only} at Type I energies. This is exactly what we prove.

\begin{Theorem}\label{thm:intro-holder}
For Diophantine frequencies \(\alpha\), the IDS at all non-critical Type~I energies is locally \(1/2\)-Hölder
continuous. This regularity is optimal in this class.
\end{Theorem}

\medskip

The following structural results provide the key framework underlying these universality theorems.

\medskip

\begin{Theorem}[Intrinsic symplectic center structure; informal]\label{thm:intro-symplectic}
Let $E$ be an energy in the spectrum of $H_{v, \alpha, \theta} $ and $k$ be the T-acceleration in the global theory. 

There exists a canonical analytic $Sp(2k, \mathbb{C})$ structure
intrinsic to the dual eigenvalue equation
$\widehat H_{v,\alpha,\theta}u=Eu$.
This structure governs the center dynamics if $v$ is a trigonometric
polynomial, and
persists, without any loss of analyticity, under passage to analytic $v$ where the dual operator has infinite range and no cocycle formulation exists.
\end{Theorem}

To the best of our knowledge, Theorem~\ref{thm:intro-symplectic} gives
the first construction of an intrinsic analytic finite-dimensional symplectic center dynamics
for  infinite-range ergodic operators, opening these operators to
symplectic dynamical methods despite the absence of dynamics for the
full eigenvalue equation.\footnote{Following earlier
versions of this work, Li--Wang--You--Zhou \cite{LWYZ}
developed a different operator-theoretic construction
of the center for small potentials over more general
base dynamics, recovering our intrinsic center bundle
in the common one-frequency regime.} In dimension two, this structure underlies the robust and
supercritical dry Ten Martini theorems for Type~I operators
\cite{gjy1,lxz}; in higher dimensions, it is a key ingredient
in the general supercritical Ten Martini theorem \cite{gjy2}.
We expect further applications.
\medskip

For Type I energies (i.e.,
\(k = 1\)), Theorem
\ref{thm:intro-symplectic} 
 effectively reduces the infinite-dimensional dual ``dynamics''  to two-dimensional ``center''. 
 For even \(v\),  the resulting center dynamics lie in \(SL(2,\mathbb{R})\), allowing standard methods to apply.
However, the presence of symmetry is highly nongeneric and unstable
under perturbations, reflecting the fact that it imposes infinitely
many independent constraints on the projective dynamics. When $v$ is
not even, the symmetry is broken: the dynamics lie in $Sp(2,\C)$, but
no longer in $\SL(2,\R)$. Thus the \emph{classical} real $\SL(2,\R)$
rotation-number--IDS and reducibility-to-localization pathways are no
longer directly available.

To address this, we introduce the concept of \emph{projectively real
  cocycles}, where the projective action, while genuinely complex, is
algebraically conjugate (up to a scalar phase) to a real
\(SL(2,\mathbb{R})\) cocycle.  For such cocycles, if their real factors are homotopic to identity, we can naturally
define \emph{rotation pairs} \((\rho_1, \rho_2)\), and it turns out
that the \(Sp(2,\mathbb{C})\) center cocycles of duals of supercritical Type I
operators are necessarily projectively real  with their  corresponding
\(SL(2,\mathbb{R})\) cocycles homotopic to the identity.  
 
The projectively real structure resolves the longstanding
obstacle to applying real-cocycle methods to the dual center
without reflection symmetry. We use it to establish the
rotation--IDS correspondence for general analytic potentials.
We expect the notion of projective reality to have further
applications in the dynamics of complex cocycles.

\begin{Theorem}[The Rotation-IDS correspondence]\label{thmids}
For supercritical Type I energies $E_0$ of $H_{v, \alpha, \theta}$ with
irrational $\alpha$ and $v\in C^\omega(\T,\R),$ there is a neighborhood $I$ of $E_0$ such that for any $E_1,E_2\in I$,
\begin{equation}\label{Nrho}
N(E_1)-N(E_2)=\rho_1(E_1)-\rho_2(E_1)-\rho_1(E_2)+\rho_2(E_2).
\end{equation} 
\end{Theorem}

\begin{Remark}\label{rem:rotation-normalization}
For continuous real lifts on interval $I$, the correspondence reads
\[
\rho_1(E)-\rho_2(E)=N(E)+c_I,
\]
with $c_I$ independent of $E$. If $I$ meets an open spectral gap, uniform
hyperbolicity of the center there and gap labelling imply
$c_I\in\mathbb Z+\alpha\mathbb Z$. A suitable choice of Hermitian-symplectic
frame and lifts then yields the absolute identity 
\begin{equation}\label{rhoN}\rho_1(E)-\rho_2(E)=N(E).\end{equation}
For Type I operators, the existence of such gaps
in every  interval follows from the robust ten martini theorem
\cite{gjy1}, so, invoking that,  \eqref{rhoN} holds.  The difference formulation used in Theorem \ref{thmids} avoids unnecessary dependence on that
result and suffices for all our applications.
\end{Remark}

Note that \eqref{Nrho} specializes to the classical $N(E_1)-N(E_2)=2(\widehat\rho(E_1)-\widehat\rho(E_2))$ in the symmetric (even) case, where $\widehat\rho$ denotes the rotation number of the (intrinsic) dual center cocycle.

These 
results underlie all proofs in the paper.

\subsection{Historical background and 
main difficulties}\label{subsec:intro-history}

\subsubsection{The almost Mathieu operator, Aubry duality, spectral transitions, and sharp arithmetic phenomena}

The central special feature underlying the spectral theory of the
almost Mathieu family is its invariance with respect to Aubry duality,
 a Fourier-type transform that exchanges the
subcritical and supercritical regimes, and relates the operator at
coupling $\lambda$ to the operator at coupling $\lambda^{-1}$. 
Aubry duality preserves the spectrum and the
IDS Aubry duality preserves the spectrum and the IDS
\cite{as83,fp84}. The
almost Mathieu operator is exceptional in that its structure allows
both regimes to be treated simultaneously,
 allowing for cooperation of subcritical and supercritical methods
 for related problems.

The action of Aubry duality on the spectral decomposition
is a lot more delicate. In the subcritical regime $|\lambda|<1$, the Lyapunov exponent
vanishes, which by the celebrated Kotani theory \cite{kotani} implies
presence  of ac spectrum and, by Deift--Simon \cite{deiftsimon}, $L^2$ generalized eigenfunctions. It is also the domain of
KAM-based reducibility methods, going back to \cite{ds,blt} and brought to
perturbative perfection by Eliasson in \cite{eli92}. They were made nonperturbative through
  duality-based conjugation \cite{bj2,aj1}, and ultimately led to purely absolutely
continuous spectrum for all $\alpha,\theta$ \cite{avila} and absolute continuity
\cite{ad} and sharp H\"older regularity \cite{aj1} of the integrated density of
states. Following the work of Fayad--Krikorian \cite{fk},
Avila--Fayad--Krikorian \cite{afk} developed a KAM scheme
for conjugacy to cocycles of rotations at arbitrary
irrational $\alpha$. For a different KAM approach
through quasiperiodic linear systems and its transfer
to cocycles, see \cite{hy,yz}.

In the supercritical regime $|\lambda|>1$, the Lyapunov exponent is
positive, leading to no ac spectrum. It is distinguishing between
 singular continuous and pure point spectrum in this regime that
 becomes arithmetic.

 The appearance of arithmetic effects in this context was itself a
gradual and conceptually nontrivial development.  The original
Aubry--André conjecture for the almost Mathieu operator \cite{aa} predicted a
sharp transition between absolutely continuous and pure point spectrum
at the self-dual point, governed solely by the vanishing or positivity
of the Lyapunov exponent so
it did not pay respect to the arithmetics of parameters.  It was soon realized, however \cite{as} that the
arithmetic properties of the frequency $\alpha$ can obstruct
localization in the regime of positive Lyapunov exponent.

It eventually became clear  that non-arithmetic AA conjecture was
correct in the subcritical regime $|\lambda|<1$ \cite{lastac,J,avila}, while in the supercritical one arithmetics does rule the game. The
sharp arithmetic transition in frequency was conjectured in \cite{jcong}
predicting a precise sharp dichotomy between localization and singular
continuous spectrum governed by the comparison of the Lyapunov exponent $L(E)$
with an arithmetic exponent $\beta(\alpha)$ that measures the approximation properties of
$\alpha$.  This conjecture, dubbed the
Aubry--André--Jitomirskaya (AAJ) conjecture in \cite{ayz}, crystallized the role of
arithmetic as a  sharp spectral mechanism rather than a
perturbative effect.

Another form of arithmetic rigidity appears when there are different fundamental
reasons for a phenomenon to hold on the Liouville and Diophantine
sides, yet those approaches can be combined to establish
an all-irrational-frequency statement as in \cite{cey,aj} or \cite{J,ad}.

Two other almost Mathieu problems of this kind were the Ten Martini Problem (Cantor
spectrum) and absolute continuity of the IDS, promoted, along with the
a.e. version of the AAJ, in
Simon's lists \cite{sim15,simXXI} as some of the central challenges in
the subject. For example, the all-frequency absolute continuity of the
IDS was established separately for the Diophantine case \cite{J} and the Liouville case
\cite{ad}. 

We note that both Ten Martini and AAJ were solved for
(arithmetically) almost all parameters \cite{Puig,J} before their
final celebrated all $\alpha$ solutions \cite{aj, ayz, JLiu}.   A similar distinction between almost-everywhere and all-parameter results appears in
 \cite{ak,avila,last,J}. 
This underscores the particular importance attached to {\it sharp
  arithmetic} results in the subject: they are treated as number theoretic problems, with no omission, approximation, or perturbative
loss allowed. The eventual resolution of those for the almost Mathieu operator
required new ideas capable of treating both Diophantine and highly
Liouville frequencies and represented  landmark achievements in
quasiperiodic spectral theory; see, e.g.,
\cite{J,ayz,JLiu,ak,jkras,ad} and references therein.  From a broader perspective, the almost Mathieu operator thus came to be
viewed as the canonical model for sharp arithmetic phenomena in
quasiperiodic spectral theory.  At the same time, this success reinforced a natural question: which aspects of these phenomena are truly universal, and which are  almost Mathieu specific?

\subsubsection{Avila’s global theory and the universality divide}\label{subsec:intro-global}

The prototypical role of the almost Mathieu operator becomes especially
transparent in Avila’s global theory \cite{avila0} of analytic
quasiperiodic Schr\"odinger operators.  In this framework, the spectrum
is divided into \emph{subcritical}, \emph{critical}, and
\emph{supercritical} regimes according to the behavior of the
complexified Lyapunov exponent, modeled on the corresponding regimes
$|\lambda|<1$, $|\lambda|=1$, and $|\lambda|>1$ of the almost Mathieu
operator.  This  picture strongly suggests that many almost
Mathieu phenomena should admit universal counterparts in the analytic
category \cite{avila0,avila1,avila2}.

In the subcritical regime this program has largely been completed.
The solution of the Almost Reducibility Conjecture (ARC)
\cite{avila1,avila2} shows that subcriticality implies almost
reducibility.  As a consequence, a broad range of spectral features
becomes universal, including absolutely continuous spectrum for all
phases and frequencies, 
and the absolute
continuity and sharp $1/2$ H\"older continuity of the
IDS \cite{avila1,avila2,aj1}.   Thus, qualitative
universality in the subcritical regime is now well understood.

The situation in the supercritical regime is markedly different.
Certain results are known to be universal.  For example,  the metal–insulator transition for the
almost Mathieu operator \cite{J} was extended (
measure-theoretically) in the seminal work of
Bourgain and Goldstein \cite{bg}, whose semi-algebraic method proved
robust and generated substantial subsequent developments
\cite{B1,jls}.  However, sharp arithmetic phenomena in the
supercritical regime—both of the “all-frequency” type and of the
sharp transition type—have proved far more delicate.

This is reflected already in the behavior of the integrated density of
states. In the subcritical regime, all-frequency absolute continuity
of the IDS is now known to be universal. In the supercritical regime,
by contrast, the corresponding sharp arithmetic behavior was
previously known {\bf only} for the almost Mathieu operator, combining the
Diophantine case
 \cite{J} and the Liouville case
\cite{ad}. For any other potential, this remained 
open.

Also, while the Lyapunov exponent for the almost Mathieu
operator is analytic (indeed constant) on the spectrum, for general
analytic potentials analyticity holds only on spectral components with
constant acceleration.  This observation led to the conjecture
by J. You \cite{You}, that certain sharp
regularity properties of the IDS should be
universal only within fixed-acceleration classes. In this sense, acceleration naturally suggests the relevant universality divide in the supercritical regime. 

More broadly, the universality of several prominent arithmetic
phenomena for almost Mathieu operators in the supercritical regime
remains a central open problem.  Existing approaches, including the
semi-algebraic set method, typically require non-arithmetic
restrictions on the frequency.  This sharp contrast between the
subcritical and supercritical regimes underscores that it is precisely
the \emph{sharp arithmetic} aspects of the theory that resist
extension beyond the cosine case.

At the same time, the difficulty is not merely technical. The known
proofs of sharp almost Mathieu arithmetic results, though very
different in method, on both the supercritical and subcritical sides,
rely on special structural features of the model, most notably
reflection symmetry. The next subsection explains how this enters on
the supercritical side in the direct localization approach followed by subsections that explain the difficulties of the duality approach.

\subsubsection{ Zero-counting and direct localization methods.}\label{zero}
The direct localization method underlying the solution of the
almost Mathieu AAJ in \cite{JLiu} , stems from the one originally
developed in
\cite{Lana94}. It has been extended throughout the positive
Lyapunov exponent regime in \cite{J}, then
enhanced to treat exponential frequency resonances in \cite{aj}, and modified in
several important technical features in \cite{ly}. Its main mechanism is a sharp control of exponential resonances through zero-counting for finite-volume determinants.

Resonances are places where box
restrictions have  exponentially close eigenvalues compared to the
distances between the boxes. For quasiperiodic operators, one kind is
so-called {\it exponential frequency resonances}: if $\dist(q\alpha,\Z)<e^{-cq}$
for infinitely many $q, $ a condition holding for an explicit dense $G_\delta$
but measure zero set of $\alpha.$ For {\it even} potentials there are also
reflection-based {\it exponential phase resonances}, where
$\dist(\theta+n\alpha,-\theta)<e^{-cn}$ for infinitely many $n,$  a
condition holding for an explicit 
topologically generic but measure zero set
of $\theta.$

 A sharp way to
treat exponential frequency resonances was developed in \cite{JLiu}, finally
solving the original AAJ, the phase part was handled in \cite{JLiu1}.

The essence of all the above proofs is in showing that for the almost Mathieu
operator there are no other types of resonances, thus when one removes
the arithmetically explicit measure zero set of $\theta$ for which
there are infinitely many phase resonances, only the frequency
resonances need to be dealt with.

Another important case where there are only frequency resonances, and that can therefore be effectively treated by the direct methods is that
of $H_{v,\alpha,\theta}$ with $v$ that is monotone on the period. 
 In \cite{jk3}, which appeared almost simultaneously with the first
version of this paper, the universality of AAJ is
established for all anti-Lipschitz monotone potentials, developing the
ideas of \cite{JLiu} in the framework of the arguments of \cite{jk, Kach}. This was recently extended significantly beyond the anti-Lipshitz setting in \cite{jls1}.

The key to all the methods that go back to \cite{Lana94} is that the set of
phases where eigenvalues of a restriction to a box of size $q$ can be
  exponentially close to $E,$ is confined to $q$ exponentially small
  intervals around zeros of the determinants of box restrictions. It is a feature that
  enables localization proofs not only in the almost
  Mathieu results \cite{J,aj,JLiu,
  JLiu1,liuresonant} where it is due to the
  fact that, $\cos$ being an {\it even} function, the determinants of   restrictions
  to boxes of size $q$ are polynomials of degree $q$ in the shifted
  $\cos$, but also in many other recent localization results that
fundamentally go back to the same zero-counting idea,
e.g. \cite{jks,bj, hdry, jy,wxyzz,JHY,jz,jk,Kach,r,omar,zhu,kz,jk3}.

In particular, the same phenomenon has been long understood to also happen for $H_{v,\alpha,\theta},$ with
 {\it even} analytic $v$ at energies with \emph{acceleration 1}. Indeed, Avila's {\it proof} of the global theory \cite{avila0}
essentially showed that, for energies with acceleration $1$, traces of transfer-matrices 
(i.e., determinants of block-restrictions with
periodic boundary conditions) of size $q_n$ effectively behave like
trigonometric polynomials of
degree $q_n,$ so also have  no more than $2q_n$
  zeros. Thus the extension of techniques of \cite{J} and even those
 of  \cite{JLiu} to the acceleration-one, even-potential setting was not expected to encounter any fundamental difficulties.
It has now been implemented in \cite{hs,lpos}, who obtained sharp
estimates on zero count through an approach different from Avila's.

If $v$ is not even, however, the situation changes qualitatively.  Even in the acceleration-one case, the effective zero count doubles to $2q_n$,
breaking the confinement mechanism that underlies the resonance
analysis of \cite{Lana94} and its successors.  For essentially the same
reason, these approaches do not currently extend beyond the
acceleration 1 case, where the effective degree again
increases. A related difficulty is visible already in perturbative regimes, where the
methods are completely different.  It is significantly more
difficult to obtain the result without requiring $v$ to be even
\cite{sinai, eli97, tv,csz2} than for even $v$ 
\cite{fsw,gyzh,csz1}. \footnote{It is claimed in \cite{tv} that evenness
    of $v$ is also de-facto required in \cite{sinai}.}

The direct approach to robust AAJ therefore requires significant new
ideas to extend beyond the {\it even} acceleration $1$ setting. Here, we
instead pursue the dual reducibility route of \cite{ayz}. However, as
we explain next, it
encounters the same symmetry barrier, and then additional ones.

\subsubsection{The classical $SL(2,\R)$ dual framework}\label{subsec:intro-sl2r}

Much of the sharp arithmetic theory of the almost Mathieu operator is
ultimately built upon a dual dynamical framework that may be
summarized by the label $SL(2,\R)$.  Three features are
simultaneously present.

{\bf The ``SL''.}  Both $H$ and its dual are finite-range operators,
thus  producing  finite-dimensional cocycles.

{\bf The 2.} The relevant dynamics are 
two-dimensional. The dual operator is also second-difference, leading
to two-dimensional dual dynamics.

{\bf The R.} The reflection symmetry makes the dual cocycle real, so
the dynamics lie in $SL(2,\R)$ rather than $SL(2,\C).$ The real structure ensures the
existence of a classical rotation number and enables the
rotation number–IDS correspondence that links dynamical invariants to
spectral quantities.

These three ingredients: finite-dimensional cocycle,
two-dimensionality, and real dynamics, form the basis of all reducibility-based
sharp arithmetic arguments for the almost Mathieu operator.

\subsubsection{The $\SL$: Absence of a finite-dimensional cocycle beyond the trigonometric case}\label{subsec:intro-nococycle}

The first part of the classical $SL(2,\R)$ framework already fails
outside
the trigonometric setting.  

For trigonometric potentials, the dual
operator is finite range, and the associated eigenvalue equation can be
encoded by a finite-dimensional cocycle. For general analytic potentials
the dual has infinite range, so this transfer-matrix formulation is unavailable, and the usual
rotation-number and reducibility methods  cannot be applied  on the
dual side.  
To use these dynamical methods within the duality approach
to general analytic potentials, one needs to somehow uncover
finite-dimensional dynamics in the infinite-range setting.

\subsubsection{The 2}\label{subsec:intro-typeI}

Many of the tools of spectral analysis of operators \eqref{sch},
especially for the
arithmetically delicate parts,  have been developed for
second-difference operators, with two-dimensional matrix cocycles. From the Wronskian-based arguments, to
celebrated Kotani theory, to power-law subordinacy, all require
second-difference for sharp formulations, and, some remarkable recent
work (e.g. \cite{wxz})  notwithstanding, have been resisting extensions
allowing applications to higher-order problems.

At the same time,
almost Mathieu is the \emph{only} operator \eqref{sch} whose dual
dynamics is two-dimensional. Moreover, by \cite{gjyz}, for trigonometric 
potentials the dual symplectic cocycle is partially hyperbolic with center of dimension 
$2k$ where 
$k$ is the T-acceleration.  Thus the classical two-dimensional framework
has any chance of surviving only in the T-acceleration-one,  i.e. Type I, regime; beyond that, the relevant dual dynamics are higher-dimensional.

\subsubsection{The $\R$: symmetry (evenness) barrier and its dual
  manifestation}\label{subsec:intro-symmetry}

An approach to AAJ, developed in \cite{ayz} and related works, proceeds
through dual reducibility.  Localization in the supercritical regime is
obtained by establishing reducibility (or almost reducibility) of the
dual cocycle in the subcritical regime and transporting this
information via Aubry duality.  A crucial  input here is
that, for the almost Mathieu operator, the dual cocycle is real.

The real $SL(2,\R)$ structure provides a classical
rotation number and, through the Johnson--Moser correspondence,
 a direct link between dynamical quantities and the integrated
density of states.  It is this rotation number–IDS relation that allows
one to convert reducibility information into sharp spectral
conclusions.

In particular, the duality-based proof of almost Mathieu localization as well as many
other almost Mathieu proofs are
hinged on the fact, going back to \cite{Puig}, that if $E$ is an eigenvalue
of $H_{cos,\alpha,\theta}$
then \begin{equation}\label{rhoe}\rho(E)=\pm\theta+k\alpha (\mod
  \Z).\end{equation}

Once reflection symmetry is removed, however, the dual dynamics are no
longer governed by real cocycles.
 Instead,
one encounters non-real Hermitian-symplectic dynamics, and the
classical $\SL(2,\R)$ rotation number is no longer directly available. As a result, the reducibility-to-localization pathway based on \eqref{rhoe} and the classical Johnson--Moser rotation-number--IDS correspondence break down at a structural level.

\medskip

This obstruction is not limited to localization.
In the almost Mathieu setting, both the all-frequency absolute
continuity of the IDS and the sharp $1/2$-Hölder regularity in the
supercritical regime are derived through the real $SL(2,\R)$ structure
of the dual cocycle and the classical rotation number–IDS
correspondence.  

Thus all known proofs of sharp arithmetic phenomena for the almost
Mathieu operator use reflection symmetry. It enters the
zero-counting arguments for localization, the real
$SL(2,\mathbb R)$ dynamics used in dual reducibility, and
the rotation number--IDS correspondence used to prove
absolute continuity and sharp regularity of the IDS.
Its absence creates a significant obstruction for all three universality
problems.

\subsubsection{Type I as the correct universality class and the dual center}
\label{subsec:intro-typeI}

Our work builds on the duality approach to global theory initiated in
\cite{gjyz}, 
which has provided powerful tools to various spectral
problems.

A key result of \cite{gjyz} is that for trigonometric polynomial $v$
of degree $d$
the dual $Sp(2d,C)$ cocycles are partially hyperbolic with center of
dimension $2k$ where $k$ is the T-acceleration. 
This serves as a foundation  to the
current work, as well as other subsequent developments.

In particular, when $k=1$ the center of this high-dimensional dynamical system remains 
two-dimensional, allowing to potentially recover some of the $2$ in
$SL(2,R)$. 
 This lead to the introduction of Type I, in \cite{gjy},
the notion which, in the supercritical case, coincides with acceleration $1.$
Before this concept was introduced, Type I was already used as a
simplifying property in
\cite{gjz} where the absolute continuity of the IDS was proved for
Diophantine $\alpha.$
Restricting to this robust (open and conjecturally dense, see Appendix
\ref{app})
class, identifies the natural setting in which one may hope to recover
the two-dimensional part of the dual picture for
the reducibility methods.
The same restriction also appears naturally from the direct side: as
discussed in \ref{zero}, acceleration one is also necessary for the
zero-counting method,
though there it must still
be coupled with symmetry.
Thus Type I is the natural common regime in which both direct and dual approaches retain part of the classical structure. The earlier preprint \cite{gjy} also developed the Kotani--Puig
approach to Cantor spectrum for even trigonometric Type~I potentials.

For trigonometric potentials, Type I gives a two-dimensional center, but it does not
resolve the difficulties caused by the absence of reflection symmetry or the infinite
range of a general analytic dual. On the direct side, acceleration one controls the zero
count only when combined with symmetry. Extending that localization method to
non-even potentials remains open.

For the dual approach, while Type~I restores the two-dimensional center dynamics for trigonometric $v$, the other two key
ingredients of the $SL(2,\R)$ framework remain absent.

These are the difficulties
we address.
 We construct intrinsic symplectic center dynamics that persist in the infinite-range setting. In dimension two,
we introduce a projectively real concept that replaces the missing real structure, and establish a generalized rotation–IDS correspondence for associated Hermitian-symplectic dynamics, leading to the universality
results of Section~\ref{subsec:intro-results}.
The extension of the Ten Martini theorem to general analytic
Type~I potentials in \cite{gjy1} uses the intrinsic center
construction developed here. This construction also enters the
supercritical dry Ten Martini theorem in the type I regime \cite{lxz} and, in higher
dimensions, the general supercritical Ten Martini theorem
\cite{gjy2}. We expect further applications.

\subsection{Setup and definitions}

Let $v \in C^\omega(\T,\R)$ and $\alpha \in \R \backslash \Q$.
Consider the analytic one-frequency Schrödinger operator
\begin{equation}\label{eq:sch}
(H_{v,\alpha,\theta} u)_n
= u_{n+1} + u_{n-1} + v(\theta + n\alpha) u_n,
\qquad n \in \Z.
\end{equation}

Let $L(E)$ denote the Lyapunov exponent of the associated Schrödinger cocycle
$(\alpha, A_E)$, $\omega(E)$ denote its acceleration, and
$\bar{\omega}(E)$ its T-acceleration. For supercritical energies, we have $\bar{\omega}(E)=\omega(E).$

\begin{definition}
An energy $E$ is called \emph{Type~I}
if its T-acceleration satisfies $\bar{\omega}(E)=1$.
\end{definition}

Let $I\subset \R$ be the set of all supercritical Type~I energies.

We define the arithmetic exponent of $\alpha$ by
\begin{equation}\label{beta}
\beta(\alpha)
= \limsup_{n\to\infty}
\frac{-\log \| n\alpha \|_{\R/\Z}}{|n|}.
\end{equation}

We say that $\alpha$ is \emph{Diophantine} if $\beta(\alpha)=0$.

Let $N(E)$ denote the integrated density of states (IDS) associated with
\eqref{eq:sch}.

An energy is called \emph{non-critical} if it lies outside the critical regime
in the sense of Avila's global theory.

Precise definitions and further background are given in Section~\ref{pre}.
\subsection{Universality Results}

\begin{Theorem}
\label{t1}
Let $v \in C^\omega(\T,\R)$ and $\alpha \in \R \backslash \Q$. Then the
sharp arithmetic transition in frequency (AAJ) is
  universal for type I energies $E$ of  operators
  $H_{v,\alpha,\theta}.$ That is
\begin{itemize}
\item For almost every $\theta$,
the operator $H_{v,\alpha,\theta}$ has pure point spectrum on $I\cap\{E: L(E)>\beta(\alpha)\}$;
\item For every $\theta$, $H_{v,\alpha,\theta}$ has purely singular
  continuous spectrum on $I\cap\{E: L(E)<\beta(\alpha)\}$.
\end{itemize}
\end{Theorem}

In particular, this applies to all type I operators thus 
to {\it all} existing
models with previously known (not necessarily sharp) arithmetic localization results, and
their analytic  neighborhoods.
\begin{Remark} 
The singular continuous part of the AAJ was established in \cite{ayz} for
all Lipschitz $v$ and all $E$, 
what we
prove is the localization statement.
\end{Remark}

For the particular case of the    almost Mathieu neighborhood, the AAJ
conjecture immediately leads to a corollary of a particularly nice form. Let
  $H^\delta_{\lambda,\alpha,x}$ be given by
\begin{equation}\label{hdelta}
 (H^\delta_{\lambda,\alpha,\theta}u)_n=u_{n+1}+u_{n-1}+(2\lambda\cos2\pi(\theta+n\alpha)+\delta f(\theta+n\alpha))u_n,\ \ n\in\Z.
\end{equation}
 \begin{Corollary}For
   $\alpha\in\R\backslash\Q$ and any 1-periodic real analytic $f\in C^\omega_h(\T,\R)$, there exists
   $\delta_0(\lambda,\alpha,\|f\|_h)$ such that if   $|\delta|<\delta_0$, we have
\begin{enumerate}
\item If $|\lambda|<1$, $H^\delta_{\lambda,\alpha,\theta}$ has purely absolutely continuous spectrum for all $\theta$;
\item If $1<|\lambda|<e^\beta$, $H^\delta_{\lambda,\alpha,\theta}$ has purely singular continuous spectrum for all $\theta$;
\item If $|\lambda|>e^\beta$, $H^\delta_{\lambda,\alpha,\theta}$ has Anderson localization for a.e. $\theta$.
\end{enumerate}
 \end{Corollary}
At fixed $\alpha$, subcritical openness and \cite{avila1,avila2} imply (1) above.
For $|\lambda|>1$, Type~I openness \cite{gjy1} and continuity of the spectrum and $L(E)$ preserve
the strict inequalities in (2) and (3), since $L(E)=\ln|\lambda|$ on the unperturbed spectrum.

Finally, using that by the global theory \cite{avila0}  for typical
operators $H_{v,\alpha,\theta}$  there are no critical energies, and
invoking the almost reducibility theorem \cite{avila1,avila2}, we obtain
\begin{Corollary}
For $\alpha\in\R\backslash\Q$ and a (measure-theoretically) typical Type I operator  $H_{v,\alpha,\theta}$, we have
\begin{enumerate}
\item $H_{v,\alpha,\theta}$ has purely absolutely continuous spectrum for all $\theta$ on $\{E;L(E)=0\}$;
\item $H_{v,\alpha,\theta}$ has purely singular continuous spectrum for all $\theta$ on $\{E:0<L(E)<\beta(\alpha)\}$;
\item $H_{v,\alpha,\theta}$ has Anderson localization for a.e. $\theta$ on $\{E:L(E)>\beta(\alpha)\}$.
\end{enumerate}
\end{Corollary}
Here ``measure-theoretically typical'' means prevalent: fixing some probability measure $\mu$ of compact support (describing a set of admissible perturbations $w$), a property is measure-theoretically typical if it is satisfied for almost every perturbation $v+w$ of every starting condition $v$.

We now move to our second main result. The  integrated density of states (IDS) is defined for Schr\"odinger operators $(H_{v,\alpha,\theta})_{\theta\in\T}$ by
$$
N(E)=\int_{\T}\mu_{\theta}(-\infty,E]d\theta,
$$
where $\mu_\theta$ is the spectral measure associated with
$H_{v,\alpha,\theta}$ and $\delta_0\in\ell^2(\Z).$ We have

\begin{Theorem}
\label{t2}
The absolute continuity of the IDS for all $\alpha$ is universal for
  all non-critical type I operators $H_{v,\alpha,\theta}$.
\end{Theorem}
\begin{Remark} \label{rem} The novelty here lies in the supercritical
  regime. The non-critical condition is
  essential \cite{last,ak}, and is conjectured to be necessary.  The type I condition,
  however, is likely not necessary even in the supercritical regime,
  but removing it would require some
  further ideas. 
\end{Remark}
\begin{Remark} The almost Mathieu proofs \cite{J,aj1} and \cite{ad}
heavily use the specifics of the $\cos.$ The Liouville argument of
\cite{ad} besides using several specific almost Mathieu facts,
transpires entirely in the subcritical range and is thus superseded by
\cite{avila1} and not useful for the supercriticality. As for  the
localization pathway, while the latter is shown to be universal for type I in Theorem
\ref{t1},
we actually use Theorem \ref{t2} to prove Theorem \ref{t1}, so this cannot be used
either.
\end{Remark}
Previous supercritical results on absolute continuity either
require Diophantine $\alpha$ \cite{gjz,xwyz} or are not arithmetic at all \cite{gs2},
needing highly implicit elimination of $\alpha$ due to the need to get
rid of the so-called ``double resonances''. We note also that the
method of \cite{gjz} has no hope to be extendable to Liouvillean
frequencies, since it is based on the homogeneity of the spectrum,
a property that can fail for Liouvillean frequencies
\cite{alsz}. The methods of \cite{xwyz, gwx} not only require the
Diophantine condition but are perturbative (work only for large couplings
with the largeness dependent on $\alpha$ \footnote{They extend however
 in some other ways.}).

As above, we also have some immediate corollaries for the neighborhood
of the almost Mathieu operator

\begin{Corollary}
For $\alpha\in\R\backslash\Q$, $|\lambda|\neq 1$ and any 1-periodic real analytic $f\in C^\omega_h(\T,\R)$, there
exists $\delta_0(\lambda,\alpha, \|f\|_h)$ such that if   $|\delta|<\delta_0$,
the integrated density of states of
$(H^\delta_{\lambda,\alpha,\theta})_{\theta\in\T}$ is absolutely
continuous.
\end{Corollary}
and for typical type I operators
\begin{Corollary}
For $\alpha\in\R\backslash\Q$ and a (measure-theoretically) typical type I operator
$H_{v,\alpha,\theta}$, the integrated density of states is absolutely
continuous.
\end{Corollary}

Our next universality result concerns the sharp H\"older exponent. The IDS of all non-critical almost Mathieu operators with
$\beta(\alpha)=0$ are {\bf exactly}
$1/2$-H\"older continuous \cite{aj1}. This statement is sharp and
optimal already in the almost Mathieu family (there
are square root singularities at gap edges \cite{Puig1}, and the result does not
hold for non-Diophantine $\alpha$ \cite{alsz} or for some Diophantine $\alpha$
at criticality (\cite{B1}, Remark after Corollary
8.6). 
In fact, a lot more delicate statement
was recently obtained about local Holder continuity \cite{lyz}, but
the overall exponent $1/2$ is sharp. This sharp $1/2$-H\"older
regularity for Diophantine $\alpha$ is universal throughout the
subcritical regime (through a combination of \cite{aj1} and
\cite{avila2}). For $\beta(\alpha)=0$, see \cite{lyzz}. In the supercritical regime, You's conjecture ties the
modulus of continuity to the acceleration, in particular, expecting $1/2$-H\"older
regularity {\it only} for Type I energies. Here we prove
it for all non-critical type I
operators, thus resolving the corresponding part of You's conjecture.

\begin{Theorem}\label{t3}
The  $1/2$-H\"older
regularity of the IDS for all $\alpha$ with $\beta(\alpha)=0,$  is universal  for
  all non-critical type I operators $H_{v,\alpha,\theta}$.

\end{Theorem}

Other than the mentioned subcritical universality result of \cite{aj1,avila2}, the exact
uniform $\frac{1}{2}$-H\"older continuity of the  IDS for Diophantine $\alpha$
was earlier obtained in
\cite{amor} 
in the perturbatively small
regime of Eliasson \cite{eli92}. More
recently it was extended to smooth perturbative almost reducibility regime
in \cite{ccyz}. Non-sharp 
$\frac{1}{2}$-H\"older continuity (that is $1/2-\e$ for all $\e>0$)
was obtained under various further conditions in \cite{b1,gs2,
  gyzh1,hs} (see also \cite{gs1} for first H\"older estimates). There have been no previous sharp $\frac{1}{2}$-H\"older results
in the supercritical regime other than for the almost Mathieu
\cite{aj1}.

Finally, as above we have the following immediate corollaries

\begin{Corollary}
For $|\lambda|\neq 1$ and any 1-periodic real analytic $f\in C^\omega_h(\T,\R)$, there
exists $\delta_0(\lambda,\alpha,\|f\|_h)$ such that if   $|\delta|<\delta_0$,
the integrated density of states of
$(H^\delta_{\lambda,\alpha,\theta})_{\theta\in\T}$ with $\beta(\alpha)=0$ is $\frac{1}{2}$-H\"older continuous.
\end{Corollary}
\begin{Corollary}
For a (measure-theoretically) typical type I
operator  $H_{v,\alpha,\theta}$ with $\beta(\alpha)=0$, the integrated density of states is $\frac{1}{2}$-H\"older continuous.
\end{Corollary}

\subsection{Structural results}
The Aubry dual of operator \eqref{sch} is the  
following quasiperiodic long-range operator (see Subsection \ref{au}):
\begin{equation}\label{long-range}
(\widehat{H}_{v,\alpha,\theta}u)_n=\sum\limits_{k\in\Z}\widehat{v}_k u_{n+k}+2\cos2\pi(\theta+n\alpha)u_n, \ \ n\in\Z.
\end{equation} Let $T: \C^\Z\to \C^\Z$ be the shift,  $(Tu)_n=u_{n+1}.$
Our first structural result shows that the symplectic center
dynamics of finite-range duals persist for general analytic
potentials. Although the  dual equation then has no
finite-dimensional transfer-matrix formulation, it still has
a distinguished finite-dimensional family of formal solutions
on which the shift induces analytic Hermitian-symplectic dynamics.
We construct this family and its symplectic structure by
trigonometric approximation, without loss of analyticity.
The resulting center is intrinsic to the dual eigenvalue
equation and makes finite-dimensional symplectic methods
available for infinite-range dual operators.
\begin{Theorem}[Intrinsic symplectic structure]\label{hss}
 Assume $v\in C^{\omega}(\T,\R)$, and $\bar{\omega}(E_0)=k$. There exist a neighborhood $I$ of $E_0$, $M\in C^\omega(I\times \T,Sp(2k,\C))$  and linearly
  independent $O^j\in C^\omega(I\times \T,\C^\Z), \;j=1,\ldots,2k$,  such
  that each $O^j(E,\theta)$ is a formal solution of the dual eigenvalue
  equation: $\widehat{H}_{v,\alpha,\theta}O^j(E,\theta)=EO^j(E,\theta),$ and
  $$(TO^1(E,\theta),\ldots,TO^{2k}(E,\theta))=(O^1(E,\theta+\alpha),\ldots,O^{2k}(E,\theta+\alpha))M(E,\theta).$$ 
  Moreover, if $L(E_0)>0$, then for any $0<h<L(E_0)/2\pi$, there is $E_0\in I_h\subset I$, such that $O^j\in C^\omega_h(I_h\times \T,\C^\Z),$
$$
L_j(M(E,\cdot+i\e))=L_j(M(E,\cdot)),\ \ j=1,\ldots,2k, \ \ \forall |\e|<h.
$$
If, in addition, $E\in I_h\cap \Sigma_{v,\alpha}$, $L(E)>0$ and $\omega(E)=k$, the
center exponents vanish:
\[
L_j(M(E,\cdot+i\e))=0,\ \  j=1,\ldots,2k, \ \ \forall |\e|<h.
\]

\end{Theorem}

\begin{Remark}\label{rem:center-frame}
Let $O(\theta)=(O^1(\theta),\ldots,O^{2k}(\theta))$ and
$\mathcal E(\theta)=\operatorname{span}\{O^1(\theta),\ldots,O^{2k}(\theta)\}$.
For a fixed solution bundle $\mathcal E$  another
analytic symplectic frame has the form
$\widetilde O(\theta)=O(\theta)B(\theta)$, with
$B(\theta)\in Sp(2k,\C)$. Its cocycle matrix is
\[
\widetilde M(\theta)=B(\theta+\alpha)^{-1}M(\theta)B(\theta).
\]
While $M$ is a matrix representative in a chosen frame, the shift action itself is independent of its coordinates.
Our statements of
convergence below and the rotation quantities require a choice of frame and lifts, as
in Remark~\ref{rem:rotation-normalization}. However, in the
Type~I setting, the increments of the rotation difference
are frame independent. Indeed, for two analytic symplectic frames
on an interval $I$, and continuous real lifts,
\[
\widetilde\rho_1(E)-\widetilde\rho_2(E)
=\rho_1(E)-\rho_2(E)+\kappa,
\qquad E\in I,
\]
where $\kappa\in\mathbb Z+\alpha\mathbb Z$ is independent of $E$.
This follows from the conjugacy rule for the projective rotation
number, since the degree of the frame change is constant on $I$.
Thus
\[
\bigl(\rho_1(E_1)-\rho_2(E_1)\bigr)
-\bigl(\rho_1(E_2)-\rho_2(E_2)\bigr)
\]
is independent of the frame and continuous lifts, making the
difference formulation of the rotation--IDS correspondence intrinsic.
\end{Remark}

For Type I energies, the two-dimensional center dynamics have an additional important feature.
We formulate this through the notion of \emph{projectively real cocycles}: after factoring
out a scalar phase, the projective action is conjugate to that of a real
$\SL(2,\mathbb R)$ cocycle. This makes it possible to define rotation data beyond the
classical symmetric setting.

\begin{Theorem}[Projectively real structure]\label{thm:projectively-real}
Let $E_0$ be a Type I energy and let $M$ be the analytic $Sp(2,\mathbb C)$ cocycle
associated with the intrinsic two-dimensional center dynamics and $h$, $I_h$ be as in Theorem \ref{hss}. Then $M$ is
projectively real. More precisely, there exist 
$
\phi\in C^\omega(I\times \mathbb T,\mathbb R),
\;
{\bf C}\in C^\omega(I\times \mathbb T,\SL(2,\mathbb R)),
$
such that
\[
M(E,\theta)=e^{2\pi i \phi(E,\theta)}\,{\bf C}(E,\theta).
\]
Moreover, if $L(E_0)>0$, the cocycle $(\alpha,{\bf C}(E,\cdot))$ is regular on the strip
$\{|\Im \theta|<h\}$ for any $E\in I_h$.
\end{Theorem}

The decomposition in Theorem~\ref{thm:projectively-real} allows one to associate to
$M$ two rotation quantities: one coming from the scalar phase and one from the
underlying real cocycle.

\begin{definition}[Rotation pair]\label{def:rotation-pair}
Let $(\alpha,A)$ be a projectively real analytic cocycle, so that
\[
A(\theta)=e^{2\pi i \phi(\theta)}\,C(\theta),
\]
where $\phi\in C^\omega(\mathbb T,\mathbb R)$ and
$C\in C^\omega(\mathbb T,\SL(2,\mathbb R))$ is homotopic to the identity.

Set
\[
\widehat\phi:=\int_{\mathbb T}\phi(\theta)\,d\theta,
\]
and let $\rho(\alpha,C)\in \mathbb R/\mathbb Z$ denote the fibered rotation number of
the real cocycle $(\alpha,C)$.

The \emph{rotation pair} of $(\alpha,A)$ is the ordered pair
$
(\rho_1(\alpha,A),\rho_2(\alpha,A))
$
defined by
\[
\rho_1(\alpha,A):=\rho(\alpha,C)+\widehat\phi,
\qquad
\rho_2(\alpha,A):=-\rho(\alpha,C)+\widehat\phi.
\]
Equivalently,
\[
\rho_1(\alpha,A)+\rho_2(\alpha,A)=2\widehat\phi,
\qquad
\rho_2(\alpha,A)-\rho_1(\alpha,A)=-2\rho(\alpha,C)
\quad \text{in } \mathbb R/\mathbb Z.
\]
\end{definition}

Rotation-number formulas for the integrated density of
states of block Jacobi operators were developed by
Schulz-Baldes \cite{schulzbaldes}.
Our final main theorem identifies the rotation pair for the cocycle $M$ from
Theorem~\ref{thm:projectively-real}, with the integrated density of states.

\begin{Theorem}[Rotation pair--IDS correspondence]\label{thm:rotation-ids}
Let $E_0$ be a supercritical Type I energy in the spectrum of $H_{v,\alpha,\theta}$, with
$\alpha\in\mathbb R\backslash\mathbb Q$ and $v\in C^\omega(\mathbb T,\mathbb R)$.
Then there exists a neighborhood $I$ of $E_0$ such that for every $E_1,E_2\in I$
\[
N(E_1)-N(E_2)=\left(\rho_1(E_1)-\rho_2(E_1)\right)-\left(\rho_1(E_2)-\rho_2(E_2)\right).
\]
\end{Theorem}

\begin{Remark}\label{rem:even-case-rotation-ids}
In the even case, the projectively real decomposition reduces to the classical
$\SL(2,\mathbb R)$ setting, and Theorem~\ref{thm:rotation-ids} specializes to the usual
formula
\[
N(E_1)-N(E_2)=2(\widehat\rho(E_1)-\widehat\rho(E_2)).
\]
\end{Remark}
\begin{Remark}
The classical rotation--IDS correspondence
\cite{johonson and moser,as83} was extended to block Jacobi
operators by Schulz--Baldes \cite{schulzbaldes}. In \cite{lw}
Li and Wu gave a fibered-rotation formulation for Hermitian-symplectic cocycles, using it to also establish the IDS correspondence. For transfer cocycles of block Jacobi operators, these rotation numbers agree after
matching orientation and normalization conventions,
up to the choice of lift. 
For projectively real cocycles, a direct computation shows that the quantity $\rho_2(E)-\rho_1(E)$ agrees with the corresponding generalized fibered rotation number of
\[
M(\theta)=e^{2\pi i\phi(\theta)}C(\theta),
\]
up to the normalization/sign convention. We thank Xianzhe
Li for this observation.

 For our purposes, however, the essential point is the full rotation pair $(\rho_1,\rho_2)$, not only the scalar quantity $\rho_2-\rho_1$. Indeed, the pair canonically separates the scalar mean phase from the real matrix dynamics.
Theorem~\ref{thm:projectively-real} gives the underlying hidden $\SL(2,\R)$ cocycle, whose
 rotation number determines the IDS locally up to an additive
constant. For the reducibility-to-localization argument we need the full pair: reducibility
produces eigenfunctions at phases $\rho_1(E)$ and $\rho_2(E)$.

\end{Remark}

\section{Main ideas}
Our starting point is the duality approach to Avila's global theory developed in
\cite{gjyz}.  For trigonometric polynomial potentials, the dual eigenvalue equation defines a
finite-dimensional cocycle; in the Type I regime it is partially hyperbolic with a two-
dimensional center \cite{gjyz,gjy}. We first work with these finite-range duals, then pass
to general analytic potentials by proving convergence of the center dynamics under
trigonometric approximation. 

The main new ingredients that may also be of independent interest are the following.

\begin{enumerate}
\item \emph{Intrinsic 
symplectic center dynamics.}
Convergence of the dual Lyapunov exponents under trigonometric
polynomial approximation was proved in \cite{gjyz}. 
 Here we
prove convergence of the center dynamics themselves, in
suitable analytic symplectic frames. The construction works
for centers of any dimension, thus for any T-acceleration $k$,
without the Type~I restriction.

We construct center sections from Green's functions and use
Aubry duality to relate the corresponding resolvents.
Avila's global theory supplies uniform hyperbolicity at
suitable complex phase shifts, ensuring that these resolvents
exist. The resolvent identity then gives convergence estimates
under trigonometric approximation. From these estimates,
we construct convergent symplectic frames and center cocycles.
The limit defines intrinsic analytic $Sp(2k,\mathbb C)$
dynamics on a $2k$-dimensional family of formal solutions
of the infinite-range dual eigenvalue equation.

In dimension two, this construction is used in the robust
and supercritical dry ten martini theorems for Type~I
operators \cite{gjy1,lxz}. Its higher-dimensional form is
a key ingredient in the general supercritical ten martini
theorem \cite{gjy2}.

\item \emph{Projectively real structure and rotation data.}
For supercritical Type I energies, the two-dimensional center admits an additional rigidity: after
factoring out a scalar phase, its projective action is conjugate to
that of a  {\it subcritical} real
$\SL(2,\mathbb R)$ cocycle. This allows us to define a rotation pair and to establish the
corresponding rotation--IDS relation. Combined with almost reducibility in the
subcritical regime, this yields the absolute continuity and sharp H\"older regularity of
the IDS.  Our proof is based on a dynamical point of view of the $m$-function, going back to Johnson-Moser \cite{johonson and moser}. The subcritical real center dynamics are also used to prove the 
ten martini results of \cite{gjy1,lxz} and the sharp
arithmetic transitions in phase of \cite{flz}.

\item \emph{A new completeness argument.}
The  reducibility-to-localization argument was first
developed in \cite{ayz} exploiting certain quantitative information on the
 almost Mathieu reducibility. A simple general argument was then
presented in \cite{jk2}, and an arithmetic in $\theta$ way was found
 in \cite{gy}. All these proofs, as well as other related 
  developments. were crucially symmetry-based: the fact that each
  eigenvalue $E$ corresponds only to two phases $\pm \theta(E),$ For
  non-even Type I family, we discover here a different phenomenon:
  each eigenvalue $E$ corresponds to two phases, $\rho_1(E)$ and
  $\rho_2(E)$ with $\rho_1(E)\neq -\rho_2(E).$ Thus the previous
  localization arguments \cite{ayz,jk2,gy} do not work. 
 We develop a new completeness argument that works in this
asymmetric setting and leads to arithmetic localization.

\item \emph{Multiplicative Jensen formula for the duals.}
We prove a multiplicative Jensen formula for dual cocycles.  The
   argument that {\it subcriticality implies dual supercriticality},
   plays an important role in giving a duality-based  proof of the almost
   reducibility conjecture \cite{ge}.
  Here we prove that {\it 
     supercriticality implies dual subcriticality},  as a direct corollary of our
  multiplicative Jensen's formula for the duals. More importantly, we  give {\it precise characterization of the
   subcritical radius} for the dual cocycles. This plays a fundamental role in obtaining
  sharp phase transition for  Type I operators.
 \item {\it Analytic symplectic basis for the center.} The sharp phase transition requires
the projectively real structure on the whole subcritical strip, without loss of
analyticity. We therefore need an analytic symplectic basis of the center on that
 entire strip.
This is highly non-trivial. Earlier constructions  \cite{dwxz, lxz} do not
guarantee preservation of the full analytic strip.  In dimension two, Avila \cite{avila} found a  way of doing this by using the Corona Theorem leading to various sharp spectral applications. However, this does not work in higher dimensions. {\it Here we develop a technique that works in any dimension.} This involves a completely different idea: abstract matrix factorization theory \cite{LSp,gks}. Although the factorization theorem starts with continuous matrix functions,
we show that combining it with analyticity and the Schwarz reflection principle allows us to extend an analytic  symplectic basis of the center to the whole strip. This technique may also be useful for other spectral problems.
\end{enumerate}

\subsection{Structure of the rest of the paper}

Section~3 contains the preliminaries. In Section~4 we prove the multiplicative Jensen
formula for the duals. In particular,  at supercritical spectral energies the dual center is subcritical. In Sections~5 and~6 we construct the intrinsic symplectic center
dynamics, prove convergence under trigonometric approximation, and obtain the
projectively real structure. In Section~7 we establish the rotation--IDS correspondence
and deduce absolute continuity and sharp H\"older regularity of the IDS. In Section~8
we develop the new reducibility-to-localization argument and prove the universality of
the sharp arithmetic transition.

\section{Preliminaries}\label{pre}
For background on ergodic Schr\"odinger operators, see e.g.
\cite{pf,df}. We will use $Sp(2k,\C)$ to denote the set of all $2k\times 2k$ Hermitian symplectic matrices, in some documents, denoted by $HSp(2k,\C)$. If $A\in C^\omega(I\times \T,Sp(2k,\C))$, then its holomorphic extensions satisfy
$$
A(\bar{z},\bar{\theta})^*J_{2k}A(z,\theta)=J_{2k}.
$$
The Lyapunov exponent of the complexified  Schr\"odinger cocycle, associated
 with operator \eqref{sch}  is defined as
\begin{align}\label{multiergodicsch}
L_\e(E)=\lim\limits_{n\rightarrow\infty}\frac{1}{n}\int_\T\ln
  \|A_E(\theta+i\e+(n-1)\alpha)\cdots A_E(\theta+i\e)\|d\theta; \ \ L(E):=L_0(E),
\end{align}
where
\begin{equation}\label{Sch}
A_E (\theta)=\begin{pmatrix}E-v(\theta)&-1\\ 1&0\end{pmatrix}.
\end{equation}

Avila showed \cite{avila0} that $L_\e(E)$ (as a function of $\e$) is an even convex piecewise
affine function with slopes in $2\pi\mathbb Z$ for real $E$.  

In particular, with acceleration defined as 
\begin{Definition}[\cite{avila0}]\label{acc}
{\rm The {acceleration} is defined by
$$
\omega(E)=\lim\limits_{\e\rightarrow 0^+}\frac{L_\e(E)-L_0(E)}{2\pi\e}.
$$}
\end{Definition}
we have \begin{equation} \label{ac}
  \omega(E)\in\N\cup\{0\}.
\end{equation}

In \cite{gjy} the T-acceleration was introduced as the slope at the first
turning point of the complexified Lyapunov exponent.

\begin{Definition}[T-acceleration]\label{gege1}
{\rm Let $h\le\infty$ be the natural boundary of analyticity of $v\in   C^{\omega}(\T,\R).$ The {\it T-acceleration} is defined by
$$
\bar{\omega}(E)=\lim\limits_{\e\rightarrow \e_1^+}\frac{L_\e(E)-L_{\e_1}(E)}{2\pi(\e-\e_1)}
$$
where $0\le\e_1<h$ is the first turning point  of the piecewise affine
function $L_\e(E)$.  If there is no
turning point, we set $\bar{\omega}(E)=0$.}
\end{Definition}

\begin{Definition}[Type I]\label{typei}
  {\rm We say $E$ is a {\it Type I energy} for operator
    $H_{v,\alpha,\theta}$ if
  $\bar{\omega}(E)= 1$. We say $H_{v,\alpha,\theta}$ is a {\it Type I operator,}
  if every $E$ in the spectrum of
  $H_{v,\alpha,\theta}$ is Type I.}
\end{Definition}

The property of T-acceleration being
equal to $1$ is stable in each $C^{\omega}_{h'}$,  $\e_1<h'<{h},$ and the following models and their analytic perturbations have Type~I energies in the appropriate regimes  
\cite{gjy1}: 
the almost Mathieu operator, the supercritical GPS model, the supercritical generalized
Harper's model, and analytic cosine type quasiperiodic operators.\footnote{Where
$v$ is a real analytic function satisfying the cosine type condition
introduced in \cite{sinai}, at nonperturbatively high coupling.} We
refer to \cite{gjy1} for more details.
\subsection{Cocycles and the Lyapunov exponents}\label{3.2}
Let $GL(m,\C)$ be the set of all $m\times m$ invertible matrices. Given $\alpha\in \R\backslash\Q$ and $A \in C^\omega(\T,GL(m,\C))$, we define the complex one-frequency cocycle $(\alpha,A)$ by:
$$
(\alpha,A)\colon \left\{
\begin{array}{rcl}
\T \times \C^{m} &\to& \T \times \C^{m}\\[1mm]
(x,v) &\mapsto& (x+\alpha,A(x)\cdot v)
\end{array}
\right.  .
$$
The iterates of $(\alpha,A)$ are of the form $(\alpha,A)^n=(n\alpha,A_n)$, where
$$
A_n(x):=
\left\{\begin{array}{l l}
A(x+(n-1)\alpha) \cdots A(x+\alpha) A(x),  & n\geq 0\\[1mm]
A^{-1}(x+n\alpha) A^{-1}(x+(n+1)\alpha) \cdots A^{-1}(x-\alpha), & n <0
\end{array}\right.    .
$$
Let $L_1(A)\geq L_2(A)\geq...\geq L_{m}(A)$ be the Lyapunov exponents
of $(\alpha,A)$ listed according to their multiplicities, i.e.,
$$
L_k(A)=\lim\limits_{n\rightarrow\infty}\frac{1}{n}\int_{\T}\ln\sigma_k(A_n(x))dx,
$$
where  $\sigma_1(A_n)\geq...\geq \sigma_m(A_n)$ denote its singular values (eigenvalues of $\sqrt{A_n^*A_n}$). Since the k-th exterior product $\Lambda^k A_n$ satisfies $\sigma_1(\Lambda^k A_n)=\|\Lambda^k A_n\|$, $L^k(A):=\sum_{j=1}^kL_j(A)$ satisfies
$$
L^k(A)=\lim\limits_{n\rightarrow \infty}\frac{1}{n}\int_{\T}\ln\|\Lambda^kA_n(x)\|dx.
$$

In particular, consider a finite-range quasiperiodic operator $\widehat{H}_{v,\alpha,\theta}$ given by \eqref{long-range} with $v(\theta)=\sum_{k=-d}^d\hat{v}_ke^{2\pi ik\theta}$. Given $E\in \R$,  the eigenequation $\widehat{H}_{v,\alpha,\theta}u=Eu$ induces the $2d$-dimensional cocycle $(\alpha,\widehat{A}_E)$ where
\begin{align}\small\label{sco}
\widehat{A}_{E}(\theta)=\frac{1}{\hat{v}_d}
\begin{pmatrix}
-\hat{v}_{d-1}&\cdots&-\hat{v}_1&E-2\cos2\pi(\theta)-\hat{v}_0&-\hat{v}_{-1}&\cdots&-\hat{v}_{-d+1}&-\hat{v}_{-d}\\
\hat{v}_d& \\
& &  \\
& & & \\
\\
\\
& & &\ddots&\\
\\
\\
& & & & \\
& & & & & \\
& & & & & &\hat{v}_{d}&
\end{pmatrix}.
\end{align}
Let \begin{equation}\label{S}
C=\begin{pmatrix}
\hat{v}_d& \cdots &\hat{v}_1\\
&\ddots&\vdots\\
&&\hat{v}_d
\end{pmatrix},\ \ S=\begin{pmatrix}0&-C^*\\ C&0\end{pmatrix}.
\end{equation}
Since $(\alpha, \widehat{A}_E)$ is Hermitian-symplectic with respect to
$S$ \cite{hp}, its top $d$ Lyapunov exponents are non-negative. We
denote them as $\widehat{L}_1(E)\geq \cdots\geq \widehat{L}_d(E)\geq 0$.

\subsection{Uniform hyperbolicity and dominated splitting}\label{3.3}
Let $Sp(2m,\C)$ be the set of all $2m\times 2m$ Hermitian-symplectic matrices. For  $A\in C^\omega(\T,Sp(2m,\C))$,  we say the  cocycle $(\alpha, A)$ is {\it uniformly hyperbolic} if  there exists a continuous splitting $\C^{2m}=E^s(x)\oplus E^u(x)$ such that for some constants $C>0,c>0$, and for every $n\geqslant 0$,
$$
\begin{aligned}
\lvert A_n(x)v\rvert \leqslant Ce^{-cn}\lvert v\rvert, \quad & v\in E^s(x),\\
\lvert A_n(x)^{-1}v\rvert \leqslant Ce^{-cn}\lvert v\rvert,  \quad & v\in E^u(x+n\alpha).
\end{aligned}
$$
This splitting is left invariant by the dynamics: for every $x \in \T$,
$$
A(x)E^{\ast}(x)=E^{\ast}(x+\alpha),\ \ \ast=s,u.
$$

For  $A\in C^\omega(\T,GL(m,\C))$, we say the cocycle  $(\alpha,A)$  is $k$-dominated (for some $1\leq k\leq m-1$) if there exists a continuous invariant  decomposition $\C^{m}=E_+(x) \oplus E_-(x)$    with $\dim E_+(x)=k$ and $n\in\N$ such that  for any unit vector $v_\pm \in E_\pm(x)\backslash \{0\}$, we have $$\|A_n(x)v_+\|>\|A_n(x)v_-\|.$$

\subsection{Global theory of one-frequency quasiperiodic cocycles}\label{gt} The key  concept for Avila's global theory \cite{avila0} is the acceleration.  If $A\in C^{\omega}(\T,GL(m,\C))$ admits a holomorphic extension to $|\Im z|<\delta$, then for $|\e|<\delta$ we can denote  $A_\e(x)=A(x+i\e)$. The accelerations of $(\alpha,A)$  are defined as
$$
\omega^k(A)=\lim\limits_{\e\rightarrow 0^+}\frac{1}{2\pi\e}(L^k(A_\e)-L^k(A)).
$$
The key ingredient of the global theory is that the acceleration  is quantized.

\begin{Theorem}[\cite{avila0,ajs}]\label{ace}
For each $k$, there exist $1\leq l\leq m$, $l\in\N$, such that $l\omega^k$ is an integer. In particular, if $A\in C^\omega(\T, SL(2,\C))$, then $\omega^1(A)$ is an integer.
\end{Theorem}

\begin{Remark}\label{rem1}
{\rm If $L_j(A)>L_{j+1}(A)$, then $\omega^j(A)$
  is an integer, as follows from the proof of Theorem 1.4 in \cite{ajs}, see also footnote 17 in \cite{ajs}.}
\end{Remark}

We say that $(\alpha,A)$ is {\it $k$-regular} if $\e\rightarrow L^k(A_\e)$ is an affine function of $\e$ in a neighborhood of $0$.

\begin{Theorem}[\cite{avila0,ajs}]\label{t2.1}
Let $\alpha\in\R\backslash\Q$ and $A\in C^\omega(\T,GL(m,\C))$. If $1\leq j\leq m-1$ is such that $L_j(A)>L_{j+1}(A)$, then $(\alpha,A)$ is $j$-regular if and only if $(\alpha,A)$ is $j$-dominated. 
\end{Theorem}

For Schr\"odinger cocycles $(\alpha, A_E)$ with $A_E$ given by
\eqref{Sch}, we have $L(E)=L(A_E)$ and $\omega(E)=\omega(A_E)$ where
$L(E)$ is given by \eqref{multiergodicsch} and $\omega(E)$ is
defined in the Definition \ref{acc}. For $E$ outside the spectrum of $H_{v,\alpha,\theta}$, we have
$L(E)>0, \omega(E)=0.$ For $E$ in the spectrum of $H_{v,\alpha,\theta}$,
Avila \cite{avila0} defines three regimes as modeled by the almost Mathieu operator:
\begin{enumerate}
\item The {\it subcritical} regime: $L(E)=\omega(E)=0$;
\item The {\it critical} regime: $L(E)=0$ and $\omega(E)>0$;
\item The {\it supercritical} regime: $L(E)>0$ and $\omega(E)>0$.
\end{enumerate}
Moreover,  an immediate corollary \cite{gjy1} of the multiplicative Jensen's formula of
\cite{gjyz}, is that for operators \eqref{sch} of Type I, with
$v(\theta)=\sum_{k=-d}^d \hat{v}_ke^{2\pi ik\theta}$ and thus the dual given
by $\widehat{H}_{v,\alpha,\theta}$ we have simplicity of the smallest dual Lyapunov exponent $\widehat L_d$
among $\widehat L_1,\ldots,\widehat L_d$
\begin{center}{\large
\begin{tabular}{c|c|c}
\hline\hline
Regime& {$H_{v,\alpha,x}$}& $\widehat{H}_{v,\alpha,\theta}$\\
\hline
subcritical&$L(E)=0,\omega(E)=0$& $L(E)=0$ and $\widehat{L}_d(E)>0$ is simple  \\
\hline
critical &$L(E)=0,\omega(E)=1$&$L(E)=0$ and $\widehat{L}_d(E)=0$ is simple \\
\hline
supercritical &$L(E)>0,\omega(E)=1$ & $L(E)>0$ and $\widehat{L}_d(E)=0$ is simple\\
\hline
\end{tabular}}
\end{center}

\subsection{The rotation
  number and the IDS}\label{secrot}
Assume $A \in C^\omega(\T, {\rm SL}(2, \R))$ is homotopic to the identity. Then $(\alpha, A)$ induces the projective skew-product $F_A\colon \T \times \mathbb{S}^1 \to \T \times \mathbb{S}^1$
$$
F_A(x,w):=\left(x+\a,\, \frac{A(x) \cdot w}{|A(x) \cdot w|}\right),
$$
which is also homotopic to the identity. Lift $F_A$ to a map $\widetilde{F}_A\colon \T \times \R \to \T \times \R$ of the form $\widetilde{F}_A(x,y)=(x+\alpha,y+\psi_x(y))$, where for every $x \in \T$, $\psi_x$ is $\Z$-periodic.

Map $\psi\colon\T \times \R  \to \R$ is called a {\it lift} of $A$. Let $\mu$ be any probability measure on $\T \times \mathbb{S}^1$ which is invariant by $F_A$, and whose projection on the first coordinate is given by Lebesgue measure.
The number
\begin{equation}\label{rot}
\rho(A):=\int_{\T \times \mathbb S^1} \psi_x(y)\,d\mu(x,y)\pmod{\mathbb Z}
\end{equation}
 depends  neither on the lift $\psi$ nor on the measure $\mu$, and is
 called the \textit{fibered rotation number} of $(\alpha,A)$ (see
 \cite{H,johonson and moser} for more details).
 
 For Schr\"odinger cocycles $(\alpha,A_E)$ with $A_E$ given by
 \eqref{Sch}, we will write
 $\rho(E):=\rho(A_E)$, when $v,\alpha$ are otherwise fixed.

It is well known that $\rho(E)\in[0,\frac{1}{2}]$ is related to the integrated density of states $N(E)$ as follows:
\begin{equation}\label{relation}
N(E)=1-2\rho(E).
\end{equation}

\subsection{Aubry duality}\label{au} 

We consider the following (possibly non-hermitian) quasiperiodic operators,
\begin{equation*}
(H^w_{v,\alpha,x}u)_n=\sum\limits_{k\in\Z} \widehat{v}_k u_{n+k}+w(x+n\alpha)u_n, \ \ n\in\Z,
\end{equation*}
where $v,w$ are two possibly complex-valued 1-periodic analytic functions. 

Consider the fiber direct integral,
$$
\mathcal{H}:=\int_{\T}^{\bigoplus}\ell^2(\Z)dx,
$$
which, as usual, is defined as the space of $\ell^2(\Z)$-valued, $L^2$-functions over the measure space $(\T,dx)$.  The extensions of the
Schr\"odinger operators  and their long-range duals to  $\mathcal{H}$ are given in terms of their direct integrals, which we now define.
Let $\alpha\in\T$ be fixed. Interpreting $H^w_{v,\alpha,x}$ as fibers of the decomposable operator,
$$
H^w_{v,\alpha}:=\int_{\T}^{\bigoplus}H^w_{v,\alpha,x}dx,
$$
the family $\{H^w_{v,\alpha,x}\}_{x\in\T}$ naturally induces an operator on the space $\mathcal{H}$, i.e.,
$$
(H^w_{v,\alpha} \Psi)(x,n)=\sum\limits_{k\in \Z} \hat{v}_k \Psi(x,n+k) + w(x+n\alpha) \Psi(x,n).
$$

Similarly,  the direct integral of long-range operator  $H^v_{w,\alpha,\theta}$,
denoted as $H^v_{w,\alpha}$, is given by
$$
(H^v_{w,\alpha}  \Psi)(\theta,n)=  \sum\limits_{k\in \Z} \hat{w}_k \Psi(\theta,n+k)+v(\theta+n\alpha) \Psi(\theta,n).
$$
Indeed, by analogy with the heuristic and classical approach to Aubry duality \cite{mz,gjls}, let  $U$ be the following operator on $\mathcal{H}:$
\begin{equation}\label{definitionU}
(U\phi)(\eta,m):=\hat{\phi}(m, \eta+m\alpha)=\sum_{n\in\Z}\int_{\T}e^{2\pi imx}e^{2\pi in(m\alpha+\eta)}\phi(x,n)dx.
\end{equation}
$U$ is clearly unitary. Let
$\check v(t):=v(-t)$. A direct computation gives
\begin{equation}\label{equ}
U H^w_{v,\alpha} U^{-1}=H^{\check v}_{w,\alpha}.
\end{equation}.

\section{Multiplicative Jensen's formula for the duals}\label{mult}
Multiplicative Jensen’s formula serves as the foundation of the
duality based quantitative global theory for analytic one-frequency
Schrödinger operators  \cite{gjyz}, having been instrumental in
resolving the almost reducibility conjecture \cite{ge} and the Ten
Martini Problem \cite{gjy1}.

In this section we establish a multiplicative Jensen formula for the duals of analytic
one-frequency Schr\"odinger operators with trigonometric polynomial potential. This
gives a precise description of the first affine segment of the individual complexified
dual Lyapunov exponents and,  in particular, implies that supercriticality of the
Schr\"odinger operator at a given energy yields subcriticality of the dual center.

We consider the following non-Hermitian finite-range quasiperiodic operator
\begin{equation}\label{long}
(\widehat{H}^\e_{v,\alpha,\theta}u)_n=\sum\limits_{k=-d}^{d} \widehat{v}_k u_{n+k}+2\cos2\pi(\theta+i\e+n\alpha)u_n, \ \ n\in\Z,
\end{equation}
where $\hat{v}_k=\overline{\hat{v}_{-k}}$. It is a complexification of the dual
operator \eqref{long-range} of the operator \eqref{sch} with trigonometric
polynomial $v(\theta)=\sum_{k=-d}^d \hat{v}_ke^{2\pi ik\theta}.$

The eigenequation $\widehat{H}^\e_{v,\alpha,\theta}u=Eu$ induces the $2d$-dimensional cocycle $(\alpha,\widehat{A}_E(\cdot+i\e))$ where
\begin{align*}\small
\widehat{A}_{E}(\theta+i\e)=\frac{1}{\hat{v}_d}
\begin{pmatrix}
-\hat{v}_{d-1}&\cdots&-\hat{v}_1&E-2\cos2\pi(\theta+i\e)-\hat{v}_0&-\hat{v}_{-1}&\cdots&-\hat{v}_{-d+1}&-\hat{v}_{-d}\\
\hat{v}_d& \\
& &  \\
& & & \\
\\
\\
& & &\ddots&\\
\\
\\
& & & & \\
& & & & & \\
& & & & & &\hat{v}_{d}&
\end{pmatrix}.
\end{align*}
Following the definition in Section \ref{3.2}, for $1\leq k\leq 2d$,
we denote $\widehat{L}^k_\e(E):=L^k(\widehat{A}_E(\cdot+i\e))$ and $\widehat{L}_k(E):=L_k(\widehat{A}_E).$ Let the multiplicity of $\widehat{L}_d(E)$ among the first $d$ exponents  be $n_d.$ If $n_d<d$ set
$$
\delta(E)=\widehat{L}_{d-n_d}(E)-\widehat{L}_d(E).
$$ 
If $n_d=d$, set $\delta(E)=0$.  Recall that $L(E):=L(A_E)$ where $A_E$ is defined in \eqref{Sch}. 
\begin{Theorem}\label{strip}
For $E\in\R, \alpha\in\R\backslash\Q$, 
\begin{align*}
\widehat{L}_\e^{d}(E)=\left\{
\begin{aligned} & \widehat{L}^d_0(E)&|\e|\in
[0,L(E)/2\pi],\\
&\widehat{L}^d_{L(E)/2\pi}(E)+(2\pi |\e|-L(E))
&|\e|\in (L(E)/2\pi,\infty).
\end{aligned}\right.
\end{align*}
Moreover, for $1\leq k\leq d-n_d$,

$$
\widehat{L}^k_\e(E)=\widehat{L}^k_0(E),\ \ \forall |\e|\in [0,(\delta(E)+L(E))/2\pi];
$$
 if $L(E)>0$, then for $1\leq k\leq d$,
$$
\widehat{L}^k_\e(E)=\widehat{L}^k_0(E),\ \ \forall |\e|\in [0,L(E)/2\pi].
$$

\end{Theorem}
\begin{pf}
For any $E\in\C$,  one always has 
$$
\widehat{L}_\e^{d}(E)\leq \sup_{\theta\in \T} \ln \|\Lambda^d \widehat{A}_E(\theta+i\e) \| \leq 2\pi|\e| + O(1).
$$
Thus, by convexity, for any $E\in\C$, the absolute value of the slope of $\widehat{L}_\e^{d}(E)$ as a function of $\e$ is less than or equal to $2\pi$.  
By direct computation,  for  sufficiently large $\e$,
\begin{equation*}
		(\widehat{A}_E)_d(\theta+i\e)=e^{2\pi |\e|}e^{-2\pi i \theta}\begin{pmatrix}
			-C^{-1}D&0\\0 &0
		\end{pmatrix}+ O(1)
	\end{equation*}
where
$$
C=\begin{pmatrix}
\hat{v}_d& \cdots &\hat{v}_1\\
&\ddots&\vdots\\
&&\hat{v}_d
\end{pmatrix},\ \ D=\begin{pmatrix}e^{-2\pi i (d-1)\alpha}&&\\ &\ddots& \\ &&1\end{pmatrix}.
$$
For $\e\to-\infty$ conjugate the phase factors in the leading term. By continuity of the Lyapunov exponent \cite{ajs}, we have
	\begin{equation*}
	\widehat{L}_\e^{d}(E)= 2\pi |\e| -\ln|\widehat{v}_d|+o(1).
	\end{equation*}
	Thus by Theorem \ref{ace}, 
	\begin{equation}\label{new3}
		\widehat{L}_\e^{d}(E)= 2\pi |\e|-\ln |\widehat{v}_d |  \quad \text{as $|\e| \rightarrow \infty$,}
	\end{equation} 
i.e. the slope of $\widehat{L}^d_\e(E)$ is $\pm 2\pi$, as $|\e| \rightarrow \infty$. 

For $1\leq k\leq 2d$, set
$$
\omega_\pm^k(A):=\lim\limits_{\e\rightarrow 0^\pm}\frac{L^k(A(\cdot+i\e))-L^k(A)}{2\pi\e},\ \ \omega_\pm^k(E)=\omega_\pm^k(\widehat{A}_E),
$$
$$
t^+_k(E):=\sup\left\{\e: \e\geq 0, \ \ L_k(\widehat{A}_E(\cdot+i\e))=\widehat{L}_k(E)\right\},
$$
$$
t^-_k(E):=\sup\left\{-\e:\e\leq 0,\ \ L_k(\widehat{A}_E(\cdot+i\e))=\widehat{L}_k(E)\right\}.
$$

For nonreal $E$, $(\alpha,\widehat A_E)$ is uniformly hyperbolic with unstable rank $d$ \cite{hp}. Thus

$$
\widehat{L}_d(E)>0>\widehat{L}_{d+1}(E)
$$
which by Remark \ref{rem1} implies that $\omega_\pm^d(E)$ is an integer \footnote{Although Remark \ref{rem1} only states the result for $\omega^d_+(E)$, all results in \cite{ajs} work for  both $\omega^d_\pm(E)$.}. On the other hand, taking exterior powers in \eqref{strip2} gives $d\widehat L^j_\e(E)\le 2\pi j|\e|+O(1)$, so convexity gives 
\begin{equation}\label{less than 1}
|\omega_\pm^k(\widehat{A}_E(\cdot+i\e))|<1,\ \ 1\leq k\leq d-1,
\end{equation}
\begin{equation}\label{equal 1}
|\omega_\pm^d(\widehat{A}_E(\cdot+i\e))|\leq 1
\end{equation}
for any $\e\in\R$. By $d$-domination the signed derivatives at zero agree. If their common normalized value is $1$ or $-1$, convexity and \eqref{equal 1} give, along the positive or negative ray, respectively,
	\begin{equation*}
		\widehat{L}_\e^{d}(E)= \widehat{L}_0^{d}(E)+2\pi |\e|,
	\end{equation*} 
which by \eqref{new3} implies that 
\begin{equation}\label{LE11}
\widehat{L}_0^{d}(E)=-\ln|\widehat{v}_d|.
\end{equation} 
By Haro and Puig \cite{hp},
we have for any $E\in\C$,
\begin{equation}\label{hp1}L(E)=\widehat{L}^d_0(E)+\ln|\hat{v}_d|.
  \end{equation}

By  \eqref{LE11} and \eqref{hp1}, we have $L(E)=\widehat{L}_0^{d}(E)+\ln|\widehat{v}_d|=0$, contradicting  $E\in\C\backslash\R$ \footnote{In this case $(\alpha,A_E)$ is uniformly hyperbolic, thus $L(E)>0$. On each outward ray the sum through $d$ is convex and initially constant. At common differentiability points every sum ending at a strict gap below $d$ has zero slope by \eqref{less than 1} and Remark~\ref{rem1}. Tied exponents have equal derivatives by ordering. Subtracting the sums at the ends of each tied block shows that only the block ending at $d$ can increase, at a common rate. By \eqref{new3} and the bound before \eqref{less than 1}, $L_d(\widehat A_E(\cdot+i\e))\ge 2\pi|\e|/d-O(1)$. Thus each $t_k^\pm$ with $k\le d$ ends an initial plateau, with positive outgoing slope by piecewise affinity. This also starts the induction, with the empty sum taken as zero.\label{blo}}. Hence $\omega_\pm^d(E)=0$. By Remark~\ref{rem1} and \eqref{less than 1}, for every  $k<d$ with $\widehat{L}_k(E)>\widehat{L}_{k+1}(E)$, $\omega_\pm^k(E)=0$.

For $E\in\C\backslash\R$ and $1\leq k\leq d-1$, we further claim 
\begin{enumerate}
\item $t^\pm_{k+1}(E)\begin{cases}
<t^\pm_{k}(E), & \widehat{L}_k(E)>\widehat{L}_{k+1}(E)\\
=t^\pm_{k}(E), & \widehat{L}_k(E)=\widehat{L}_{k+1}(E)
\end{cases}$;
\item $\omega^{k}_\pm(E)=0$;
\item $\sigma\bigl[\omega_\sigma^{k+1}(\widehat A_E(\cdot+\sigma i t^\sigma_{k+1}(E)))-\omega_\sigma^k(\widehat A_E(\cdot+\sigma i t^\sigma_{k+1}(E)))\bigr]>0,\quad\sigma=\pm1$.
\end{enumerate}
We only prove (1)-(3) for the $+$ case (the $-$
case is proved in exactly the same way) and omit $+$ in the notations for
simplicity.

Assume for $k\leq m-1$, the (1)-(3) above are true. Then for $k=m< d$, we distinguish two cases:

{\bf Case I:} $\widehat{L}_m(E)>\widehat{L}_{m+1}(E)$. For (1), if $t_{m+1}(E)\geq t_m(E)$, then by the definition of $t_m(E)$ and $t_{m+1}(E)$,
$$
L_m(\widehat{A}_E(\cdot+it_m(E)))=\widehat{L}_m(E)>\widehat{L}_{m+1}(E)=L_{m+1}(\widehat{A}_E(\cdot+i t_m(E))).
$$
Thus by Remark \ref{rem1}, we have $\omega^m(\widehat{A}_E(\cdot+i t_m(E)))\in\Z$. By \eqref{less than 1}, we have $|\omega^m(\widehat{A}_E(\cdot+it_m(E)))|<1$. Thus $\omega^m(\widehat{A}_E(\cdot+i t_m(E)))=0$. On the other hand, by induction ((3) is true for $k=m-1$) and convexity of $\widehat{L}^{m-1}_\e(E)$, we have 
$$
\omega^m(\widehat{A}_E(\cdot+i t_m(E)))>\omega^{m-1}(\widehat{A}_E(\cdot+i t_m(E)))\geq 0
$$
which is a contradiction. Thus $t_{m+1}(E)< t_m(E)$. 

For (2), since $\widehat{L}_m(E)>\widehat{L}_{m+1}(E)$, by Remark \ref{rem1} and \eqref{less than 1}, we have $\omega^m(E)\in\Z$ and $|\omega^m(E)|<1$. Therefore $\omega^m(E)=0$. 

For (3), note first that (1), (2) and the definition of $t_m(E)$ imply that  $\omega^{m}(\widehat{A}_E(\cdot+it_{m+1}(E)))=0$. By the convexity of $\widehat{L}^{m+1}_\e(E)$ and the definition of $t_{m+1}(E)$, we have
$$
\omega^{m+1}(\widehat{A}_E(\cdot+it_{m+1}(E)))-\omega^{m}(\widehat{A}_E(\cdot+it_{m+1}(E)))=\omega^{m+1}(\widehat{A}_E(\cdot+it_{m+1}(E)))> 0.
$$

{\bf Case II:} $\widehat{L}_m(E)=\widehat{L}_{m+1}(E)$. For (2), if there is a sequence $\e_n\rightarrow 0^+$  such that 
$$
L_m(\widehat{A}_E(\cdot+i\e_n))\neq L_{m+1}(\widehat{A}_E(\cdot+i\e_n)),
$$ 
then $\omega^m(E)\in\Z$, thus by the same argument as above, we have $\omega^m(E)=0,$ which  implies (2). Otherwise, there is $\delta_0>0$ such that 
\begin{equation}\label{new5}
L_m(\widehat{A}_E(\cdot+i\e))=L_{m+1}(\widehat{A}_E(\cdot+i\e)),\ \ 0\leq \e\leq \delta_0.
\end{equation}
In this case, if there is a sequence $\e_n\rightarrow 0^+$ such that 
$$
L_{m+1}(\widehat{A}_E(\cdot+i\e_n))\neq L_{m+2}(\widehat{A}_E(\cdot+i\e_n)),
$$ 
then $\omega^{m+1}(E)=0,$ which together with \eqref{new5} implies by induction that 
$$
\omega^{m}(E)=\omega^{m-1}(E)+\frac{\omega^{m+1}(E)-\omega^{m-1}(E)}{2}=0.
$$ 
Otherwise, there is $0<\delta_1<\delta_0$ such that 
$$
L_m(\widehat{A}_E(\cdot+i\e))=L_{m+1}(\widehat{A}_E(\cdot+i\e))=L_{m+2}(\widehat{A}_E(\cdot+i\e)),\ \ 0\leq \e\leq \delta_1.
$$ 
Continuing this process, we  finally arrive at the case
$$
L_m(\widehat{A}_E(\cdot+i\e))=\cdots=L_{d}(\widehat{A}_E(\cdot+i\e)),\ \ 0\leq \e\leq \delta'.
$$ 
 Then, by  induction,
$$
\omega^{m}(E)=\omega^{m-1}(E)+\frac{\omega^{d}(E)-\omega^{m-1}(E)}{d-m+1}=0
$$
which completes the proof.

For (1), if $t_{m+1}(E)>t_m(E)$, take $t_{m}(E)<\e_1<t_{m+1}(E)$ sufficiently
close to $t_m(E)$. By induction ((3) is true for $k=m-1$), we have
$$
L_m(\widehat{A}_E(\cdot+i\e_1))>L_m(\widehat{A}_E(\cdot+it_m(E))).
$$
By the definition of $t_m(E)$ and $t_{m+1}(E)$, 
$$
L_m(\widehat{A}_E(\cdot+it_m(E)))=\widehat{L}_m(E)=\widehat{L}_{m+1}(E)=L_{m+1}(\widehat{A}_E(\cdot+i\e_1)),
$$
hence
$$
 L_{m}(\widehat{A}_E(\cdot+i\e_1))>L_{m+1}(\widehat{A}_E(\cdot+i\e_1)).
 $$ 
By a similar argument to the above, we have $\omega^m(\widehat{A}_E(\cdot+i\e_1))=0$. Again, by induction on (3) and convexity of $\widehat{L}^{m-1}_\e(E)$, we have 
$$
\omega^m(\widehat{A}_E(\cdot+i \e_1))=\omega^m(\widehat{A}_E(\cdot+it_m(E)))>\omega^{m-1}(\widehat{A}_E(\cdot+i t_m(E)))\geq 0
$$
which is a contradiction.

If $t_{m+1}(E)<t_m(E)$, then $L_{m+1}(\widehat{A}_E(\cdot+i\e))>\widehat{L}_{m+1}(E)=L_{m}(\widehat{A}_E(\cdot+i\e))$ just after $t_{m+1}(E)$, contradicting the ordering. Hence $t_{m+1}(E)=t_m(E)$.

For (3), note that, by induction, we have 
\begin{equation}\label{new1}
\omega^{m}(\widehat{A}_E(\cdot+it_{m}(E)))-\omega^{m-1}(\widehat{A}_E(\cdot+it_{m}(E)))>0.
\end{equation}
If we further assume
\begin{equation}\label{new2}
\omega^{m+1}(\widehat{A}_E(\cdot+it_{m+1}(E)))-\omega^{m}(\widehat{A}_E(\cdot+it_{m+1}(E))) \leq 0,
\end{equation}
then we take $\e_1=t_m(E)+\delta=t_{m+1}(E)+\delta$ where $\delta>0$
is sufficiently small. Then \eqref{new1} and \eqref{new2} imply
$$
L_m(\widehat{A}_E(\cdot+i\e_1))>L_m(\widehat{A}_E(\cdot+it_{m}(E)))\geq L_{m+1}(\widehat{A}_E(\cdot+i\e_1)).
$$
Thus $\omega^{m}(\widehat{A}_E(\cdot+i\e_1))=0$ which contradicts the definition of $t_m(E)$ and convexity of $L^{m}_\e(E)$.

For the remaining results, we distinguish two cases:

{\bf Case I:} $L(E)=0.$ Then $|\omega_\pm^d(E)|=1$. Indeed, otherwise, if $|\omega_\pm^d(E)|<1$, by convexity of $\widehat{L}^d_\e(E)$ and Theorem \ref{ace} \footnote{I.e. $\widehat L_\e^d(E)$ is a convex piecewise affine function with final outward slope $2\pi$.}, there is $f_\pm(E)>0$ on that ray such that for $\e$ sufficiently large, we have
$$
\widehat{L}_\e^d(E)=\widehat{L}_0^d(E)-f_\pm(E)\pm 2\pi \e.
$$
Together with \eqref{new3}, we have $\widehat{L}_0^d(E)-f_\pm(E)=-\ln|\hat{v}_d|$ and hence $L(E)=\widehat{L}^d_0(E)+\ln|\hat{v}_d|=f_\pm(E)>0,$ which is a contradiction.

 For real $E$, \eqref{new111} makes $\widehat L^d_\e(E)$ even. Then by convexity, $\omega^d(\widehat{A}_E(\cdot\pm i\e))=\pm 1$ for any $\e>0$. Thus,
	\begin{equation}\label{fgj1}
		\widehat{L}_\e^{d}(E)= \widehat{L}_0^{d}(E)+2\pi |\e|.
	\end{equation} 

{\bf Case II:} If $E\in\R$ and $L(E)>0$, there is a sequence
$E_n\rightarrow E$ with $E_n\in\C_+$, such that $\omega_\pm^d(E_n)=\omega_\pm^d(\overline{E_n})=0$ and (1)-(3) are true for all $E_n$ and $\overline{E_n}$. Note that  we always have
\begin{equation}\label{new111}
\widehat{A}_{E}(\theta+i\e)^*S\widehat{A}_{\bar{E}}(\theta-i\e)=S.
\end{equation}
 In particular $\widehat L^d_\e(z)=\widehat L^d_{-\e}(\bar z)$, since $|\det\widehat A_z|=1$. {By (1), both sums end their initial plateaux when their respective $d$th exponents do. Hence}
 \begin{equation}\label{lana1}
 t_d^+(E_n)=t_d^-(\overline{E_n}).
 \end{equation}  
 By \eqref{new111} we have 
$$
L_{d+1}(\widehat{A}_{E_n}(\cdot+it_d^+(E_n)))=-L_d(\widehat{A}_{\overline{E_n}}(\cdot-it_d^+(E_n)))\leq 0
$$ 
Hence 
$$
L_{d}(\widehat{A}_{E_n}(\cdot+it^+_d(E_n)))>0\geq L_{d+1}(\widehat{A}_{E_n}(\cdot+it^+_d(E_n))).
$$ 
It follows that $\omega^d(\widehat{A}_{E_n}(\cdot+it^+_d(E_n)))\in\Z.$
By  convexity of $\widehat{L}^d_\e(E_n)$, the definition of
$t^+_d(E_n)$ and \eqref{equal 1},
$0<\omega^d(\widehat{A}_{E_n}(\cdot+it^+_d(E_n)))\leq 1$, so we have $\omega^d(\widehat{A}_{E_n}(\cdot+it^+_d(E_n)))=1$. Hence  by  convexity of $\widehat{L}_\e^d(E_n)$, we have $\omega^d(\widehat{A}_{E_n}(\cdot+i\e))=1$ for any $\e\geq t^+_d(E_n)$. Together with \eqref{new3}, it follows that
$$
\widehat{L}^d_\e(E_n)=\widehat{L}^d_0(E_n)+2\pi (\e-t^+_d(E_n))=2\pi \e-\ln|\widehat{v}_d|,
$$
thus $t^+_d(E_n)=L(E_n)/(2\pi)$. On the other hand,
$\widehat{L}^d_{-\e}(\overline{E_n})$ is a piecewise convex affine
function of $\e\geq 0$ with final slope $2\pi$, so there is $f(\overline{E_n})>0$ such that for $\e$ sufficiently large, we have
$$
\widehat{L}_{-\e}^d(\overline{E_n})=\widehat{L}_0^d(\overline{E_n})-f(\overline{E_n})+2\pi \e,
$$
Moreover, by (1)--(3) for $\overline{E_n}$, we have $f(\overline{E_n})\ge2\pi t_d^-(\overline{E_n})$\footnote{The equality holds if and only if $-t_d^-(\overline{E_n})$ is the only turning point of $\widehat L^d_\e(\overline{E_n})$ when $\e\le0$.}. Hence by \eqref{lana1}, we have $L(E_n)/\mathopen{(}2\pi\mathclose{)}=t_d^+(E_n)\le t_d^-(\overline{E_n})\le f(\overline{E_n})/\mathopen{(}2\pi\mathclose{)}=L(\overline{E_n})/\mathopen{(}2\pi\mathclose{)}$. The positive-ray formula, including its initial constant part, now passes to the limit at each fixed $\e\ge0$. By continuity \cite{ajs,bj} and the reflected identity,
\[
 \widehat L^d_{\pm\e}(E)=\widehat L^d_0(E)\mathbin{+}\max\{2\pi\e-L(E),0\} .
\]
In both cases the sum is nondecreasing on each ray. By the block argument of Footnote \ref{blo} its plateau ends at $t_d^+(E)=t_d^-(E)=L(E)/\mathopen{(}2\pi\mathclose{)}$ and all first $d$ exponents are constant there. If $n_d=d$, this proves the theorem. Otherwise the gap at $d-n_d$ persists beyond both endpoints. By \eqref{less than 1} and the block argument, the first $d-n_d$ exponents remain constant up to $t_{d-n_d}^\pm(E)$. Thus, for $1\le k\le d-n_d$,
\begin{equation}\label{fgj2}
\widehat{L}^k_\e(E)=\widehat{L}^k_0(E),\ \ \forall \e\in[-t^-_{d-n_d}(E),t^+_{d-n_d}(E)].
\end{equation}
If $t^\pm_{d-n_d}(E)<\frac{L(E)+\widehat L_{d-n_d}(E)-\widehat L_d(E)}{2\pi}=t^\pm_d(E)+\frac{\widehat L_{d-n_d}(E)-\widehat L_d(E)}{2\pi}$, then
\[
L_{d-n_d+1}(\widehat A_E(\cdot\mathbin{\pm} it^\pm_{d-n_d}(E)))\le L_{d-n_d+1}(\widehat A_E(\cdot\mathbin{\pm} it^\pm_d(E)))+2\pi(t^\pm_{d-n_d}(E)-t^\pm_d(E))<\widehat L_{d-n_d}(E).
\]
We therefore have that $\omega^{d-n_d}_\pm(\widehat A_E(\cdot\mathbin{\pm} it^\pm_{d-n_d}(E))\mathclose{)}$ is an integer, thus, by \eqref{less than 1} it must be $0$. However this contradicts the definition of $t^\pm_{d-n_d}(E)$. This completes the proof.
\end{pf}

\section{Symplectic structure and convergence of the dual center}
In this section, we use the multiplicative Jensen formula from Section 4 to construct analytic symplectic frames for the centers of finite-range dual cocycles. We then prove convergence
of the center dynamics under trigonometric polynomial approximation. The limit gives
the intrinsic dual center for general analytic potentials and
serves as the key step in the proof of Theorem
\ref{hss}. 

We define 
$$
(H^\epsilon_{v(\cdot+i\e),\alpha} \Psi)(x,n)= e^{-2\pi\epsilon}\Psi(x,n+1)+ e^{2\pi\epsilon}\Psi(x,n-1) +  v(x+i\e+n\alpha) \Psi(x,n).
$$
and its Aubry dual
$$
(\widehat{H}^\epsilon_{v(\cdot+i\e),\alpha}  \Psi)(\theta,n)=  \sum\limits_{k\in \Z} e^{-2\pi k\e}\hat{v}_k \Psi(\theta,n+k)+  2\cos2\pi (\theta+i\epsilon+n\alpha) \Psi(\theta,n).
$$
Indeed, it follows by a direct computation that
\begin{equation}\label{equ}
U H^\epsilon_{v(\cdot+i\e),\alpha} U^{-1}=  \widehat{H}^{-\epsilon}_{v(\cdot+i\e),\alpha}.
\end{equation}
where $U$ is given by \eqref{definitionU}.

\begin{Definition}[Type k]\label{typeik}
  {\rm For any $k\in\N$,  we say $E$ is a {\it Type k energy} for operator
    $H_{v,\alpha,\theta}$ if
  $\bar{\omega}(E)= k$. We say $H_{v,\alpha,\theta}$ is a {\it Type k operator,}
  if every $E$ in the spectrum of
  $H_{v,\alpha,\theta}$ is Type k.}
\end{Definition}
For any $E\in\R$, $\delta$,  $h>0$, let 
\begin{equation}\label{ha4}
R_{\delta}(E)\times\Omega_h=\{z\in\C:\max\{|\Re z-E|,|\Im z|\}<\delta\}\times \{z\in\C:e^{-2\pi h}<|z|<e^{2\pi h}\}.
\end{equation} 
Let $V$ be a complex
holomorphic vector bundle over $\R_\delta(E)\times \Omega_h$ and let $\pi:V\rightarrow R_\delta(E)\times \Omega_h$ be the bundle projection. A holomorphic vector bundle $V$ of rank $r$ over $R_\delta(E)\times\Omega_h$ is called trivial if it is isomorphic to the bundle $(R_\delta(E)\times\Omega_h)\times \C^r$. This is equivalent to the existence of a global frame $v_1,\cdots,v_r$ of holomorphic sections  in $V$ over $R_\delta(E)\times\Omega_h$, such that for each $(z_1,z_2)\in R_\delta(E)\times \Omega_h$,  the elements $v_1(z_1,z_2),\cdots,v_r(z_1,z_2)\in \pi^{-1}(z_1,z_2)$ are linearly independent.
\begin{Theorem}\label{trivial}
Any holomorphic vector bundle $V$ over $R_\delta(E)\times\Omega_h$ is trivial.
\end{Theorem}
\begin{pf}
Since $R_\delta(E)$ is contractible and $\Omega_h$ is homotopy equivalent to $\mathbb{S}^1$, then $R_\delta(E)\times \Omega_h$ is homotopy equivalent to $\mathbb{S}^1$, which implies that $V$ is topologically trivial. On the other hand,  $R_\delta(E)\times \Omega_h$ is a product of Stein manifold, thus is a Stein manifold itself. It follows from Theorem 5.3.1 in \cite{stein} that $V$ is trivial.
\end{pf}

\subsection{The trigonometric polynomial case}\label{4.21}
For this subsection,  let 
$$
v(x)=\sum_{m=-d}^d \hat{v}_m e^{2\pi imx},\ \ \hat{v}_m=\overline{\hat{v}_{-m}}
$$ 
be a real trigonometric polynomial of degree $d$.  For simplicity, we denote the dual of $H_{v,\alpha,x}$ by $\widehat{H}_{v,\alpha,\theta}$. 
Let $\Sigma_{v,\alpha}$ be the spectrum of $H_{v,\alpha,x}$.
The following is the dual characterization of Type $k$ operators, we denote $\widehat{L}_0(E)=\infty$,
\begin{Proposition}\label{simple le}
For any $1\leq k\leq d$, $H_{v,\alpha,x}$ is a Type $k$  operator if and
only if $\widehat{H}_{v,\alpha,\theta}$ is $\mathcal{PH}_{2k}$, in the sense that for all
 $E\in\Sigma_{v,\alpha}$,
\begin{enumerate}
\item{$\widehat{L}_{d-k}(E)>\widehat{L}_{d-k+1}(E)=\cdots=\widehat{L}_d(E)$;}
\item {$(\alpha,\widehat{A}_{E})$ is $(d-k)$ and $(d+k)$-dominated.}
\end{enumerate}
\end{Proposition}
\begin{pf}
If $k=d$, the results follow from the finite Jensen formula, both domination conditions being vacuous. For $k<d$, by Theorem 1.1 in \cite{gjyz}, $\bar{\omega}(E)=k$ if and only if $\widehat{L}_{d-k}(E)>\widehat{L}_{d-k+1}(E)=\cdots=\widehat{L}_d(E)$.

We let
$$
C=\begin{pmatrix}
\hat{v}_d&\cdots&\hat{v}_1\\
0&\ddots&\vdots\\
0&0&\hat{v}_d
\end{pmatrix},\ \ B(\theta)=\hat{v}_0 I_d+\begin{pmatrix}
2\cos2\pi(\theta_{d-1})&\hat{v}_{-1}&\cdots&\hat{v}_{-d+1}\\
\hat{v}_1&\ddots&\ddots&\vdots\\
\vdots&\ddots&2\cos2\pi(\theta_1)&\hat{v}_{-1}\\
\hat{v}_{d-1}&\cdots&\hat{v}_1&2\cos2\pi(\theta)
\end{pmatrix}
$$
where $\theta_j=\theta+j\alpha$.  Then one can check that
\begin{equation}\label{strip2}
\widehat{A}_{E}(\theta+(d-1)\alpha)\cdots \widehat{A}_{E}(\theta)=:\widehat{A}_{d,E}(\theta)=\begin{pmatrix}
C^{-1}(EI-B(\theta))& -C^{-1}C^*\\
I_d&O_d
\end{pmatrix}
\end{equation}
where $I_d$ and $O_d$ are the $d$-dimensional identity and zero matrices, respectively.

Notice that \eqref{strip2} implies that we always have
$$
dL^{d-k}(\widehat{A}_{E})=L^{d-k}(\widehat{A}_{d,E}).
$$
Thus by the definition of regularity, $(\alpha,\widehat{A}_{E})$ is
($d-k$)-regular if and only if  $(d\alpha,\widehat{A}_{d,E})$ is
($d-k$)-regular. Let $\left(\ell_{ij}\right):=(\widehat{A}_{d,E})_n(\theta)$. It is easy to check that each
$\ell_{ij}$ is a Laurent polynomial  with degree $\leq
n$. Similarly,  let $L_{ij}$ be the $ij$-th entry of
$\Lambda^{d-k}(\widehat{A}_{d,E})_n(\theta).$ By the definition of wedge
product, each $L_{ij}$ is a Laurent polynomial  of degree $\leq n(d-k)$. Taking $n=1$ bounds the Lyapunov sum by $2\pi(d-k)|\e|+O(1)$. Hence by convexity
\begin{align*}
&|\omega^{d-k}(\widehat{A}_{d,E})|=\left|\lim\limits_{\e\rightarrow 0^+}\frac{1}{2\pi\e}(L^{d-k}(\widehat{A}_{d,E}(\cdot+i\e))-L^{d-k}(\widehat{A}_{d,E}))\right|\\
=&\lim\limits_{\e\rightarrow 0^+}\frac{1}{2\pi\e}\left|\lim\limits_{n\rightarrow \infty}\frac{1}{n}\int_{\T}\ln\|\Lambda^{d-k}\left(\widehat{A}_{d,E}\right)_n(\theta+i\e)\|d\theta-\lim\limits_{n\rightarrow \infty}\frac{1}{n}\int_{\T}\ln\|\Lambda^{d-k}\left(\widehat{A}_{d,E}\right)_n(\theta)\|d\theta\right|\\
\leq& d-k.
\end{align*}
It follows that
$$
|\omega^{d-k}(\widehat{A}_{E})|=\left|\frac{\omega^{d-k}(\widehat{A}_{d,E})}{d}\right|\leq \frac{d-k}{d}<1.
$$
Note that $\widehat{L}_{d-k+1}(E)<\widehat{L}_{d-k}(E).$ Thus by
Remark  \ref{rem1}, $\omega^{d-k}(\widehat{A}_{E})$ is an
integer. Thus, since $|\omega^{d-k}(\widehat{A}_{E})|$ is strictly smaller than $1$, we have $\omega^{d-k}(\widehat{A}_{E})=0$.  This implies that
$$
L^{d-k}(\widehat{A}_{E}(\cdot+i\e))=L^{d-k}(\widehat{A}_{E})
$$
for sufficiently small $\e>0$. A similar  argument  works for  $\e<0$.
This means  $(\alpha,\widehat{A}_{E})$ is ($d-k$)-regular, hence, by Theorem
\ref{t2.1}, $(\alpha,\widehat{A}_{E})$ is ($d-k$)-dominated. Since $(\alpha,\widehat{A}_{E})$
is Hermitian-symplectic, we have $(\alpha,\widehat{A}_{E})$ is ($d+k$)-dominated.

\end{pf}
In the following, we explore further properties of $\widehat{H}_{v,\alpha,\theta}$ which will be important for our applications. As a corollary of Theorem \ref{strip}, we have
\begin{Corollary}\label{dominate1}
Assume $\alpha\in \R\backslash\Q,$ $k\geq 1$,  and $E\in\R$ is a Type k energy
of operator $H_{v,\alpha,x}.$ Then, for $k<d$, $(\alpha,\widehat{A}_E(\cdot+i\e))$ is $(d-k)$ and $(d+k)$-dominated for any $|\e|<(\delta(E)+L(E))/2\pi$ where $\delta(E)=\widehat{L}_{d-k}(E)-\widehat{L}_d(E)$. For $k=d$, both conditions are vacuous for every $\e$.
\end{Corollary}
\begin{pf}
It follows directly from
Theorem \ref{strip} that $(\alpha,\widehat{A}_E(\cdot+i\e))$ is $(d-k)$
and $(d+k)$-regular by reflection \eqref{new111} for any $|\e|<(\delta(E)+L(E))/2\pi.$ For $k<d$, the bound on $L_{d-k+1}(\alpha,\widehat{A}_E(\cdot+i\e))$ in its proof gives
$L_{d-k}(\alpha,\widehat{A}_E(\cdot+i\e))>L_{d-k+1}(\alpha,\widehat{A}_E(\cdot+i\e))$. Reflection gives the gap for $\e<0.$ For $k=d$ both conditions are vacuous.  For $k<d$, by  (1) of Proposition \ref{simple le} and Theorem \ref{t2.1}, $(\alpha,\widehat{A}_E(\cdot+i\e))$ is $(d-k)$ and $(d+k)$-dominated for any $|\e|<(\delta(E)+L(E))/2\pi$.
\end{pf}

Hence for $E$ with $\bar{\omega}(E)=k$  there exists a continuous invariant decomposition
\begin{equation}\label{eqnew1}
\C^{2d}=E_s(\theta)\oplus E_c(\theta)\oplus E_u(\theta),\ \ \forall \theta\in\T,
\end{equation}
where $E_c(\theta)$ is the  $2k$-dimensional invariant subspace
corresponding to the exponents $\pm\widehat L_d(E)$.  We
denote  by $C^\pm(\T, *)$ the set of all $(*)$-valued functions such
that $f$ is holomorphic inside/outside the unit circle and can be extended
continuously to the unit circle. Let also $C_h^\pm(\T,*)$ be the set of all $(*)$-valued functions that are analytic on $\{0<\pm\Im \theta<h\}$ and can be extended continuously to $\T$.  Fix $0<h<(\widehat{L}_{d-k}(E)-\widehat{L}_d(E)+L(E))/2\pi$.   Involving  the Hermitian-symplectic structure of the bundles, we  actually have a {\it symplectic}  invariant decomposition of $\C^{2d}$, depending analytically on $\theta$.

\begin{Lemma}\label{uvec}
There are $2k$ linearly independent $u^i_E(\theta)\in E_c(\theta)$ with $u^i_E\in C^\omega_h(\T,\C^{2d})$, $1\leq i\leq 2k$,  such that
\begin{equation}\label{uveceq}
O_E(\theta)^*SO_E(\theta)=
\begin{pmatrix}
0&I_k\\
-I_k&0
\end{pmatrix}:=J_{2k}.
\end{equation}
where $O_E(\theta)=\begin{pmatrix}u_E^1(\theta)&u_E^2(\theta)&\cdots&u_E^{2k}(\theta)\end{pmatrix}$.
\end{Lemma}
\begin{pf}
Note that by Theorem 6.1 in \cite{ajs}, the complex vector bundle $E_c(\theta)$,   $E_s(\theta)$ and $E_u(\theta)$ are holomorphic over $\Omega_h$. By Theorem \ref{trivial}, there are global holomorphic frames $\{f_E^j(\theta)\}_{j=1}^{d-k}\in E_s(\theta)$,  $\{\tilde{u}^j_E(\theta)\}_{j=1}^{2k}\in E_c(\theta)$  and $\{g_E^j(\theta)\}_{j=1}^{d-k}\in E_u(\theta)$ respectively.  

Let
 $$
 \tilde{O}_E(\theta)=\begin{pmatrix}\tilde{u}_E^1(\theta)&\tilde{u}_E^2(\theta)&\cdots&\tilde{u}_E^{2k}(\theta)\end{pmatrix},
 $$
 $$
 P_E(\theta)=\begin{pmatrix}f_E^1(\theta)&\cdots&f_E^{d-k}(\theta)&\tilde{u}_E^1(\theta)&\tilde{u}_E^2(\theta)&\cdots&\tilde{u}_E^{2k}(\theta)&g_E^1(\theta)&\cdots&g^{d-k}_E(\theta)\end{pmatrix},
 $$
 \begin{equation}\label{eiga}
\tilde{\Omega}_E(\theta)=\tilde{O}_E(\theta)^*S\tilde{O}_E(\theta), \ \ \Lambda_E(\theta)=P_E(\theta)^*SP_E(\theta)=\begin{pmatrix}&&A_E(\theta)\\
&\tilde{\Omega}_E(\theta)&\\ -A_E(\theta)^*&&\end{pmatrix} \footnote{$\Lambda_E(\theta)$ has this form because of the symplectic orthogonality.}.
\end{equation}
Since $iS$ has $d$ positive and $d$ negative eigenvalues, by \eqref{eiga},  we have $i\tilde{\Omega}_E(\theta)$ has $k$-positive and $k$-negative eigenvalues for all $\theta\in\T$.

By Theorem 1.15 and Theorem 2.12 in \cite{gks}, there exists $G_E^+\in C^+(\T,GL(2k,\C))$ such that
 for real $\theta$,
\[
G_E^+(\theta)^*i\tilde\Omega_E(\theta)G_E^+(\theta)=D_E(\theta),
\]
where
\[
D_E(\theta)=
\begin{pmatrix}
&&&&&&&e^{2\pi i k_1\theta}I_{m_1}\\
&&&&&&\begin{sideways}$\ddots$\end{sideways}&\\
&&&&&e^{2\pi i k_{p-1}\theta}I_{m_{p-1}}&&\\
&&&I_{m_p}&0&&&\\
&&&0&-I_{m_p}&&&\\
&&e^{-2\pi i k_{p-1}\theta}I_{m_{p-1}}&&&&&\\
&\begin{sideways}$\ddots$\end{sideways}&&&&&&\\
e^{-2\pi i k_1\theta}I_{m_1}&&&&&&&
\end{pmatrix}.
\]
Here $k_1>\cdots>k_{p-1}>0$ are the distinct  partial
indices, $m_j$ is the multiplicity of each of $k_j$ and $-k_j$
for $j<p$, and $2m_p$ is the multiplicity of the zero partial index.
In particular,
$\sum_{j=1}^{p}m_j=k.
$
The signature $(k,k)$ gives the two central blocks the same size
$m_p$. \footnote{Zero determinant winding fixes the sum
of the partial indices, and in principle they need not 
vanish individually.}

For complex $\theta$ we use the holomorphic extension
\[
\tilde\Omega_E(\theta)
=\tilde O_E(\bar\theta)^*S\tilde O_E(\theta).
\]
It is invertible on the full strip: the reflected version
$P_E(\bar\theta)^*SP_E(\theta)$ of \eqref{eiga} is invertible
there, and its off-diagonal zero blocks persist by analyticity.

Let
$$
G_E(\theta)=\begin{cases}
G_E^+(\theta)& \Im\theta\geq 0\\
(i\tilde{\Omega}_E(\theta))^{-1}(G_E^+(\bar{\theta})^*)^{-1}D_E(\theta) &-h<\Im \theta<0
\end{cases}.
$$
The two expressions for $G_E$ agree on $\mathbb T$ by the
factorization identity. They are holomorphic and invertible on
their respective half-strips, so by Schwarz reflection principle, they glue to
$G_E\in C_h^\omega(\mathbb T,GL(2k,\mathbb C))$.

 Let 
\[
P_E(\theta)=
\begin{pmatrix}
I_{m_1}&&&&&&&\\
&\ddots&&&&&&\\
&&I_{m_{p-1}}&&&&&\\
&&&I_{m_p}&&&&\\
&&&&I_{m_p}&&&\\
&&&&&e^{-2\pi i k_{p-1}\theta}I_{m_{p-1}}&&\\
&&&&&&\ddots&\\
&&&&&&&e^{-2\pi i k_1\theta}I_{m_1}
\end{pmatrix}.
\]
This matrix and its inverse are holomorphic on the full strip.
For real $\theta$,
$
P_E(\theta)^*D_E(\theta)P_E(\theta)=D_E(0).
$
It is therefore a constant Hermitian matrix with $k$-positive and $k$-negative eigenvalues. Thus there is a constant matrix $M\in GL(2k,\C)$ such that 
$$
-iM^*P^*_E(\theta)D_E(\theta)P_E(\theta)M=\begin{pmatrix}
0&I_k\\
-I_k&0
\end{pmatrix}.
$$

Finally, we define
\begin{align*}
O_E(\theta)=\tilde{O}_E(\theta)G_E(\theta)P_E(\theta)M
\end{align*}
then one can check that  $O_E\in C^\omega_h(\T,GL_{2d\times 2k}(\C))$ \footnote{$
GL_{2d\times 2k}(\C):=\{F\in M_{2d\times 2k}(\C): {\rm Rank}(F)=2k\}.$.}. Moreover, we have 
$$
O_E(\theta)^*SO_E(\theta)=
\begin{pmatrix}
0&I_k\\
-I_k&0
\end{pmatrix}.
$$
Finally, taking $u_E^j$ to be the $j$-th column of $O_E$, we complete the proof.
\end{pf}

By invariance of $E_c(\theta)$, there exists $M_E\in C_h^\omega(\T,GL(2k,\C))$ be such that
\begin{equation}\label{ff6}
\widehat{A}_E(\theta)O_E(\theta)=O_E(\theta+\alpha)M_E(\theta).
\end{equation}
\begin{Proposition}\label{p1}
For any $E\in\R$ and $\theta\in\T$, we have
$$
M_E(\theta)^*J_{2k}M_E(\theta)=J_{2k}.
$$
\end{Proposition}
\begin{pf}
Taking the adjoint of both sides of \eqref{ff6}, we have
\begin{equation}\label{m}
O_E(\theta)^*\widehat{A}_E(\theta)^*=M_E(\theta)^*O_E(\theta+\alpha)^*.
\end{equation}
Multiplying both sides of the above equation by $S$ , one has
$$
O_E(\theta)^*\widehat{A}_E(\theta)^*S=M_E(\theta)^*O_E(\theta+\alpha)^*S.
$$
Involving the fact that
$$
\widehat{A}_E(\theta)^*S=S\widehat{A}_E(\theta)^{-1},
$$
it follows
$$
O_E(\theta)^*S=M_E(\theta)^*O_E(\theta+\alpha)^*S\widehat{A}_E(\theta).
$$
Multiplying by  $O_E(\theta)$ from the right  side, we obtain
\begin{align*}
&O_E(\theta)^*SO_E(\theta)=M_E(\theta)^*O_E(\theta+\alpha)^*S\widehat{A}_E(\theta)O_E(\theta)\\
=&M_E(\theta)^*O_E(\theta+\alpha)^*SO_E(\theta+\alpha)M_E(\theta).
\end{align*}
Together with \eqref{uveceq}, this completes the proof.
\end{pf}

We define
$$
Sp_{2d\times 2k}(\C)=\{F\in M_{2d\times 2k}(\C): F^*S F=J_{2k}\}.
$$
\begin{Corollary} \label{c4general}For $\alpha\in \R\backslash\Q$,
  $E\in\R$ and $\bar{\omega}(E)=k.$, there exist  $O_E\in C_h^\omega(\T,Sp_{2d\times 2k}(\C))$ and $M_E\in C_h^\omega(\T,Sp(2k,\C))$ such that
\begin{equation}\label{aaa18}
\widehat{A}_E(\theta)O_E(\theta)=O_E(\theta+\alpha)M_E(\theta).
\end{equation}
Moreover if $L(E)>0$, then 
$$
L_1(M_E(\cdot+i\e))=L_1(M_E),\cdots,L_{2k}(M_E(\cdot+i\e))=L_{2k}(M_E),\ \ \forall |\e|<h'<L(E)/2\pi.
$$
\end{Corollary}
\begin{Remark}
$O_E\in C_h^\omega(\T,Sp_{2d\times 2k}(\C))$ means $O_E$ can be
extended to $|\Im \theta|<h$ and for any $\theta$ with $|\Im \theta|<h$, one has $O_E(\bar{\theta})^*SO_E(\theta)=J_{2k}$. 
\end{Remark}
\begin{pf}
Note that by the property of invariant decomposition, Theorem \ref{strip}, Proposition \ref{p1}, and analyticity, if $L(E)>0$, for $|\e|<h'<L(E)/(2\pi)$, we have
\begin{align*}
L_i(M_E(\cdot+i\e))=L_{d-k+i}(\widehat{A}_E(\cdot+i\e))=L_{d-k+i}(\widehat{A}_E),\ \ 1\leq i\leq 2k.
\end{align*}
If also $E\in\Sigma_{v,\alpha}$, these exponents are zero.
\end{pf}
Thus, for  trigonometric polynomial $v$ and Type $k$ energy $E$,  we obtain an
$Sp(2k,\C)$-cocycle of the form $(\alpha,M_E),$ where
$M_E$ are as in Corollary \ref{c4general}, corresponding to the $2k$-dimensional center of $(\alpha,\widehat{A}_E).$ 

\subsection{The general case}
In this subsection, $v$ will be a real analytic function with
$\widehat{H}_{v,\alpha,\theta}$ an infinite-range operator. Therefore
we cannot define a corresponding cocycle, which means the methods in
Section \ref{4.21} are not applicable. We will instead proceed to
define the $Sp(2k,\C)$ cocycle corresponding to the ``center'' of the
long-range operator using the 
trigonometric polynomial approximation.

Assume $v\in C^\omega_{h_1}(\T,\R)$. Throughout this section, we use approximants of
exact degree $l(n)$. Namely, let $n$ be sufficiently large and 
\begin{equation}
v_n(x)=\sum_{m=-l(n)}^{l(n)}\hat v_m e^{2\pi i mx},\ \ l(n)=\max\{|m|\leq n:\hat{v}_m\neq 0\}.
\end{equation}
Every transfer matrix and symplectic form
associated with $v_n$ uses these coefficients, including at the leading
coefficient. For every $0\leq h<h_1$, we have
\begin{align}\label{exact-degree-convergence}
|v_n-v|_h+|v_{n+1}-v_n|_h
&\leq C_s e^{-2\pi(h_1-h)n}.
\end{align}

 To specify
the dependence on $v_n$, in this case, we rewrite
$A_z(x)/\widehat{A}_z(x)$ for $v_n$ as $S_z^{v_n}(x)/\widehat{S}_z^{v_n}(x)$, and denote the corresponding non-negative Lyapunov exponents by $L^{v_n}(z)/\{\gamma^{v_n}_j(z)\}_{j=1}^{l(n)}$ (where $0\leq \gamma_1^{v_n}(z)\leq \gamma^{v_n}_2(z)\leq \cdots \leq \gamma_{l(n)}^{v_n}(z)$), respectively and rewrite $S$ as $S_n$. Denote by $\e_1(z)\geq 0$ the first turning point of $L^v_\e(z)$.

The following multiplicative Jensen's formula is proved  in \cite{gjyz}.
\begin{Theorem}[Theorem 1.2 in \cite{gjyz}]\label{1general}
For $\alpha\in\R\backslash\Q$ and $v\in C^\omega_{h_1}(\T,\R)$, there exist  non-negative $\{\gamma^v_i(E)\}_{i=1}^m$ such that for any $E\in\R$
$$
\gamma^v_i(E)=\lim\limits_{n\rightarrow \infty}\gamma^{v_n}_i(E), \ \ 1\leq i\leq m.
$$
Moreover, \begin{align*}\label{gne1}
L^v_{\e}(E)= L^v_0(E) -\sum_{\{i:\gamma^v_i(E)< 2\pi|\e|\}}  \gamma^v_i(E)+2\pi(\#\{i:\gamma^v_i(E)<2\pi|\e|\})|\e|
\end{align*} for $|\e|<h_1$.
\end{Theorem}
\begin{Remark}
{\rm Note that although Theorem \ref{1general}  is only stated in \cite{gjyz} for real $E$,  the proof holds for all complex $z\in \C$. But due to the lack of evenness, one needs to state the formula for $\e\geq 0$ and $\e\leq 0$ individually. Here for simplicity, we only state half of the formula, i.e., $\e\geq 0$ case.  Namely, there exist non-negative $\{\gamma_i(\bar{z})\}_{i=1}^m$ such that    
$$
\gamma_i(\bar{z})=\lim\limits_{n\rightarrow \infty}\gamma^{n}_i(\bar{z}), \ \ 1\leq i\leq m.
$$
Moreover, \begin{align*}
L_{\e}(z)= L_0(z) -\sum_{\{i:\gamma_i(\bar{z})< 2\pi\e\}}  \gamma_i(\bar{z})+2\pi(\#\{i:\gamma_i(\bar{z})<2\pi\e\})\e
\end{align*} for $0\leq \e<h_1$.}
\end{Remark}
Fix $E_0\in \R$. For $k\ge1$ choose $h_1>\e_1(E_0)$ inside the analytic strip of $v$ ($k=0$ gives the zero center). For functions of $z,\theta$, we write
\[
|F|_{r,h}
:=\sup_{\substack{|z-E_0|<r\\
                  \theta\in\mathbb C/\mathbb Z,\ |\Im\theta|<h}}
\|F(z,\theta)\|.
\]
Denote by $\e_2(E)> 0$ the second non-negative turning point of $L^v_\e(E)$ in $[0,h_1)$. If there is only one non-negative turning point there, set $\e_2(E)=h_1.$ For any $h<h_0<\frac{\e_2(E_0)-\e_1(E_0)}{200\pi}+L(E_0)/2\pi$, we have

\begin{Theorem} \label{theorem-main1general}
For $\alpha\in \R\backslash\Q$, $v\in C_{h_1}^\omega(\T,\R),$ and $E_0\in\R$ with $\bar{\omega}(E_0)=k$, there exist a neighborhood  $I$ of $E_0$ and, for all large $n$, $O_n\in C^\omega(I\times \T,Sp_{2l(n)\times 2k}(\C))$ and $M_n\in C^\omega(I\times \T,Sp(2k,\C))$ such that  for any $E\in I$, we have
$$
\widehat{S}_E^{v_n}(\theta)O_n(E,\theta)=O_n(E,\theta+\alpha)M_n(E,\theta).
$$
Moreover, there is $M\in C^\omega(I\times \T,Sp(2k,\C))$  such that
$$
|M_n-M|_{r,h}\leq Ce^{-cn}
$$
for some $C,c,r>0$. If $L(E_0)>0$, then for any $0<h'<L(E_0)/2\pi$, there is $E_0\in I_{h'}\subset I$,  such that for any $E\in I_{h'}$,
$$
L_1(M(E,\cdot+i\e))=L_1(M(E,\cdot)),\cdots, L_{k}(M(E,\cdot+i\e))=L_k(M(E,\cdot)),\ \ \forall |\e|<h'.
$$
Also, denoting
$$
O_n(E,\theta)=\begin{pmatrix}O_n(E,\theta,l(n)-1)\\ O_{n}(E,\theta,l(n)-2)\\ \vdots\\ O_n(E,\theta,-l(n))\end{pmatrix}
$$ 
 there is $O(\cdot,\cdot,0)\in C^\omega(I\times \T,\C^{2k})$ such that 
$$
|O_n(0)-O(0)|_{r,h}\leq Ce^{-cn}.
$$ 
\end{Theorem}
\begin{pf}
We first work with $L(E_0)/(2\pi)<h<h_0$
Let $\e=\frac{\e_1(E_0)+\e_2(E_0)}{2}$. At $E_0$ both shifted cocycles are regular with $L^v_{\pm\e}(E_0)>2\pi h_0$, hence uniformly hyperbolic by  \cite{avila0}. This persists for nearby complex $z$ by openness and continuity. Thus there exists a continuous invariant
splitting $\C^{2}=E_\e^s(x)\oplus E_\e^u(x)$, such that  any nonzero
$u_\e(x)\in E_\e^u(x)$ and nonzero $s_\e(x)\in E_\e^s(x)$, setting
$$
\begin{pmatrix}
f^m_1(x)\\
f^{m-1}_1(x)
\end{pmatrix}=(S_z^{v})_m(x+i\e)\cdot s_\e(x),
$$
$$
\begin{pmatrix}
f^m_2(x)\\
 f^{m-1}_2(x)
\end{pmatrix}=(S_z^{v})_m(x+i\e)\cdot u_\e(x).
$$
For any $x\in\T$, we have
\begin{equation}\label{new13}
\limsup\limits_{m\rightarrow\infty}\frac{\ln(|f^m_1(x)|^2+|f^{m-1}_1(x)|^2)}{2m}=-L_\e^{v}(z),
\end{equation}
\begin{equation}\label{new14}
\limsup\limits_{m\rightarrow\infty}\frac{\ln(|f^m_2(x)|^2+|f^{m-1}_2(x)|^2)}{2m}=L_\e^{v}(z).
\end{equation}

After shrinking the complex neighborhood of $E_0$ if necessary, we
have $L^v_{\pm\e}(z)>2\pi h_0$ for both signs. 
Thus, for every $|\epsilon|\leq h_0$, the cocycle
$(\alpha,e^{2\pi\epsilon}D^{-1}S_z^v(\cdot\pm i\e)D)$, where
$D={\rm diag}\{1,e^{2\pi\epsilon}\}$, is uniformly hyperbolic. Indeed, the scalar factor shifts
the two exponents to $2\pi\epsilon\pm L^v_{\sigma\e}(z)$ for
each  $\sigma\in\{+1,-1\}$, and these have opposite signs.
By the complex-valued version of
Johnson's theorem \cite{gjyz}
the 
resolvent exists for every $x\in\mathbb{T}$. Uniform growth on the invariant lines gives
\[
\big|\langle\delta_p,(H_{v(\cdot\pm i\varepsilon),x,\alpha}-z)^{-1}\delta_q\rangle\big|\le Ce^{-c|p-q|},
\]
with fixed $2\pi h_0<c<\min_\sigma L^v_{\sigma\varepsilon}(z)$ on a smaller neighborhood of $E_0$. For $|\epsilon|\le h_0$, $e^{2\pi\epsilon(p-q)}\langle\delta_p,(H_{v(\cdot\pm i\varepsilon),x,\alpha}-z)^{-1}\delta_q\rangle$ is summable in rows and columns and satisfies both inverse equations for $H^\epsilon_{v(\cdot\pm i\varepsilon),x,\alpha}-z$.

Let $v\in C^\omega_{h_1}(\T,\R)$ \footnote{By our assumption and continuity of $\e_1(E)$,
  $h_1>\e_1(E)$.} and $\delta$ be any fixed constant such that $0<\delta<\frac{\e_2(E_0)-\e_1(E_0)}{10^8\pi}$. By the resolvent identity and compactness argument,  there is an analytic neighborhood $U_v$ of $v$ and a complex neighborhood $D_{E_0}$ of $E_0$, such that for any $v'\in U_v$ and $z\in D_{E_0}$, we have
\begin{equation}\label{gjne1}
|(H^\epsilon_{v(\cdot\pm i\e),x,\alpha}-z)^{-1}-(H^\epsilon_{v'(\cdot\pm i\e),x,\alpha}-z)^{-1}|\leq C(E_0,v,\delta)|v-v'|_{h_1-\delta}.
\end{equation} 
for $(x,\epsilon)\in\T\times[-h_0,h_0]$ \footnote{Note $\e$ is fixed at the beginning.}.

Fix $\ell\in\Z$ and let $\delta_\ell(x,m):=\delta_\ell(m)$ for any $x\in\T$
and $m\in\Z.$ Then, by \eqref{definitionU},
$$
(U^{-1}\delta_\ell)(x,m)=e^{-2\pi i \ell x}\delta_0(x,m).
$$
\begin{equation}\label{keyeq}
\begin{aligned}
\left(U(H^{-\epsilon}_{v'(\cdot\pm i\e),\alpha}-z)^{-1}U^{-1}\delta_\ell\right)(\theta,m)=\left((\widehat{H}^{\epsilon}_{v'(\cdot\pm i\e),\alpha}-z)^{-1}\delta_\ell\right)(\theta,m).
\end{aligned}  \footnote{This is the key equality which implies the convergence of the center.}
\end{equation}
We denote
\begin{align}\label{zx1}
f^\epsilon_{v'(\cdot\pm i\e)}(z,x,m)&=\left((H^\epsilon_{v'(\cdot\pm i\e),\alpha}-z)^{-1}U^{-1}\delta_\ell\right)(x,m)\\ \nonumber
&=\langle \delta_m,e^{-2\pi i\ell x}(H^\epsilon_{v'(\cdot\pm i\e),\alpha,x}-z)^{-1}\delta_0\rangle,
\end{align}
\begin{equation}\label{final3}
g_{v'(\cdot\pm i\e)}(z,\theta+i\epsilon,m)=\left(U(H^{-\epsilon}_{v'(\cdot\pm i\e),\alpha}-z)^{-1}U^{-1}\delta_\ell\right)(\theta,m).
\end{equation}
By \eqref{definitionU}, we have
$$
g_{v'(\cdot\pm i\e)}(z,\theta+i\epsilon,m)=\sum_{p\in\Z}\int_{\T}e^{2\pi im\eta}e^{2\pi ip(m\alpha+\theta)}f^{-\epsilon}_{v'(\cdot\pm i\e)}(z,\eta,p)d\eta,
$$
By \eqref{gjne1} at $-\epsilon$, \eqref{keyeq}, unitarity and norm continuity of the dual fibers,
\begin{equation}\label{new15}
|g_{v'(\cdot\pm i\e)}(z,\theta+i\epsilon,m)-g_{v(\cdot\pm i\e)}(z,\theta+i\epsilon,m)|\leq C(E_0,v,\delta)e^{2\pi|m-\ell|\e}e^{\pm2\pi  (m-\ell)\e}|v'-v|_{h_1-\delta}.
\end{equation}

For any $\theta\in\T$ and $|\epsilon|\leq h_0$, we define
\begin{equation}\label{lana4}
u^\ell_n(z,\theta+i\epsilon,m):=e^{-2\pi (m-\ell)\e}g_{v_n(\cdot+i\e)}(z,\theta+i\epsilon,m)-e^{2\pi (m-\ell)\e}g_{v_n(\cdot-i\e)}(z,\theta+i\epsilon,m).
\end{equation}
By \eqref{final3}, $e^{\mp2\pi (m-\ell)\e}g_{v_n(\cdot\pm i\e)}(z,\theta+i\epsilon,m)$ are solutions of $(\widehat{H}^\epsilon_{v_n,\alpha,\theta}-z)u=\delta_\ell$, hence $\{u^\ell_n(z,\theta+i\epsilon,m)\}_{m\in\Z}$ are solutions of $(\widehat{H}^\epsilon_{v_n,\alpha,\theta}-z)u=0$. Moreover, by \eqref{new15}, we have $u^\ell_n(m)\in C^\omega(I\times\T,\C)$,
\begin{equation}\label{new200}
|u^\ell_n(m)-u^\ell_{n+1}(m)|_{r,h_0}\leq Ce^{2\pi|m-\ell|\e}e^{-2\pi \delta n}
\end{equation}
for some $r>0$ and all $m,\ell\in\mathbb Z$, by unitary duality and norm continuity of the fibers.

We next need the following lemma,
\begin{Lemma}[\cite{gjyz}]\label{basic}
Consider the following $2d$ order difference operator,
$$
(Lu)(n)=\sum\limits_{k=-d}^da_ku(n+k)+V(n)u(n).
$$
If the eigenequation $Lu=Eu$ has $2d$ linearly independent solutions $\{\phi_i\}_{i=1}^{2d}$ satisfying
$$
\phi_i\in\ell^2(\Z^+)(i=1,\cdots,m),\ \ \phi_i\in\ell^2(\Z^-)(i=m+1,\cdots,2d),
$$
and $L-EI$ is invertible. Then we have
$$
\langle\delta_p,(L-EI)^{-1}\delta_q\rangle=\begin{cases}
\frac{\sum\limits_{i=1}^m\phi_i(p)\Phi_{1,i}(q)}{a_d\det{\Phi(q)}} &\text{$p\geq q+1$}\\
-\frac{\sum\limits_{i=m+1}^{2d}\phi_i(p)\Phi_{1,i}(q)}{a_d\det{\Phi(q)}} &\text{$p\leq q$}
\end{cases},
$$
where
$$
\Phi(q)=\begin{pmatrix}\phi_1(q+d)&\phi_2(q+d)&\cdots&\phi_{2d}(q+d)\\ \phi_1(q+d-1)&\phi_2(q+d-1)&\cdots&\phi_{2d}(q+d-1)\\ \vdots&\vdots& &\vdots\\
\phi_1(q-d+1)&\phi_2(q-d+1)&\cdots&\phi_{2d}(q-d+1)\end{pmatrix}
$$
and $\Phi_{i,j}(q)$ is the $(i,j)$-th cofactor of $\Phi(q)$.
\end{Lemma}
\begin{Lemma}\label{add1}
For any $-n\leq \ell\leq n-1$, $E\in\R$ with $|E-E_0|$ sufficiently small, and  $|\Im \theta|<h_0$,
\begin{itemize}
\item we have
$$
u^\ell_n(E,\theta)=\begin{pmatrix}u^\ell_n(E,\theta,l(n)-1)\\ u^\ell_n(E,\theta,l(n)-2)\\ \vdots \\u^\ell_n(E,\theta,-l(n))\end{pmatrix}\in E^n_c(\theta)
$$
where $E^n_c(\theta)$ is the center of the continuous decomposition of $\C^{2l(n)}$ corresponding to $(\alpha, \widehat{S}_E^{v_n}(\theta))$;
\item $G_n(E,\theta):=\begin{pmatrix} -u_n^{l(n)-1}(E,\theta)&\cdots &-u_{n}^{-l(n)}(E,\theta)\end{pmatrix}$ satisfies $G_n(E,\bar\theta)^*=-G_n(E,\theta)$ and ${\rm Rank} (G_n(E,\theta))=2k$.
\end{itemize}
\end{Lemma}
\begin{pf}
Note that $\bar{\omega}(E_0)=k$, thus by Theorem \ref{1general}, we have
\begin{equation}\label{zx2}
\liminf\limits_{n\rightarrow\infty} \gamma_{k+1}^{v_n}(E_0)\geq 2\pi \e_2(E_0)>\lim\limits_{n\rightarrow\infty} \gamma_{k}^{v_n}(E_0)=2\pi \e_1(E_0).
\end{equation}
Set $a=2\pi\e_1(E_0)$, $b=2\pi\e_2(E_0)$ and $\eta=(\e_2(E_0)-\e_1(E_0))/75$.
For all large $n$ and real $E$ close to $E_0$, the finite Jensen formula gives
\[
 \gamma_k^{v_n}(E)<a+\eta/8,\qquad \gamma_{k+1}^{v_n}(E)>b-\eta/8 .
\]
Indeed, on the interior of the affine segment $(a/(2\pi),b/(2\pi))$ the slopes of $L^{v_n}_s(E)$ on any fixed compact subinterval equal $2\pi k$, by continuity, convexity and quantization.

Write  $t=L^{v_n}(E)/(2\pi)$, and let $L_j(s)$ denote the ordered exponents of $\widehat S_E^{v_n}(\cdot+is)$. By Theorem~\ref{strip} and its proof,
\[
  0\le L_j(s)-L_j(0)\le 2\pi(|s|-t)_+,\qquad 1\le j\le l(n).
\]
For $|s|<h_0$ the right side is less than $7\eta/8$ by continuity and the choice of $h_0$. Thus $L_{l(n)-k}(s)>b-\eta>a+\eta>L_{l(n)-k+1}(s)$. Integrality at this gap and \eqref{less than 1} give $(l(n)-k)$-regularity, hence $(l(n)-k)$-domination by Theorem~\ref{t2.1}. Reflection gives $(l(n)+k)$-domination on the other side of the center.

Thus there exist a holomorphic invariant decomposition of $\C^{2l(n)}$,
$$
\C^{2l(n)}=E^n_s(\theta)\oplus E^n_c(\theta)\oplus E^n_u(\theta),
$$
which by Theorem \ref{trivial} and Lemma \ref{uvec} implies that there are linearly independent $\{f^n_j(E,\theta)\}_{j=1}^{l(n)-k}\in E_s^n(\theta)$, $\{v_j^n(E,\theta)\}_{j=1}^{2k}\in E^n_c(\theta)$ and $\{g_j^n(E,\theta)\}_{j=1}^{l(n)-k}\in E_u^n(\theta)$ depending analytically on $\theta$ on $\Omega_{h_0}$, such that
\begin{equation}\label{lana2}
\widetilde{O}_n(E,\theta)^*S_n\widetilde{O}_n(E,\theta)=J_{2k}.
\end{equation}
for any $\theta\in\T$ where $\widetilde{O}_n(E,\theta)=\begin{pmatrix}v_1^n(E,\theta)&v_2^n(E,\theta)&\cdots&v_{2k}^{n}(E,\theta)\end{pmatrix}$.

Let $\{f_j^n(E,\theta,\ell)\}_{j=1}^{l(n)-k}$, $\{g_j^n(E,\theta,\ell)\}_{j=1}^{l(n)-k}$, $\{v_j^n(E,\theta,\ell)\}_{j=1}^{2k}$ be  $2l(n)$ linearly independent solutions of $\widehat{H}_{v_n,\alpha,\theta}u=Eu$ with the above corresponding initial datum. The preceding bounds, reflection at $-\Im\theta$, and the uniform subadditive upper bound
give for $|\Im \theta|<h_0$,
$$
\limsup\limits_{m\rightarrow \infty}\frac{1}{2m}\ln\sum\limits_{\ell=-l(n)}^{l(n)-1} |f_j^n(E,\theta,m+\ell)|^2<-2\pi \e_2(E_0)+\eta,\ \  1 \leq j\leq l(n)-k,
$$
$$
\limsup\limits_{m\rightarrow \infty}\frac{1}{2m}\ln\sum\limits_{\ell=-l(n)}^{l(n)-1} |g_j^n(E,\theta,-m+\ell)|^2<-2\pi \e_2(E_0)+\eta,\ \  1 \leq j\leq l(n)-k,
$$
$$
\limsup\limits_{m\rightarrow \infty}\frac{1}{2m}\ln\sum\limits_{\ell=-l(n)}^{l(n)-1} |v_j^n(E,\theta,\pm m+\ell)|^2\leq 2\pi \e_1(E_0)+\eta,\ \ 1\leq j\leq 2k.
$$
Since 
$$
\int v_n(\theta+i\e)e^{2\pi ik\theta} d\theta=e^{2\pi k\e}\int v_n(\theta)e^{2\pi ik\theta} d\theta,
$$ 
we have that $\{e^{\pm 2\pi \ell\e}f_j^n(E,\theta,\ell)\}_{j=1}^{l(n)-k}$, $\{e^{\pm 2\pi \ell\e}g_j^n(E,\theta,\ell)\}_{j=1}^{l(n)-k}$, $\{e^{\pm 2\pi \ell\e}v_j^n(E,\theta,\ell)\}_{j=1}^{2k}$ are $2l(n)$ independent solutions of $\widehat{H}_{v_n(\cdot\pm i\e),\alpha,\theta}u=Eu$. Thus we have
\begin{align}\label{so1}
e^{\pm 2\pi \ell\e}f_j^n(E,\theta,\ell)\in\ell^2(\Z^+), \ \ 1\leq j\leq l(n)-k,
\end{align}
\begin{align}\label{so2}
e^{\pm 2\pi \ell\e}g_j^n(E,\theta,\ell)\in\ell^2(\Z^-), \ \ 1\leq j\leq l(n)-k,
\end{align}
\begin{align}\label{so3}
e^{\pm 2\pi \ell\e}v_j^n(E,\theta,\ell) \in\ell^2(\Z^\mp),\ \ 1\leq j\leq 2k.
\end{align}

By \eqref{so1}-\eqref{so3} and Lemma \ref{basic}, we have for $\theta\in \Omega_{h_0}$, 
\begin{align}\label{g44}
&g_{v_n(\cdot+ i\e)}(E,\theta,m)=\langle \delta_m, (\widehat{H}_{v_n(\cdot+ i\e),\alpha,\theta}-E)^{-1}\delta_\ell\rangle\\ \nonumber
&=\begin{cases}
\frac{e^{2\pi (m-\ell)\e}\left(\sum\limits_{j=1}^{l(n)-k}f^n_j(E,\theta,m)\Phi^n_{1,j}(E,\theta,\ell)\right)}{\hat{v}_{l(n)}\det{\Phi^n(E,\theta,\ell)}} &\text{$m\geq \ell+1$}\\
-\frac{e^{2\pi (m-\ell)\e}\left(\sum\limits_{j=1}^{2k}v_j^n(E,\theta,m)\Phi^n_{1,l(n)-k+j}(E,\theta,\ell)+\sum\limits_{j=1}^{l(n)-k}g^n_j(E,\theta,m)\Phi^n_{1,l(n)+k+j}(E,\theta,\ell)\right)}{\hat{v}_{l(n)}\det{\Phi^n(E,\theta,\ell)}} &\text{$m\leq \ell$}
\end{cases},
\end{align}
\begin{align}\label{g45}
&g_{v_n(\cdot-i\e)}(E,\theta,m)=\langle \delta_m, (\widehat{H}_{v_n(\cdot-i\e),\alpha,\theta}-E)^{-1}\delta_\ell\rangle\\ \nonumber
&=\begin{cases}
\frac{e^{-2\pi (m-\ell)\e}\left(\sum\limits_{j=1}^{2k}v_j^n(E,\theta,m)\Phi^n_{1,l(n)-k+j}(E,\theta,\ell)+\sum\limits_{j=1}^{l(n)-k}f^n_j(E,\theta,m)\Phi^n_{1,j}(E,\theta,\ell)\right)}{\hat{v}_{l(n)}\det{\Phi^n(E,\theta,\ell)}} &\text{$m\geq \ell+1$}\\
-\frac{\e^{-2\pi (m-\ell)\e}\left(\sum\limits_{j=1}^{l(n)-k}g^n_j(E,\theta,m)\Phi^n_{1,l(n)+k+j}(E,\theta,\ell)\right)}{\hat{v}_{l(n)}\det{\Phi^n(E,\theta,\ell)}} &\text{$m\leq \ell$}
\end{cases},
\end{align}
where
{\tiny $$
\Phi^n(E,\theta,\ell):=\begin{pmatrix}f^n_1(E,\theta,\ell+l(n))&\cdots&v_1^n(E,\theta,\ell+l(n))&\cdots&v_{2k}^n(E,\theta,\ell+l(n))&\cdots g_{l(n)-k}^n(E,\theta,\ell+l(n))\\ f^n_1(E,\theta,\ell+l(n)-1)&\cdots&v_1^n(E,\theta,\ell+l(n)-1)&\cdots&v_{2k}^n(E,\theta,\ell+l(n)-1)&\cdots g_{l(n)-k}^n(E,\theta,\ell+l(n)-1)\\ \vdots& &\vdots &&\vdots &\vdots\\
f^n_1(E,\theta,\ell-l(n)+1)&\cdots&v_1^n(E,\theta,\ell-l(n)+1)&\cdots&v_{2k}^n(E,\theta,\ell-l(n)+1)&\cdots g_{l(n)-k}^n(E,\theta,\ell-l(n)+1)\end{pmatrix}.
$$}

Let 
$$
F_n(E,\theta)=\begin{pmatrix}f^n_1(E,\theta)&f^n_2(E,\theta)&\cdots&f_{l(n)-k}^n(E,\theta)\end{pmatrix},
$$
$$
K_n(E,\theta)=\begin{pmatrix}g^n_1(E,\theta)&g^n_2(E,\theta)&\cdots&g_{l(n)-k}^n(E,\theta)\end{pmatrix}.
$$
By symplectic orthogonality, we have
$$
F_n(E,\bar{\theta})^*S_nF_n(E,\theta)=K_n(E,\bar{\theta})^*S_nK_n(E,\theta)=0,
$$
$$
F_n(E,\bar{\theta})^*S_n\widetilde{O}_n(E,\theta)=K_n(E,\bar{\theta})^*S_n\widetilde{O}_n(E,\theta)=0.
$$
On the other hand, by \eqref{lana2} and symplectic invariance, we further have
\begin{align}\label{fur}
\nonumber &\Phi^n(E,\bar{\theta},\ell)^*S_n\Phi^n(E,\theta,\ell)=\Phi^n(E,\bar{\theta},-1)^*S_n\Phi^n(E,\theta,-1)\\
=&\begin{pmatrix} &&A_n(E,\theta,\ell)\\ & J_{2k}&\\  -A_n(E,\bar{\theta},\ell)^*&& \end{pmatrix}
\end{align}
for some $A_n(E,\theta,\ell)$.

We need the following lemma
\begin{Lemma}\label{lana3}
Assume $T_1=\begin{pmatrix}
t_l& \cdots &t_1\\
&\ddots&\vdots\\
&&t_l
\end{pmatrix}$ and $T=\begin{pmatrix}& -T^*_1\\ T_1&\end{pmatrix}$, $A=(a_{i,j})$ is an $2l\times2l$ matrix such that 
$$
A^*TA=\begin{pmatrix}
&&A_1\\
&J_{2l-2l_1}&\\
-A^*_1&&
\end{pmatrix}
$$
for some invertible $A_1$ and $J_{2l-2l_1}=\begin{pmatrix}&I_{l-l_1}\\ -I_{l-l_1}&\end{pmatrix}$. We have that
$$
\begin{pmatrix}
A_{1,1}\\
A_{1,2}\\
\vdots\\
A_{1,2l}
\end{pmatrix}={t_l}{\det{A}}\begin{pmatrix} -(A_1^*)^{-1}\alpha^*_3\\ -J_{2l-2l_1}\alpha^*_2
\\ (A_1)^{-1}\alpha^*_1\end{pmatrix}
$$
where $A_{ij}$ are the $(i,j)$-th cofactor of $A$ and
$$
\alpha_1=\begin{pmatrix}a_{l+1,1}&a_{l+1,2}&\cdots &a_{l+1.l_1}\end{pmatrix}\ \ \alpha_2=\begin{pmatrix}a_{l+1,l_1+1}&a_{l+1,l_1+2}&\cdots &a_{l+1.2l-l_1}\end{pmatrix},$$
$$
\alpha_3=\begin{pmatrix}a_{l+1,2l-l_1+1}&a_{l+1,2l-l_1+2}&\cdots &a_{l+1.2l}\end{pmatrix}.
$$
\end{Lemma}
\begin{pf}
By the assumption, we have that
$$
A^{-1}=\begin{pmatrix}
&&(-A^*_1)^{-1}\\
&-J_{2l-2l_1}&\\
(A_1)^{-1}&&
\end{pmatrix}
A^*T,
$$
thus 
$$
\begin{pmatrix}
A_{1,1}\\
A_{1,2}\\
\vdots\\
A_{1,2l}
\end{pmatrix}=\det{A}\begin{pmatrix}
&&(-A^*_1)^{-1}\\
&-J_{2l-2l_1}&\\
(A_1)^{-1}&&
\end{pmatrix}
A^*T\delta_1.
$$
We have that $A^*T\delta_1$ is $t_l A^*\delta_{l+1},$ so
$$
A^*T\delta_1=t_l\begin{pmatrix}
\alpha_1^*\\ \alpha_2^*\\ \alpha_3^*
\end{pmatrix}
$$
It follows that
$$
\begin{pmatrix}
A_{1,1}\\
A_{1,2}\\
\vdots\\
A_{1,2l}
\end{pmatrix}={t_l}{\det{A}}\begin{pmatrix} -(A_1^*)^{-1}\alpha^*_3\\ -J_{2l-2l_1}\alpha^*_2
\\ (A_1)^{-1}\alpha^*_1\end{pmatrix}
$$
\end{pf}

By \eqref{S}, \eqref{fur} and  Lemma \ref{lana3} with $l=l(n)$, $l_1=l(n)-k$ \footnote{So that only the middle term matters.}, we have that
\begin{align}\label{g42}
\Phi_{1,l(n)-k+j}^n(E,\theta,\ell)=-\hat{v}_{l(n)}\det{\Phi^n(E,\theta,\ell)}\overline{v_{k+j}^n(E,\bar{\theta},\ell)}, \ \ 1\leq j\leq k,
\end{align}
\begin{align}\label{g43}
\Phi_{1,l(n)-k+j}^n(\theta,\ell)=\hat{v}_{l(n)}\det{\Phi^n(E,\theta,\ell)}\overline{v_{j-k}^n(E,\bar{\theta},\ell)},\ \ k+1\leq j\leq 2k.
\end{align}
By \eqref{lana4}, \eqref{g44}, \eqref{g45}, \eqref{g42} and \eqref{g43}, we have
\begin{align}\label{final10}
u^\ell_n(E,\theta,m)&=e^{-2\pi (m-\ell)\e}g_{v_n(\cdot+i\e)}(E,\theta,m)-e^{2\pi (m-\ell)\e}g_{v_n(\cdot-i\e)}(E,\theta,m)\\ \nonumber
&=\sum_{j=1}^{k}v_j^n(E,\theta,m)\overline{v_{k+j}^n(E,\bar{\theta},\ell)}-\sum_{j=k+1}^{2k}v_j^n(E,\theta,m)\overline{v_{j-k}^n(E,\bar{\theta},\ell)}.
\end{align}
Thus $u^\ell_n(E,\theta)\in E^n_c(\theta)$.

We denote
$$
H_n(E,\theta)=\begin{pmatrix}
v_1^n(E,\theta,l(n)-1)&v_2^n(E,\theta,l(n)-1)&\cdots&v_{2k}^n(E,\theta,l(n)-1)\\
v_1^n(E,\theta,l(n)-2)&v_2^n(E,\theta,l(n)-2)&\cdots& v_{2k}^n(E,\theta,l(n)-2)\\
\vdots&\vdots&&\vdots\\
v_1^n(E,\theta,-l(n))&v_2^n(E,\theta,-l(n))&\cdots&v_{2k}^n(E,\theta,-l(n))
\end{pmatrix}.
$$
By \eqref{final10} and the
definition,
\begin{equation}\label{lanag}G_n(E,\theta)=-H_n(E,\theta)J_{2k}H_n(E,\bar{\theta})^*.\end{equation}
It follows that $G_n(E,\theta)$ is skew-Hermitian for $\theta\in\T$ and ${\rm Rank}(G_n(E,\theta))=2k$ for any $\theta$ with $|\Im \theta|<h_0$.
\end{pf}
Let
\[
D_n:=\operatorname{diag}
\bigl(e^{-2\pi(\varepsilon+\delta)|j|}\bigr)_
{j=l(n)-1,l(n)-2,\ldots,-l(n)},\ \ 
\widetilde G_n(E,\theta):=D_nG_n(E,\theta)D_n^*.
\]
Here $\varepsilon=(\varepsilon_1(E_0)+\varepsilon_2(E_0))/2$ is
the fixed parameter chosen at the beginning of the proof.

\begin{Lemma}\label{eigen0}
For $G_{n}(E,\theta)$ and $\widetilde{G}_n(E,\theta)$  as above and any $\theta\in\T$, $iG_n(E,\theta)$ and $i\widetilde{G}_n(E,\theta)$ have $k$ positive eigenvalues  $\mu^{n}_j(E,\theta)$ and $\tilde{\mu}^{n}_j(E,\theta)$, and $k$ negative eigenvalues $\kappa^n_j(E,\theta)$ and $\tilde{\kappa}^n_j(E,\theta)$, $1\leq j\leq k$ with 
$$
|\mu^{n}_j(E,\theta)|, |\kappa_j^n(E,\theta)|\geq c_0>0
$$
$$
|\tilde{\mu}^{n}_j(E,\theta)|, |\tilde{\kappa}_j^n(E,\theta)|\geq \tilde{c}_0>0
$$
where $c_0=\frac{1}{2\sup_{x\in \T}|v(x)|}$ and $\tilde{c}_0=\frac{1}{2|v(x)|_{h_1}}$.
\end{Lemma}
\begin{pf}
We first consider $G_n(E,\theta)$. Note that for any $\theta\in\T$,  $H^*_n(E,\theta)H_{n}(E,\theta)$ is
positive definite and $J_{2k}$ is skew-Hermitian.  There  is an
invertible matrix $U_n(E,\theta)$ and {$k$} positive values $g_j(E,\theta)$
and {$k$} negative values $h_j(E,\theta)$ \footnote{We do not care about
  the regularity of $U_n$, $g_j$ and $h_j$.  $g_j$ and $h_j$ also
  depend on $n$ which we omit from notation for convenience.}  such that 
\begin{equation}\label{lanah}
U^*(E,\theta)H^*_n(E,\theta)H_{n}(E,\theta)U(E,\theta)=I,
\end{equation}
\begin{equation}\label{lanah'}
U(E,\theta)^*J_{2k}U(E,\theta)={\begin{pmatrix}
ig_1(E,\theta)\\ & \ddots\\ && ig_k(E,\theta)\\ & &&ih_1(E,\theta)\\&&&&\ddots\\&&&&&ih_k(E,\theta)
\end{pmatrix}}.
\end{equation}

We have by \eqref{lana2} that
$H_n^*(E,\theta)S_nH_n(E,\theta)=J_{2k}$. Therefore, \eqref{lanah} and \eqref{lanah'} implies
$ig_j(E,\theta)=w_j^n(E,\theta)^*S_nw_j^n(E,\theta)$ and
$ih_j(E,\theta)=w_{k+j}^n(E,\theta)^*S_nw_{k+j}^n(E,\theta)$ for $1\leq j\leq k$ where $w_j(E,\theta)$ is the $j$-th column of $H_n(E,\theta)U_n(E,\theta),$ and $\|w_j(E.\theta)\|=1$ for $1\leq j\leq 2k$.
By the definition of $S_n$, it follows that

\[
|g_j(E,\theta)|^{-1},\ |h_j(E,\theta)|^{-1}
\geq\frac{1}{\|S_n\|}
\geq\frac{1}{2\sup_{x\in\T}|v(x)|},\qquad 1\leq j\leq k.
\]
\par

Indeed $C_n$ is a block of the Laurent operator with symbol $v_n$, so $\|S_n\|=\|C_n\|=\|v_n\|_0\le2\|v\|_0$ for large $n$. For $\widetilde G_n$ replace $H_n,S_n$ by $D_nH_n,D_n^{-1}S_nD_n^{-1}$. The weighted block has symbol $v_n(\cdot-i(\e+\delta))$, hence norm at most $2\|v\|_{h_1}$ since $\e+\delta<h_1$.
Finally, by \eqref{lanag}, we have 
$$
G_n(E,\theta)=H_n(E,\theta)U_n(E,\theta)\left(U_n(E,\theta)^*J_{2k}U_n(E,\theta)\right)^{-1}U_n(E,\theta)^*H_n(E,\theta)^*
$$
since $J_{2k}^{-1}=-J_{2k}$.

Note that by \eqref{lanah} and \eqref{lanah'}
\begin{align*}
G'_{n}(E,\theta)&:=\left(U_n(E,\theta)^*J_{2k}U_n(E,\theta)\right)^{-1} \\
&=U_n(E,\theta)^*H_n(E,\theta)^*H_{n}(E,\theta)U_n(E,\theta)\left(U_n(E,\theta)^*J_{2k}U_n(E,\theta)\right)^{-1}
\end{align*} 
so  $G'_{n}(E,\theta)$  and $G_n(E,\theta)$ have the same non-zero
eigenvalues, and the result follows. \end{pf}

For any $m\leq l(n)$, we  denote 
$$
H^{m}_n(E,\theta)=\begin{pmatrix}
v_1^n(E,\theta,m-1)&v_2^n(E,\theta,m-1)&\cdots&v_{2k}^n(E,\theta,m-1)\\
v_1^n(E,\theta,m-2)&v_2^n(E,\theta,m-2)&\cdots& v_{2k}^n(E,\theta,m-2)\\
\vdots&\vdots&&\vdots\\
v_1^n(E,\theta,-m)&v_2^n(E,\theta,-m)&\cdots&v_{2k}^n(E,\theta,-m)
\end{pmatrix},$$
the matrix consisting of the middle $2m$ rows of $H_n.$ Let
$G_n^{m}(E,\theta):=-H_n^{m}(E,\theta)J_{2k}H_{n}^{m}(E,\bar{\theta})^*$,
$L_n^{m}(E,\theta):=-H_{n}(E,\theta)J_{2k}H_{n}^{m}(E,\bar{\theta})^*$. It
is easily seen that $L^m_n(E,\theta)$ is a $2l(n)$ by $2m$ matrix consisting of $2m$
middle columns of $G_n(E,\theta),$ so $L^m_n(E,\theta)=\begin{pmatrix}
  -u_n^{m-1}(E,\theta)&\cdots &-u_{n}^{-m}(E,\theta)\end{pmatrix},$ and
$G_{n}^m(E,\theta)$ is a $2m$ by $2m$ matrix consisting $2m$ middle rows
of $L_n^m(E,\theta)$. In particular, $G_n^m(E,\theta)$ is a submatrix of
$L_n^m(E,\theta),$ and $L_n^m(E,\theta)$ is a submatrix of
$G_n(E,\theta)$. Let $L_n^m(E,\theta,j),j=-l(n),\ldots,l(n)-1,$ denote the $j$-th
row of $L_n^m(E,\theta)$ and
\begin{equation*}
D={\rm diag}\{e^{-2\pi(m-1)(\e+\delta)},\cdots, e^{-2\pi (\e+\delta)},\cdots, e^{-2\pi m(\e+\delta)}\},\ \ \widetilde{G}^m_n(E,\theta)=DG^m_n(E,\theta)D^*.
\end{equation*}

The following lemma is  key to the convergence of the center
\begin{Lemma}\label{add2}
There is  $n_0$ sufficiently large with $l(n_0)=n_0$, such that for all sufficiently large $n$ we have that
\begin{itemize}
\item $|G_{n}^{n_0}-G_{n+1}^{n_0}|_{r,h},
  |L_{n}^{n_0}(j)-L_{n+1}^{n_0}(j)|_{r,h} \leq
  C_{n_0}e^{-\delta n}, j=-n_0,\ldots,n_0-1$;
{\item there are vectors $e_j\in C^\omega(I\times \T,\C^{2n_0})$, $1\leq j\leq 2k$, such that for any $|z-E_0|<r$ and $|\Im \theta|<h$, we have
\begin{equation}\label{gn1'}
{\rm Rank}\left(L_{n}^{n_0}(z,\theta)(e_1(z,\theta),e_2(z,\theta)\cdots, e_{2k}(z,\theta))\right)=2k,\ \ L_{n}^{n_0}(z,\theta)e_j(z,\theta)\in E_c^n(z.\theta).
\end{equation}}
\end{itemize}
\end{Lemma}
\begin{pf}
Note that by \eqref{new200} and the definition of $G_n^{m}$, $\widetilde{G}_n^m$ and $L_n^m$, for $m$ sufficiently large and any $l(n)\geq m$, we have that
\begin{equation}\label{ext1}
 |L_{n}^{m}(j)-L_{n+1}^{m}(j)|_{r,h} \leq
  C_m e^{2\pi |j|\e}e^{-\delta n}, j=-l(n),\ldots,l(n)-1,
\end{equation}
\begin{equation}\label{ext1'}
|G^{m}_n-G^{m}_{n+1}|_{r,h}\leq   C_me^{-\delta n}, \ \ |\widetilde{G}^{m}_n-\widetilde{G}^{m}_{n+1}|_{r,h} \leq Ce^{-\delta n}.
\end{equation}

For any $|z-E_0|<r$ and $|\Im \theta|<h$, by \eqref{ext1'}, we let $G^{m}(z,\theta)=\lim_{n\rightarrow \infty}G_n^{m}(z,\theta)$ and $\widetilde{G}^{m}(z,\theta)=\lim_{n\rightarrow \infty}\widetilde{G}_n^{m}(z,\theta).$ We have
that $G^m(E,\theta)$ and $\widetilde{G}^m(E,\theta)$ are analytic in both $z$ and $\theta$. Note
that for $l(m)=m$, by  \eqref{ext1'},  
$$
|\widetilde{G}^{m}_m-\widetilde{G}^{m}|_{r,h} \leq Ce^{-\delta m}.
$$
Thus by Lemma \ref{eigen0},  and eigenvalue perturbation
theory, for $m$ sufficiently large with $l(m)=m$, we have ${\rm Rank}(\widetilde{G}^{m}(E_0,\theta))=2k$ for $\theta\in \T$. Since $\widetilde{G}^{m}(E_0,\theta)=DG^{m}(E_0,\theta)D^*$, we have ${\rm Rank}(G^{m}(E_0,\theta))=2k$ for $\theta\in \T$ if
$m$ is sufficiently large.  Fix such $m$  Fix also $h_*$ with $h<h_*<h_0$. We carry out the following
rank argument on the closed substrip $|\Im\theta|\leq h_*$.
Since $G^m(E_0,\cdot)$ is holomorphic on a neighborhood
of this substrip and has rank $2k$ on $\T$, there are
only finitely many points in this substrip at which its
rank is less than $2k$. Then, by analyticity, there
are at most finitely many points $\theta_1,\ldots,\theta_s$ with $|\Im\theta_i|\leq h_*$ such that
${\rm Rank}(G^{m}(E_0,\theta_i))\leq 2k-1$. For each $\theta_i$, by the
minimality of irrational rotation, if ${\rm
  Rank}(G^{m}(E_0,\theta_i+j2m\alpha))\leq 2k-1$ for all $j\in\Z$, then
${\rm Rank}(G^{m}(E_0,\theta))\leq 2k-1$  on the strip $|\Im\theta|\leq h_*$, so
we get a contradiction. Hence there are $j_i$ such that
$$
{\rm Rank}(G^{m}(E_0,\theta_i+j_i 2m\alpha))=2k.
$$
Taking $n_0$ sufficiently large, so that $n_0>2\max\{(|j_i|+1)2m\}$ and $l(n_0)=n_0$,
notice that by \eqref{lana4}, shift covariance of the Green's
function, and the definition of $G^{m}_n(E_0,\theta)$, we have that
$G^{m}(E_0,\theta_i+j_i 2m\alpha)$ is a submatrix of $G^{n_0}.$ Therefore
$$
{\rm Rank}(G^{n_0}(E_0,\theta))\geq \max_{j\in\{0,j_1,\ldots,j_s\}}{\rm Rank}(G^{m}(E_0,\theta+j2m\alpha))=2k
$$
on  $|\Im\theta|<h_*$.

By \eqref{lana2} and \eqref{lanag}, $G_nS_nG_n=G_n$ and $\operatorname{tr}(G_nS_n)=2k$. These identities and the invariance of
$\operatorname{Ran}G_n=E_n^c$ established in Lemma~\ref{add1},
$(I-G_n(z,\theta+\alpha)S_n)\widehat S_z^{v_n}(\theta)G_n(z,\theta)=0$, extend analytically from real parameters. All $(2k+1)$-minors vanish also in the limit. By compactness, after decreasing $r>0$ if necessary,
$G^{n_0}(z,\theta)$ has constant rank $2k$ on
$D_{2r}\times\Omega_{h_*}$.
We construct the holomorphic sections below on this
larger domain and then restrict them to
$D_r\times\Omega_h$.

By Theorem \ref{trivial}, ${\rm Ker}G^{n_0}$ and $E/{\rm Ker}G^{n_0}$ where $E=(D_r\times \Omega_h)\times \C^{2n_0}$  are trivial. The quotient frame lifts holomorphically on this Stein domain. Thus there are globally defined linearly independent holomorphic functions $v_j(z,\theta)\in {\rm Ker}G^{n_0} (1\leq j\leq 2n_0-2k)$ and $e_j(z,\theta) (1\leq j\leq 2k)$ such that they form a basis of $C^{2n_0}$ for each $(z,\theta)\in D_r\times \Omega_h$.

On the compact set
$\overline{D_r}\times\{|\Im\theta|\leq h\}$, the smallest
singular value of
$G^{n_0}(e_1,\ldots,e_{2k})$ has a positive lower bound.
Thus \eqref{ext1} implies that
$G_n^{n_0}(e_1,\ldots,e_{2k})$ has rank $2k$ throughout
this set for all sufficiently large $n$.
 Thus
\begin{align*}
 &{\rm Rank}(L_n^{n_0}(z,\theta)(e_1(z,\theta),e_2(z,\theta),\cdots,e_{2k}(z,\theta)))\\
 &\geq  {\rm Rank}(G_n^{n_0}(\theta)(e_1(z,\theta),e_2(z,\theta),\cdots,e_{2k}(z,\theta)))=2k.
\end{align*}
By Lemma \ref{add1} columns of $G_n(z,\theta)$ belong to
$E_c^n(\theta)$, so by the definition of $L_n^{n_0}(z,\theta)$ so do all
its columns. Therefore we have that $L_{n}^{n_0}(z,\theta)e_j(z,\theta)\in E_c^n(\theta)$.
\end{pf}

For any $|z-E_0|<r$ and $|\Im \theta|<h$, we let
\begin{align}\label{so01}
\nonumber O_n'(z,\theta)&=-H_{n}(z,\theta)J_{2k}H_{n}^{n_0}(\bar{z},\bar{\theta})^*(e_1(z,\theta),e_2(z,\theta), \cdots,e_{2k}(z,\theta))\\
&=L_n^{n_0}(z,\theta)(e_1(z,\theta),e_2(z,\theta),\cdots,e_{2k}(z,\theta)).
\end{align}
By Lemma \ref{add2}, the columns of $O_n'(z,\theta)$ form a basis of $E_c^n(z,\theta)$. Moreover, by \eqref{ext1} and \eqref{so01},
\begin{equation}\label{final100}
|O'_n (j)-O'_{n+1}(j)|_{r,h} \leq C_{n_0}e^{-\delta n}, \ \  j=-n_0,\ldots,n_0-1
\end{equation}
where $O'_{n+1}(j)$ is the $j$-th row of $O'_{n+1}$.

On the other hand, by Lemma \ref{uvec} and direct calculation we have
\begin{align*}
O_n'(E_0,\bar{\theta})^*S_n O_n'(E_0,\theta)=&-\begin{pmatrix}e_1(E_0,\bar{\theta})&e_2(E_0,\bar{\theta)}&\cdots& e_{2k}(E_0,\bar{\theta})\end{pmatrix}^*\\
&\cdot G_n^{n_0}(E_0,\theta)\begin{pmatrix}e_1(E_0,\theta)&e_2(E_0,\theta)& \cdots&e_{2k}(E_0,\theta)\end{pmatrix}
\end{align*}
By a similar argument as in Lemma \ref{uvec}, 
$iO_n'(E_0,\bar{\theta})^*S_n O_n'(E_0,\theta)$ has $k$ positive and $k$ negative eigenvalues for all $\theta\in\T$.

Let 
\begin{align}\label{lana16}
\nonumber \Omega'(E_0,\theta)=&-\begin{pmatrix}e_1(E_0,\bar{\theta})&e_2(E_0,\bar{\theta)}&\cdots& e_{2k}(E_0,\bar{\theta})\end{pmatrix}^*
G^{n_0}(E_0,\theta)\\
&\cdot\begin{pmatrix}e_1(E_0,\theta)&e_2(E_0,\theta)& \cdots&e_{2k}(E_0,\theta)\end{pmatrix}=\lim_{n\rightarrow\infty}O_n'(E_0,\bar{\theta})^*S_n O_n'(E_0,\theta).
\end{align}
By Lemma \ref{eigen0} and convergence,  
$i\Omega'(E_0,\theta)$ has $k$ positive and $k$
negative eigenvalues for all $\theta\in\T$. A vector in the range of $G^{n_0}(E_0,\theta)$ is orthogonal to $\ker G^{n_0}(E_0,\bar\theta)$. If also orthogonal to the $e_j(E_0,\bar\theta)$, it vanishes. Thus \eqref{lana16} is nondegenerate on the strip. By the same argument as in Lemma \ref{uvec}, there is $P\in C^\omega_h(\T,GL(2k,\C))$ such that
\begin{equation}\label{lana15}
P(E_0,\bar{\theta})^*\Omega'(E_0,\theta)P(E_0,\theta)=J_{2k}.
\end{equation}
The calculation at real energies extends analytically to $z$, with adjoints at $(\bar z,\bar\theta)$. The forms converge exponentially by \eqref{ext1'}. One can view $P(\bar{z},\bar{\theta})^*O_n'(\bar{z},\bar{\theta})^*S_n O_n'(z,\theta)P(z,\theta)$ as perturbations of $J_{2k}$. Set
\[
X_n(z,\theta)=J_{2k}^{-1}P(\bar{z},\bar{\theta})^*O_n'(\bar{z},\bar{\theta})^*S_n O_n'(z,\theta)P(z,\theta),\qquad Q_n(z,\theta)=X_n^{-1/2}(z,\theta),
\]
using the power-series branch around the identity. Here $P(z,\theta)=P(E_0,\theta)$.

For sufficiently large \(n\) and sufficiently small $r$, \(Q_n(z,\theta)\) is analytic and
\[
|Q_n-I|_{r,h}\leq C|P^*(\bar{\cdot},\bar{\cdot})O_n'^*(\bar{\cdot},\bar{\cdot})S_n O_n'(\cdot,\cdot)P(\cdot,\cdot)-J_{2k}|_{r,h}.
\]
Since \(J_{2k}^*=-J_{2k}\) and \(J_{2k}^2=-I\),
\[
X_n(\bar{z},\bar{\theta})^*J_{2k}=J_{2k}X_n(z,\theta).
\]
The power-series construction then gives
\[
Q_n(\bar{z},\bar{\theta})^*J_{2k}=J_{2k}Q_n(z,\theta).
\]
Consequently,
\begin{align}\label{gjsep1}
\nonumber &Q_n(\bar{z},\bar{\theta})^{*}P(\bar{z},\bar{\theta})^*O_n'(\bar{z},\bar{\theta})^*S_n O_n'(z,\theta)P(z,\theta)Q_n(z,\theta)\\
=&
Q_n(\bar{z},\bar{\theta})^*J_{2k}X_n(z,\theta)Q_n(z,\theta)
=
J_{2k}Q_n(z,\theta)X_n(z,\theta)Q_n(z,\theta)
=
J_{2k}.
\end{align}
We define
\begin{equation}\label{so2}
O_n(z,\theta)=O_n'(z,\theta)P(E_0,\theta)Q_n(z,\theta)
\end{equation}
then by \eqref{gjsep1}.
$$
O_n(\bar{z},\bar{\theta})^*S_nO_n(z,\theta)=J_{2k},
$$
The power series defining $Q_n$ is Lipschitz on a fixed ball about $I$, so $Q_n$ converges exponentially. By \eqref{ext1}, \eqref{so01} and \eqref{final100}, the $j$th row has a holomorphic limit $O(j)$ and, for each fixed $m$, 
\begin{equation}\label{new20211}
|O_n(j)-O(j)|_{r,h}\leq C_me^{-\delta n},\ \  -m\leq j\leq m-1.
\end{equation}

Finally we let $M_n(z,\theta),$ existing by invariance of $E_c(z,\theta)$
be defined by 
\begin{equation}\label{new20222}
\begin{pmatrix}O_n(z,\theta,l(n))\\ O_n(z,\theta,l(n)-1)\\ \vdots\\ O_n(z,\theta,-l(n)+1)\end{pmatrix}=\widehat{S}_z^{v_n}(\theta)O_n(z,\theta)=O_n(z,\theta+\alpha)M_n(z,\theta).
\end{equation}
The middle $2n_0$ rows of $O_n(\theta+\alpha)$ have uniformly bounded left inverses by the lower singular-value bound in \eqref{gn1'} and the normalization. Restricting to these rows, \eqref{new20211} gives a limit $M$, Hermitian-symplectic by Proposition~\ref{p1}, with 
\begin{equation}\label{lana19}
|M_n-M|_{r,h}\leq Ce^{-\frac{\delta}{2}n}.
\end{equation}

If $0<h'<L(E_0)/(2\pi)$, continuity gives an interval $I_{h'}\ni E_0$ with $L^{v_n}(E)>2\pi h'$ for $E\in I_{h'}$ and all sufficiently large $n$.
By Theorem \ref{strip} and continuity of the Lyapunov exponent, for any $E\in I_{h'}$, we have
\begin{align*}
L_i(M(E,\cdot+i\e))&=\lim_{n\rightarrow\infty}L_{i}(M_n(E,\cdot+i\e))=\lim_{n\rightarrow\infty}L_{l(n)-k+i}(\widehat{S}_E^{v_n}(\cdot+i\e))\\
&=\lim_{n\rightarrow\infty}L_{l(n)-k+i}(\widehat{S}_E^{v_n}),\ \ \forall |\e|<h'.
\end{align*}
for $1\leq i\leq k$. At $\e=0$ this is $L_i(M)$. The remaining exponents follow by reflection, $L_{2k+1-i}(M(\cdot+i\e))=-L_i(M(\cdot-i\e))$.
\end{pf}

\noindent
{\bf Proof of Theorem \ref{hss}:} For trigonometric $v$, project $O_{E_0}$ onto the nearby holomorphic center along its dominated complement on a closed substrip, using the identity if $k=d$, and normalize as in \eqref{gjsep1} by the holomorphic inverse square root near the identity. The recurrence gives the full solutions, and Theorem~\ref{strip}, reflection and the finite Jensen formula give the stated exponents for $2\pi h<\inf_{E\in I_h}L(E)$. For infinite-range $v$, set $O_j^i(E,\theta)$ to be the $i$th component of vector $O(E,\theta,j)\in\mathbb{C}^{2k}$ defined by \eqref{new20211}. Keeping $n_0$ and the $e_j$ fixed, \eqref{final3}, \eqref{lana4}, \eqref{so01} and \eqref{so2} give $\|O_n(j)\|\le Ce^{2\pi\varepsilon|j|}$ for the full solution sequences. Since $\sum_{\ell\in\mathbb{Z}}|\widehat v_\ell|e^{2\pi\varepsilon|j+\ell|}<\infty$, \eqref{new20211} gives the dual equation by dominated convergence. Equations \eqref{new20222}--\eqref{lana19} give the shift identity.

For $E\in I_{h'}\cap\Sigma_{v,\alpha}$ with $L(E)>0$ and $\omega(E)=k$, the first turning point is $\varepsilon_1(E)=0$: otherwise the positive-Lyapunov-exponent cocycle would be regular, hence uniformly hyperbolic, contrary to $E\in\Sigma_{v,\alpha}$. The multiplicative Jensen formula (Theorem~\ref{1general}) then gives $\gamma_1^v(E)=\cdots=\gamma_k^v(E)=0$. The convergence of the center exponents above, together with Hermitian-symplectic pairing, gives $L_i(M)=0$ for $1\le i\le 2k$.\qed

\section{Further characterization of the duals of Type I operators}\label{section6}

For Type~I energies, the intrinsic dual center constructed
in the preceding section is two-dimensional. We now show
that its dynamics has a canonical form: it is \emph{projectively real}, that is in a suitable
analytic frame, the center cocycle is a scalar phase times
a real $SL(2,\mathbb R)$ cocycle. We first prove this for
trigonometric potentials and then establish its persistence
 under trigonometric approximation.

\subsection{ Projectively Real Cocycles}\label{7.1} We define the general notion of projectively real
cocycles and give a criterion for analytic conjugacy to this form. The decomposition
will be used to define the rotation pair
in Section \ref{8}.

\begin{Definition}[Projectively Real]\label{proj_real}An
  analytic cocycle $(\alpha, A)$ with $A \in C^\omega(\mathbb{T},
  GL(2, \mathbb{C}))$ is called \textbf{projectively real} if it
  admits a decomposition of the
  form:
  \begin{equation}\label{eq:decomp}A(\theta) = e^{2\pi i\phi(\theta)}
    M(\theta),
\end{equation}
where $\phi: \mathbb{T} \to \mathbb{R}$ is
  a real-analytic phase function and $M: \mathbb{T} \to SL(2,
  \mathbb{R})$ is a real-analytic cocycle.
\end{Definition}

Clearly, the scalar factor, which encodes all the complex nature of a projectively real system, disappears on projectivization:
$A(\theta)$ and $M(\theta)$ induce the same map on
$\mathbb{CP}^1$. Since $M(\theta)$ is real, this map
preserves $\mathbb{RP}^1$. The next theorem characterizes
cocycles that admit this form after an analytic change
of coordinates.

\begin{Theorem}\label{thm:characterization}Let $(\alpha, A)$ be an
  analytic cocycle. It is \textbf{projectively real} up to an (analytic) conjugation if and only if it
  satisfies the following two conditions:
\begin{enumerate}\item 
The scalar
  $GL(1,\C)$ cocycle $(\alpha, \det A)$ is subcritical. That is, there exists
  $h>0$ such that the complexified Lyapunov exponent of the
  determinant vanishes on the strip:
\begin{equation}\label{ds}L_\e(\alpha, \det A) := \int_{\mathbb{T}} \ln
  |\det A(\theta+i\epsilon)| d\theta = 0 \quad \text{for all }
  |\epsilon| < h.
\end{equation}
\item 
The cocycle preserves an analytic anti-Hermitian form of signature $(1,1)$: there
  exists an analytic  map $iH: \mathbb{T} \to \text{Herm}(2,
  \mathbb{C})$ with indefinite signature such that for all $\theta$:
\begin{equation}\label{herm}A(\theta)^* H(\theta + \alpha) A(\theta) = H(\theta).
\end{equation}
\end{enumerate}
\begin{Remark}\label{gjlast}
If $H(\theta)=J_2$, then analytic conjugation is $I_2$.
\end{Remark}
\end{Theorem}\begin{pf}If $A$ is projectively real up to an analytic conjugation, it admits the
  decomposition $B(\theta+\alpha)^{-1}A(\theta)B(\theta)= e^{2\pi i\phi(\theta)} M(\theta)$ with $\phi\in C^\omega(\T,\R)$, $M \in C^\omega(\T,SL(2, \mathbb{R}))$ and $B\in C^\omega(\T,GL(2,\C))$. We analyze the determinant $\det A(z)$ on the
  strip $z = \theta+i\epsilon$. Since $\det M(z) \equiv 1$ and $\int_\T\phi(\theta)d\theta$ is real valued, by
  analyticity, we have
\begin{equation}
L_\e(\alpha, \det A)=-4\pi \Im \int_{\mathbb{T}} \phi(\theta+i\e)d\theta+\int_{\mathbb{T}} \ln |\det
  M(\theta+i\epsilon)|d\theta=0.
\end{equation}
Furthermore, $H(\theta)=B^{-1}(\theta)^*JB^{-1}(\theta)$ induces an indefinite anti-Hermitian form for $A(\theta)$.

Conversely, by the same argument as in Lemma \ref{uvec}, there exists $B\in C^\omega(\T,GL(2,\C))$ such that
\[
B(\theta)^{*}H(\theta)B(\theta)=J.
\]
The conjugated matrix \(A_{0}(\theta)=B(\theta+\alpha)^{-1}A(\theta)B(\theta)\) then satisfies \(A_{0}(\theta)^{*}JA_{0}(\theta)=J\). the
  subcriticality condition $\int_{\mathbb{T}} \ln |\det A(\theta+i\e)|d\theta= 0$ implies $\int_{\mathbb{T}} \ln |\det A_0(\theta+i\e)|d\theta= 0$. Hence, via 
  Jensen's formula, the degree of the determinant map is
  zero. Thus we can  choose an analytic square root \(s(\theta)=e^{2\pi i\phi(\theta)}\) of \(\det A_{0}\) where $\phi\in C^\omega(\T,\R)$, unimodular on the real axis. Then
\[
C(\theta)=s(\theta)^{-1}A_{0}(\theta),\qquad
\det C(\theta)=1,\qquad
C(\theta)^{*}JC(\theta)=J.
\]
In dimension two these imply that \(C(\theta)\) is real. \end{pf}
\begin{Remark}This characterization separates the scalar phase from the real matrix
dynamics. In Section \ref{8}, we will utilize this decomposition to
  define the rotation pair for the dual center using the average phase $\hat{\phi}$ and the standard rotation number of the real component $M$.\end{Remark}
\begin{Remark}\label{eliex}
Signature $(1,1)$ is essential in (2). Eliasson's analytic uniquely ergodic $SO(3)$
skew flows \cite{eli-so3} yield, upon lifting  to $SU(2)$, analytic
two-dimensional complex cocycles with vanishing Lyapunov
exponents that are not projectively real, even up to analytic
conjugacy.

Indeed, let $U$ be such an $SU(2)$ cocycle. A projectively real
conjugacy would provide an invariant Hermitian form $K(x)$ of
signature $(1,1)$ satisfying
\[
U(x)^*K(x+\alpha)U(x)=K(x).
\]
Since $U(x)$ is unitary, the nonzero traceless part
\[
K_0(x)=K(x)-\frac12\operatorname{tr}K(x)\,I_2
\]
satisfies
\[
K_0(x+\alpha)=U(x)K_0(x)U(x)^*.
\]
Conjugation by $SU(2)$ matrices acts as rotations on the
three-dimensional real space of traceless Hermitian matrices.
Consequently, $K_0$ defines a nonzero invariant vector field
for the associated $SO(3)$ cocycle, contradicting ergodicity
of its group extension.
\end{Remark}

\subsection{The trigonometric polynomial case}\label{4.2}
We give a more precise  characterizations of the center of duals
of Type I cocycles. 
\begin{Corollary} \label{c41}Assume $\alpha\in \R\backslash\Q$,
    and $E\in\R$ is Type I. Then $M$ given by  \eqref{ff6}
  is projectively real.\end{Corollary}
\begin{Remark}
If $v$ is even, $\hat{H}_{v,\alpha,\theta}$ has real coefficients, thus
$\widehat{A}_E(\theta)$ is real valued, so $O(E,\theta)$, $M(E,\theta)$
in \eqref{ff6} can be also chosen real valued, and there is nothing to prove.
\end{Remark}

\begin{pf}
We verify that $M_E$ satisfies the two criteria of Theorem
  \ref{thm:characterization}. Let
\[
L_1(\varepsilon)\geq\cdots\geq L_{2d}(\varepsilon)
\]
be the Lyapunov exponents of
$(\alpha,\widehat{A}_E(\cdot+i\varepsilon))$, and write \(L^{d-1}(\varepsilon)=\sum_{j=1}^{d-1}L_j(\varepsilon).\)
 By Theorem~\ref{strip} (with $L^0=0$ if $d=1$),
\(
L^{d-1}(\varepsilon)=L^{d-1}(0)
\)
for all sufficiently small $|\varepsilon|$.

The Hermitian-symplectic symmetry gives
\(
L_j(\varepsilon)
=
-L_{2d+1-j}(-\varepsilon),
\)
and the constant unimodular determinant of $\widehat A_E$ gives
\(
\sum_{j=1}^{2d}L_j(\varepsilon)=0.
\)
The two Lyapunov exponents of the center cocycle
$M_E(\cdot+i\varepsilon)$ are $L_d(\varepsilon)$ and 
$L_{d+1}(\varepsilon)$. Hence
\[
\int_{\mathbb T}
\log|\det M(E,x+i\varepsilon)|dx=
L_d(\varepsilon)+L_{d+1}(\varepsilon)=
L^{d-1}(-\varepsilon)-L^{d-1}(\varepsilon)=0.
\]
Thus we have \eqref{ds}.

By Proposition \ref{p1}, $M$ preserves $J_2.$ Define the induced
anti-Hermitian form by $H(u,v) =u^*J_2 v$. It has signature
$(1,1)$ and constant determinant, so we also have \eqref{herm}. 
\end{pf}

\begin{Corollary} \label{c4}
Assume $\alpha\in \R\backslash\Q$,
  and $E_0\in\R$ is Type I. There exist  a neighborhood $I$ of $E_0$, $F\in C^\omega(I\times \T,Sp_{2d\times 2}(\C))$,  $\phi\in C^\omega(I\times \T,\R),$ and ${\bf C}\in C^\omega(I\times \T,SL(2,\R))$ such that
\begin{equation}\label{aaa18}
\widehat{A}_E(\theta)F(E,\theta)=F(E,\theta+\alpha)e^{2\pi i\phi(E,\theta)}{\bf C}(E,\theta).
\end{equation}
\\ If $L(E_0)>0$, for any $0<h'<L(E_0)/2\pi$, there is $E_0\in I_{h'}\subset I$, such that for any $E\in I_{h'}$, 
$$
L_1({\bf C}(E,\cdot+i\e))=L_1({\bf C}(E,\cdot)),\ \ \forall |\e|<h'.
$$

\end{Corollary}

\begin{pf}
Take the jointly holomorphic finite-range frame $F$ from the proof of Theorem~\ref{hss} with $k=1$, and $M$ from \eqref{ff6}. By Corollary~\ref{c41}, $\det M(E_0,\cdot)$ has winding zero. Thus $4\pi i\phi=\log\det M$ exists on the product after shrinking $I$, and $e^{-2\pi i\phi}M\in SL(2,\R)$ 
 for real $E\in I$ and $\theta\in\mathbb T,$ by Remark~\ref{gjlast}. Therefore, by \eqref{ff6}, 
\eqref{aaa18} holds.

Finally, by the property of invariant decomposition and Theorem \ref{strip}, if $L(E_0)>0$, for any $0<h'<L(E_0)/2\pi$, there is $E_0\in I_{h'}\subset I$, such that for any $E\in I_{h'}$ and $|\e|<h'$, we have
\begin{align*}
L_1({\bf C}(E,\cdot+i\e))&=L_1(M(E,\cdot+i\e))=L_d(\widehat{A}_E(\cdot+i\e))=L_d(\widehat{A}_E)=L_1(M(E,\cdot))=L_1({\bf C}(E,\cdot)).
\end{align*}

\end{pf}
Thus, for  trigonometric polynomial $v$ and a neighborhood $I$ of Type I energy $E_0$,  we obtain a family of
$Sp(2,\C)$-cocycle, corresponding to the 2-dimensional center of
$(\alpha,\widehat{A}_E),$  that is of a more special form: $(\alpha,e^{2\pi i\phi(E,\cdot)}{\bf C}(E,\cdot)),$ where
$\phi(E,\cdot)$, ${\bf C}(E,\cdot)$ are as in Corollary \ref{c4}.
\subsection{The general case}
In this subsection, $v\in C^\omega(\T,\R)$ is a real analytic function. For any
$$
0<h<\delta(E_0)=\frac{\e_2(E_0)-\e_1(E_0)}{200\pi}+L(E_0)/2\pi.
$$
\begin{Theorem} \label{theorem-main1}
Assume $\alpha\in \R\backslash\Q$, $v\in C^\omega(\T,\R)$ and
$E_0$ is Type I. There exist a neighborhood $I$ of $E_0$ and, for all sufficiently large $n$, $(U_n, V_n)\in C^\omega(I\times\T,Sp_{2l(n)\times2}(\C))$, $\phi_n\in C^\omega(I\times \T,\R)$ and ${\bf C}_n\in C^\omega(I\times \T,SL(2,\R))$, satisfying Corollary \ref{c4}. I.e., 
$$
\widehat{S}_E^{v_n}(\theta)\begin{pmatrix}U_n(E,\theta)&V_n(E,\theta)\end{pmatrix}=\begin{pmatrix}U_n(E,\theta+\alpha)&V_n(E,\theta+\alpha)\end{pmatrix}e^{2\pi i\phi_n(E,\theta)}{\bf C}_n(E,\theta).
$$
Moreover, we have
$$
|\phi_n-\phi|_{r,h}, \ \ |{\bf C}_n-{\bf C}|_{r,h}\leq Ce^{-cn},
$$
for some $C,c,r>0$,  $\phi\in C^\omega(I\times \T,\R)$ and $C\in
C^\omega(I\times \T,SL(2,\R))$. If $L(E_0)>0,$ for any $0<h'<L(E_0)/2\pi$, there is $E_0\in I_{h'}\subset I$, such that for any $E\in I_{h'}$, we have
$$
L_1({\bf C}(E,\cdot+i\e))=L_1({\bf C}(E,\cdot)),\ \ \forall |\e|<h'.
$$
Moreover, for
$$
U_n(E,\theta)=:\begin{pmatrix}U_n(E,\theta,l(n)-1)\\ U_{n}(E,\theta,l(n)-2)\\ \vdots\\ U_n(E,\theta,-l(n))\end{pmatrix},\ \ V_n(E,\theta)=:\begin{pmatrix}V_n(E,\theta,l(n)-1)\\ V_{n}(E,\theta,l(n)-2)\\ \vdots\\ V_n(E,\theta,-l(n))\end{pmatrix},
$$ 
there is $U(E,\theta,0), V(E,\theta,0)\in C^\omega(I\times \T,\C)$ such that
$$
|U_n(0)-U(0)|_{r,h},\ \ |V_n(0)-V(0)|_{r,h}\leq Ce^{-cn}.
$$
\end{Theorem}
\begin{pf}
By Theorem \ref{theorem-main1general}, one can find $U_n, V_n\in C^\omega(I\times \T,\C^{2l(n)})$ and $M_n\in C^\omega(I\times \T, Sp(2,\C))$ such that
$$
\begin{pmatrix}U^*_n(E,\theta)\\ V^*_n(E,\theta)\end{pmatrix}S_n\begin{pmatrix}U_n(E,\theta)&V_n(E,\theta)\end{pmatrix}=J_{2},
$$
\begin{equation}\label{new202}
\widehat{S}_E^{v_n}(\theta)(U_n(E,\theta),V_n(E,\theta))=(U_n(E,\theta+\alpha),V_n(E,\theta+\alpha))M_n(E,\theta),
\end{equation}
and there are $U(0), V(0)\in C^\omega(I\times \T,\C)$ and $M\in C^\omega(I\times \T, Sp(2,\C))$ such that
\begin{equation}\label{new2021}
|U_n(0)-U(0)|_{r,h},\ \ |V_n(0)-V(0)|_{r,h}\leq Ce^{-cn},\ \  |M_n-M|_{r,h}\leq Ce^{-cn}
\end{equation}
As in the Corollary \ref{c4}, we can define 
\begin{equation}\label{xz}
{\bf C}_n(E,\theta)=\frac{1}{\sqrt{\det{M_n}(E,\theta)}}M_n(E,\theta),\ \ \sqrt{\det{M_n(E,\theta)}}=e^{2\pi i\phi_n(E,\theta)}.
\end{equation}
By Type~I stability and Corollary~\ref{c41}, $\deg\det M_n=0$ for all large $n$, hence $\deg\det M=0$ by convergence. Choose $4\pi i\phi_n=\log\det M+\log(\det M_n/\det M)$ with the second logarithm chosen so that $\log 1=0$. Then we have
$$
\widehat{S}_E^{v_n}(\theta)(U_n(E,\theta),V_n(E,\theta))=(U_n(E,\theta+\alpha),V_n(E,\theta+\alpha))e^{2\pi i\phi_n(E,\theta)}{\bf C}_n(E,\theta).
$$
By \eqref{new2021} and \eqref{xz}, we have
that there exist $\phi\in C^\omega(I\times \T,\R)$ and ${\bf C}\in C^\omega(I\times \T,SL(2,\R))$ such that
$$
|\phi_n-\phi|_{r,h}, \ \ |{\bf C}_n-{\bf C}|_{r,h}\le Ce^{-cn}.
$$
By Theorem \ref{theorem-main1general}, if $L(E_0)>0$, for any $0<h'<L(E_0)/2\pi$, there is $E_0\in I_{h'}\subset I$, such that for any $E\in I_{h'}$, we have
$$
L_1(C(E,\cdot+i\e))=L_1(M(E,\cdot+i\e))-\frac{1}{2}\int_{\T}\ln|\det{M(E,\theta+i\e)}|d\theta=L_1(M(E,\cdot)), \ \  \forall |\e|<h'.
$$
\end{pf}
Finally, this immediately implies
\begin{Theorem} \label{hsspr}
  For any $\alpha\in \R\backslash\Q$, $v\in C^\omega(\T,\R)$ and
Type I energy $E$ the symplectic structure $M$ defined in Theorem
\ref{hss} is projectively real.
\end{Theorem}

\section{The fibered rotation pair. Proof of Theorem
  \ref{t2} and  \ref{t3}}\label{8}

In this section, we extend the concept of the fibered rotation number to the duals
of Type I operators, overcoming both the absence of symmetry and the absence of a cocycle.

We use the projectively real decomposition from
Section~\ref{section6},
\[
M(E,\theta)=e^{2\pi i\phi(E,\theta)}\mathbf C(E,\theta),
\qquad \mathbf C(E,\theta)\in \mathrm{SL}(2,\mathbb R),
\]
to define the rotation pair. For a fixed frame and continuous lifts, we define
\begin{align*}
\rho_1(E)
&=\int_{\mathbb T}\phi(E,\theta)\,d\theta
  +\rho(\alpha,\mathbf C(E,\cdot)),\\
\rho_2(E)
&=\int_{\mathbb T}\phi(E,\theta)\,d\theta
  -\rho(\alpha,\mathbf C(E,\cdot)).
\end{align*}
For even $v$, one can take $\phi=0$, giving the classical pair
$(\rho,-\rho)$. In general, the pair also records the mean scalar
phase, and both components enter the localization argument.

 However, in order for this to allow to take advantage of the $\SL(2,\R)$ nature
of ${\bf C}(E,\cdot)$ the key issue becomes to show that this rotation data has the
correct spectral meaning,
namely that it satisfies the generalized rotation-number--IDS correspondence
\begin{equation}\label{nrho}
N(E_1)-N(E_2)=\rho_1(E_1)-\rho_2(E_1)-\left(\rho_1(E_2)-\rho_2(E_2)\right).
\end{equation}

We first prove it for finite-range duals, and then pass to general analytic potentials using convergence
of the center cocycles and their rotation numbers under
trigonometric approximation.

Unifying these results we obtain \eqref{nrho} for the
  general case.
This presents the generalized rotation-number--IDS relation in
this projectively real setting and serves as a key step in proving the universality of
absolute continuity and sharp H\"older regularity of the integrated
density of states through the subcriticality of the
recovered $\SL(2,\R)$ cocycle $(\alpha,{\bf C}(E,\cdot))$.

\subsection{The finite-range case}
In this subsection, we assume $v$ is a trigonometric polynomial of
degree $d$. By Corollary \ref{c4}, for any type I energy $E_0\in \R$ with $L(E_0)>0$,  we obtain a
$Sp(2,\C)$-cocycle $(\alpha,e^{2\pi i\phi(E,\cdot)}{\bf C}(E,\cdot))$ corresponding to the
2-dimensional center of $(\alpha,\widehat{A}_E)$ where $(\alpha,{\bf C}(E,\cdot))$
is regular and is homotopic to the identity for any $E\in I$ \footnote{At zero exponent use almost reducibility \cite{avila1,avila2}. At positive exponent use uniform hyperbolicity (Theorem~\ref{t2.1}): invariance preserves the degree of a projective line, forcing cocycle degree zero.}. Consider the
$SL(2,\R)$ cocycle $(\alpha,{\bf C}(E,\cdot))$ and set $\hat{\rho}(E):=\rho({\bf C}(E,\cdot)),$
where $\rho({\bf C}(E,\cdot))$ is given by \eqref{rot}. We now define two functions
that can be viewed as ``fibered rotation pair'' of
$(\alpha,\widehat{A}_E)$ and will play a key role in establishing
localization and absolute continuity of the IDS,
\begin{equation}\label{lana6}
\rho_1(E):=\int_\T \phi(E,\theta)d\theta+\hat{\rho}(E),\ \  \rho_2(E):=\int_\T \phi(E,\theta)d\theta-\hat{\rho}(E).
\end{equation}

In the following, we establish the relation between the fibered
rotation pair and integrated density of states for finite-range
operators, akin to the relation \eqref{relation} in the Schr\"odinger case.

By Corollary~\ref{dominate1} and openness of domination, there is $r>0$, such that if $|z-E_0|<r$,   $(\alpha,\widehat{A}_z)$ is still PH2. It is known that for any $z\in
\mathbb{H}_\delta=\{z:|z-E_0|<r\}\cap\mathbb{H}$, the cocycle $(\alpha,\widehat{A}_z)$ is uniformly hyperbolic, thus $d$-dominated. Hence $(\alpha,\widehat{A}_z)$ is $(d-1)$, $d$, $(d+1)$-dominated. As a consequence of dominated splitting, for any $z\in\mathbb{H}_\delta$, there exists a continuous invariant decomposition
$$
\C^{2d}=E_z^s(\theta)\oplus E_z^+(\theta)\oplus E_z^-(\theta)\oplus E_z^u(\theta),\ \ \forall \theta\in\T,
$$
which by Theorem 6.1 in \cite{ajs} and Theorem \ref{trivial} implies that there are linearly independent $\{\tilde{u}_z^i(\theta)\}_{i=1}^{d-1}\in E_z^s(\theta)$, $\tilde{u}_z(\theta)\in E_z^+(\theta)$, $\tilde{v}_z(\theta)\in E_z^-(\theta)$ and $\{\tilde{v}_z^i(\theta)\}_{i=1}^{d-1}\in E_z^u(\theta)$ depending analytically on $\theta$ and  $z$, $M_z^\pm \in C^\omega(\T,GL_{d-1}(\C))$ and $m_z^\pm \in C^\omega(\T,\C\backslash\{0\})$ such that
\begin{equation}\label{aaa11}
\widehat{A}_z(\theta)\begin{pmatrix}\widetilde{U}_z(\theta)&\tilde{u}_z(\theta)\end{pmatrix}=\begin{pmatrix}\widetilde{U}_z(\theta+\alpha)&\tilde{u}_z(\theta+\alpha)\end{pmatrix}{\rm diag}\{M^+_z(\theta),m_z^+(\theta)\},
\end{equation}
\begin{equation}\label{aaa15}
\widehat{A}_z(\theta)\begin{pmatrix}\widetilde{V}_z(\theta)&\tilde{v}_z(\theta)\end{pmatrix}=\begin{pmatrix}\widetilde{V}_z(\theta+\alpha)&\tilde{v}_z(\theta+\alpha)\end{pmatrix}{\rm diag}\{M^-_z(\theta),m_z^-(\theta)\},
\end{equation}
where
\begin{align*}
\widetilde{U}_z(\theta)=\begin{pmatrix}\tilde{u}_z^1(\theta),\cdots, \tilde{u}_z^{d-1}(\theta)\end{pmatrix},\ \ \widetilde{V}_z(\theta)=\begin{pmatrix}\tilde{v}_z^{1}(\theta), \cdots, \tilde{v}_z^{d-1}(\theta)\end{pmatrix}.
\end{align*}

Since the outer dominated splittings defining $E_z^s$ and $E_z^u$
persist for $|z-E_0|<r$, we choose their frames $\widetilde U_z,\widetilde V_z$ holomorphically on this disk. Moreover, by the invariance of each subspace we have
$$
\widehat{L}_{d+1}(z)=\Re \int_\T \ln m_z^+(\theta)d\theta, \ \ \widehat{L}_{d}(z)=\Re \int_\T \ln m_z^-(\theta)d\theta.
$$
Now, we define
\begin{equation}\label{aaa12}
\begin{pmatrix}F_+(m,\theta,z)\\
F_+(m-d,\theta,z)\end{pmatrix}=(\widehat{A}_{z})_{m}(\theta)(\widetilde{U}_z(\theta),\tilde{u}_z(\theta)),
\end{equation}
$$
\begin{pmatrix}F_-(m,\theta,z)\\
F_-(m-d,\theta,z)\end{pmatrix}=(\widehat{A}_{z})_{m}(\theta)(\widetilde{V}_z(\theta),\tilde{v}_z(\theta)),
$$
and the M matrices (as in \cite{ks,gjyz}), that we denote by $M^\pm(z,\theta)$,
\begin{equation}\label{aaa13}
M^+(z,\theta)=F_+(1,\theta,z)F^{-1}_+(0,\theta,z),
\end{equation}
$$
M^-(z,\theta)=F_-(1,\theta,z)F^{-1}_-(0,\theta,z).
$$
Finally, we define the Floquet exponents,
\begin{equation}\label{thou}
w^\pm(z)=\int_\T \ln \det{M^\pm(z,\theta)}d\theta.
\end{equation}

Note that $M^\pm(z,\theta)$ can be conjugated to ${\rm
  diag}\{M^\pm_z(\theta),m_z^\pm(\theta)\}$. We have
\begin{Proposition}\label{f1}
There exist $B_\pm\in C^\omega(\T,GL(d,\C))$  such that
$$
B_\pm^{-1}(\theta+\alpha)M^\pm(z,\theta)B_\pm(\theta)={\rm diag} \{M^\pm_z(\theta),m_z^\pm(\theta)\}.
$$
\end{Proposition}
\begin{pf}
We only prove the ``$+$'' case, the ``$-$'' case follows similarly. By \eqref{aaa11} and \eqref{aaa12}, we have
\begin{align*}
\begin{pmatrix}F_+(1,\theta,z)\\
F_+(1-d,\theta,z)\end{pmatrix}&=(\widetilde{U}_z(\theta+\alpha),\tilde{u}_z(\theta+\alpha)){\rm diag}\{M^+_z(\theta),m_z^+(\theta)\}\\
&=\begin{pmatrix}F_+(0,\theta+\alpha,z)\\
F_+(-d,\theta+\alpha,z)\end{pmatrix}{\rm diag}\{M^+_z(\theta),m_z^+(\theta)\},
\end{align*}
$$
\begin{pmatrix}F_+(0,\theta,z)\\
F_+(-d,\theta,z)\end{pmatrix}=(\widetilde{U}_z(\theta),\tilde{u}_z(\theta)).
$$
Let $B_+(\theta)=F_+(0,\theta,z).$ By \eqref{aaa13}, we have
$$
B_+^{-1}(\theta+\alpha)M^+(z,\theta)B_+(\theta)={\rm diag} \{M^+_z(\theta),m_z^+(\theta)\}.
$$
\end{pf}
By \eqref{thou} and Proposition \ref{f1}, we have for some $k_0\in\Z$,
\begin{align*}
\frac{1}{2\pi}\Im \left(w^+(z)-w^-(z)\right)=&\frac{1}{2\pi}\Im\left(\int_\T\ln \det{M_z^+(\theta)}d\theta-\int_\T\ln \det{M_z^-(\theta)}d\theta\right)\\
&+\frac{1}{2\pi}\Im\left(\int_\T\ln m_z^+(\theta)d\theta-\int_\T\ln m_z^-(\theta)d\theta\right)+k_0\alpha (\mod\Z).
\end{align*}
\begin{Lemma}\label{l51}
There exists $k_1\in\Z$ such that
$$
\frac{1}{2\pi}\lim\limits_{\Im z\rightarrow 0^+}\Im\left(\int_\T\ln \det{M_z^+(\theta)}d\theta-\int_\T\ln \det{M_z^-(\theta)}d\theta\right)=k_1\alpha (\mod\Z).
$$
\end{Lemma}
\begin{pf}
Note that $M_z^\pm(\theta)$ depend analytically on $z$ for $z\in \{z:|z-E_0|<\delta\}$.  Assume
$z=E+i\e.$ We have
\begin{align*}
\lim\limits_{\Im z\rightarrow 0^+}&\left(\int_\T\ln \det{M_z^+(\theta)}d\theta-\int_\T\ln \det{M_z^-(\theta)}d\theta\right)\\
&=\int_\T\ln \det{M_E^+(\theta)}d\theta-\int_\T\ln \det{M_E^-(\theta)}d\theta.
\end{align*}
Note that
\begin{align}\label{sym1}
&\widehat{A}_E(\theta)\begin{pmatrix}\widetilde{U}_E(\theta)&\widetilde{V}_E(\theta)\end{pmatrix}=\begin{pmatrix}\widetilde{U}_E(\theta+\alpha)&\widetilde{V}_E(\theta+\alpha)\end{pmatrix}{\rm diag}\{M^+_E(\theta), M^-_E(\theta)\}
\end{align}
Taking the adjoint, we obtain
\begin{align*}
\begin{pmatrix}\widetilde{U}^*_E(\theta)\\ \widetilde{V}^*_E(\theta)\end{pmatrix}\widehat{A}_E(\theta)^*
={\rm diag}\{M^+_E(\theta), M^-_E(\theta)\}^*\begin{pmatrix}\widetilde{U}^*_E(\theta+\alpha)\\ \widetilde{V}^*_E(\theta+\alpha)\end{pmatrix}.
\end{align*}
Therefore,
\begin{align*}
\begin{pmatrix}\widetilde{U}^*_E(\theta)\\ \widetilde{V}^*_E(\theta)\end{pmatrix}\widehat{A}_E(\theta)^*S
={\rm diag}\{M^+_E(\theta), M^-_E(\theta)\}^*\begin{pmatrix}\widetilde{U}^*_E(\theta+\alpha)\\ \widetilde{V}^*_E(\theta+\alpha)\end{pmatrix}S.
\end{align*}

Since
$$
\widehat{A}_E(\theta)^*S=S\widehat{A}_E(\theta)^{-1},
$$
it follows
\begin{align*}
\begin{pmatrix}\widetilde{U}^*_E(\theta)\\ \widetilde{V}^*_E(\theta)\end{pmatrix}S
={\rm diag}\{M^+_E(\theta), M^-_E(\theta)\}^*\begin{pmatrix}\widetilde{U}^*_E(\theta+\alpha)\\ \widetilde{V}^*_E(\theta+\alpha)\end{pmatrix}S\widehat{A}_E(\theta).
\end{align*}
Therefore,
\begin{align*}
&\begin{pmatrix}\widetilde{U}^*_E(\theta)\\ \widetilde{V}^*_E(\theta)\end{pmatrix}S\begin{pmatrix}\widetilde{U}_E(\theta)&\widetilde{V}_E(\theta)\end{pmatrix}\\
=&{\rm diag}\{M^+_E(\theta), M^-_E(\theta)\}^*\begin{pmatrix}\widetilde{U}^*_E(\theta+\alpha)\\ \widetilde{V}^*_E(\theta+\alpha)\end{pmatrix}S\widehat{A}_E(\theta)\begin{pmatrix}\widetilde{U}_E(\theta)&\widetilde{V}_E(\theta)\end{pmatrix}\\
=&{\rm diag}\{M^+_E(\theta), M^-_E(\theta)\}^*\begin{pmatrix}\widetilde{U}^*_E(\theta+\alpha)\\ \widetilde{V}^*_E(\theta+\alpha)\end{pmatrix}S\begin{pmatrix}\widetilde{U}_E(\theta+\alpha)&\widetilde{V}_E(\theta+\alpha)\end{pmatrix}{\rm diag}\{M^+_E(\theta), M^-_E(\theta)\}.
\end{align*}
By symplectic orthogonality, we may assume
\begin{align*}
\begin{pmatrix}\widetilde{U}^*_E(\theta)\\ \widetilde{V}^*_E(\theta)\end{pmatrix}S\begin{pmatrix}\widetilde{U}_E(\theta)&\widetilde{V}_E(\theta)\end{pmatrix}=\begin{pmatrix}
0&D(\theta)\\ -D^*(\theta)&0
\end{pmatrix}.
\end{align*}
Note that
$$
D(\theta)=M_E^+(\theta)^*D(\theta+\alpha)M_E^-(\theta).
$$
This implies the result.
\end{pf}

Let $r>0$ be the radius of the complex energy neighborhood
of $E_0$ fixed above, decreased if necessary.

\begin{Lemma}\label{l52}
There exist integers $k_+,k_-$, such that for every
$E\in(E_0-r,E_0+r)$,
\begin{equation}\label{eq:stable-phase-corrected}
 \lim_{\eta\rightarrow0^+}\frac1{2\pi}\operatorname{Im}
 \int_{\mathbb T}\log m_{E+i\eta}^+(\theta)\,d\theta
 =\rho_1(E)-k_+\alpha\pmod{\mathbb Z},
\end{equation}
\begin{equation}\label{eq:unstable-phase-corrected}
 \lim_{\eta\rightarrow 0^+}\frac1{2\pi}\operatorname{Im}
 \int_{\mathbb T}\log m_{E+i\eta}^-(\theta)\,d\theta
 =\rho_2(E)-k_-\alpha\pmod{\mathbb Z}.
\end{equation}
The convergence is uniform on compact
subintervals of $(E_0-r,E_0+r)$. 
\end{Lemma}

\begin{pf}
We only prove the result for $m_+$, the other case follows exactly the same way. We denote 
\[
 F_z(\theta)=(u_z(\theta),v_z(\theta)),\ \ 
 \widehat A_z(\theta)F_z(\theta)=F_z(\theta+\alpha)M(z,\theta),
 \ \  M(z,\theta)=e^{2\pi i\phi(z,\theta)}{\bf C}(z,\theta),
\]
where $F_E(\theta)^*SF_E(\theta)=J$, ${\bf C}(E,\cdot)\in
C^\omega(\mathbb T,SL(2,\mathbb R))$ is homotopic to the identity, and
\[
 J=\begin{pmatrix}0&1\\-1&0\end{pmatrix},\qquad
 Q=-\frac1{1+i}\begin{pmatrix}1&-i\\1&i\end{pmatrix}.
\]
For $z=E+i\eta$, $\eta>0$, let $\tilde u_z\in E_z^+(\theta)$ and
$\tilde{v}_z\in E_z^-(\theta)$ be from \eqref{aaa11} and \eqref{aaa15}. Let
\[
 (\tilde u_z(\theta),\tilde v_z(\theta))=F_z(\theta)B_z(\theta),\qquad
 QB_z(\theta)=\begin{pmatrix}a(z,\theta)&b(z,\theta)\\
 c(z,\theta)&d(z,\theta)\end{pmatrix}.
\]
 A direct computation gives
\[
\widehat A_z(\theta)^*S\widehat A_z(\theta)-S
=2i\eta\,\delta_0\delta_0^*,
\]
where $\delta_0=e_d\in\C^{2d}$  in the ordering $d-1,\ldots,0,-1,\ldots,-d$.
Set
\[
U_k(z,\theta):=(\widehat A_z)_k(\theta)\tilde u_z(\theta).
\]
Then
\[
U_{k+1}(z,\theta)^*SU_{k+1}(z,\theta)
-U_k(z,\theta)^*SU_k(z,\theta)
=2i\eta\,|\langle\delta_0,U_k(z,\theta)\rangle|^2.
\]
Since $\tilde u_z(\theta)\in E_z^+(\theta)$,
$U_k(z,\theta)$ decays exponentially as $k\to+\infty$.
Therefore
\begin{equation}\label{20242gjy}
\tilde u_z(\theta)^*S\tilde u_z(\theta)
=-2i\eta\sum_{k=0}^{\infty}
|\langle\delta_0,U_k(z,\theta)\rangle|^2.
\end{equation}
By analyticity of $F_z(\theta)$ and direct calculation, we have
\begin{equation}
\label{norm}
    F_{E+i\eta}^*SF_{E+i\eta}=J+O(\eta),\ \  QJQ^{-1}=\begin{pmatrix}i&0\\0&-i\end{pmatrix}.
\end{equation}

By \eqref{20242gjy},
$-i\tilde u_z^*S\tilde u_z\leq0$.
The same  argument in the negative direction,
using the backward decay of $\tilde v_z$, gives
$-i\tilde v_z^*S\tilde v_z\geq0$.
Since $Q$ is unitary and
\[
\tilde u_z=F_zQ^*\binom{a}{c},
\qquad
\tilde v_z=F_zQ^*\binom{b}{d},
\]
by \eqref{norm} we have
\begin{align*}
-i\tilde u_z^*S\tilde u_z
&=|a|^2-|c|^2+O(\eta)(|a|^2+|c|^2),\\
-i\tilde v_z^*S\tilde v_z
&=|b|^2-|d|^2+O(\eta)(|b|^2+|d|^2).
\end{align*}

Consequently, uniformly for $E$ in a compact subinterval and $\theta\in\mathbb T$,
\[
 |a(z,\theta)|^2-|c(z,\theta)|^2\le K\eta(|a(z,\theta)|^2+|c(z,\theta)|^2),\ \
 |b(z,\theta)|^2-|d(z,\theta)|^2\ge-K\eta(|b(z,\theta)|^2+|d(z,\theta)|^2).
\]
For $K\eta<1$, these inequalities imply that $c=0$
would force $a=0$, and $b=0$ would force $d=0$,
contradicting the invertibility of $QB_z$.
Define
$
w_z^+(\theta):=\frac{a(z,\theta)}{c(z,\theta)}.
$
We therefore have
\[
|w_z^+(\theta)|^2\leq\frac{1+K\eta}{1-K\eta}.
\]
Now define 
\[
s_z^+(\theta):=
F_z(\theta)Q^{-1}\binom{w_z^+(\theta)}{1}.
\]
Then
$
\tilde u_z(\theta)=c(z,\theta)s_z^+(\theta).
$ By invariance of the stable line, there is a unique
nonzero scalar $\mu_z^+(\theta)$ such that  
\[
 \widehat A_z(\theta)s_z^+(\theta)=\mu_z^+(\theta) s_z^+(\theta+\alpha).
\]
 Set

\[
 \widetilde {\bf C}(z,\theta)=Q{\bf C}(z,\theta)Q^{-1}
 =\begin{pmatrix}a_1(z,\theta)&b_1(z,\theta)\\
 c_1(z,\theta)&d_1(z,\theta)\end{pmatrix}.
\]
For real $E$, and $\theta\in\T$, by our choice of $Q$ we have $\widetilde{\mathbf C}(E,\theta)\in SU(1,1)$;
in particular,
$|d_1(E,\theta)|^2-|c_1(E,\theta)|^2
=|a_1(E,\theta)|^2-|b_1(E,\theta)|^2=1$. Using the  phase $\phi$ from
$M=e^{2\pi i\phi}\mathbf C$, define
\[
 \tau_z^+(\theta,w)=e^{2\pi i\phi(z,\theta)}
       (c_1(z,\theta)w+d_1(z,\theta)).
\]
On a compact real energy interval these functions are nonzero on one fixed
disk $|w|\le1+\delta$, and remain nonzero there for sufficiently small
$|\eta|$.  Since ${\bf C}(E,\cdot)$ is homotopic to the identity we have that
$a_1(E,\cdot)$ and $d_1(E,\cdot)$ have degree zero. Since
$\phi(E,\theta)$ is real valued and periodic,  $\tau_z^+$ admits
continuous periodic logarithm, which is jointly continuous in the parameters after
shrinking the neighborhood. We choose the branch compatible with the lift of
$C_E$ used to define $\hat\rho(E)$, and put
\[
a_z^+(\theta,w)
 =\frac1{2\pi}\Im\ln\tau_z^+(\theta,w).
\]
By invariance of the stable direction, comparing the second
components in the $Q$-coordinates we obtain
\[
 \mu_z^+(\theta)=\tau_z^+(\theta,w_z^+(\theta)).
\]

Fix a real energy $E$. Since $\widetilde{\mathbf C}(E,\theta)\in SU(1,1)$,
the associated M\"obius transformation

\[
 P_E^+(\theta,w)=
 \frac{a_1(E,\theta)w+b_1(E,\theta)}
      {c_1(E,\theta)w+d_1(E,\theta)},
\]
preserves $\overline{\mathbb D}$. Let $\nu$ be any invariant
probability measure for $(\theta,w)\mapsto(\theta+\alpha,P_E^+(\theta,w))$
on $\mathbb T\times\overline{\mathbb D}$ whose projection to $\mathbb T$
is Lebesgue measure. Denote by $a_{E,n}^+(\theta,w)$ the  $n$-step
sum of $a_E^+(\theta,w)$ along this skew product.  Denote 
\[
Q\mathbf C_n(E,\theta)Q^{-1}
=:\begin{pmatrix}
a_n(E,\theta)&b_n(E,\theta)\\
c_n(E,\theta)&d_n(E,\theta)
\end{pmatrix}.\]
Then
\[
a_{E,n}^+(\theta,w)=\frac{1}{2\pi}\Im \ln \left(e^{2\pi i\varphi_n(E,\theta)}
       (c_n(E,\theta)w+d_n(E,\theta))\right)
\]
where $\varphi_n(E,\theta):=
\sum_{j=0}^{n-1}\phi(E,\theta+j\alpha).$
Since $Q\mathbf C_n(E,\theta)Q^{-1}\in SU(1,1)$,
we have $|d_n(E,\theta)|^2-|c_n(E,\theta)|^2=1$. Therefore
 the disk
$d_n+c_n\overline{\mathbb D}$ avoids zero and has argument oscillation
less than $\pi$. Thus, for any fixed $w_0\in\partial\mathbb D$,
\[
 \left|a_{E,n}^+(\theta,w)
       -a_{E,n}^+(\theta,w_0)\right|\le\frac12
 \qquad(w\in\overline{\mathbb D}).
\]
Invariance of $\nu$ therefore implies
\[
 \left|\int a_E^+(\theta,w)\,d\nu
 -\frac1n\int_{\mathbb T}a_{E,n}^+(\theta,w_0)\,d\theta\right|
 \le\frac1{2n}.
\]
Let $w_0=e^{-4\pi i t_0}$,
and write $v(t)=(\cos 2\pi t,\sin 2\pi t)^T$.
Using the lift defining $\hat\rho(E)$, write
\[
\mathbf C_n(E,\theta)v(t_0)
=r_n(\theta)v(t_n(\theta)),
\qquad r_n(\theta)>0.
\]
Our choice of $Q$ gives
\[
c_n(E,\theta)w_0+d_n(E,\theta)
=r_n(\theta)e^{2\pi i(t_n(\theta)-t_0)}.
\]
Thus, 
by shift invariance of Lebesgue
measure, we have
\[
 \lim_{n\to\infty}\frac1n\int_{\mathbb T}
 a_{E,n}^+(\theta,w_0)\,d\theta
 =\int_{\mathbb T}\phi(E,\theta)\,d\theta+\hat\rho(E)
 =\rho_1(E).
\]
It follows that every such invariant measure satisfies
\[
 \int_{\T\times \overline{\mathbb{D}}} a_E^+(\theta,w)\,d\nu=\rho_1(E).
\]

For $z=E+i\eta$, define the graph measures
\[
 \nu_z^+=\int_{\mathbb T}
       \delta_{(\theta,w_z^+(\theta))}\,d\theta.
\]
 By the bound on $w_z^+$ the supports of $\nu_z^+$
lie in a common compact subset of $\T\times\C$, giving weak compactness. Any subsequential limit as
$\eta\downarrow0$ is supported on
$\mathbb T\times\overline{\mathbb D}$ and projects to Lebesgue measure \footnote{The weak limit of the invariant graph measure remains invariant under skew product. This follows by applying graph invariance to continuous test functions and using uniform convergence of the skew product map on a compact neighborhood of the supports. }.
 Since $a_z^+(\theta,w)$
converge uniformly on a neighborhood of the supports, every subsequential
limit of their averages equals $\rho_1(E).$ 
Consequently,
\[
 \lim_{\eta\rightarrow0^+}\frac1{2\pi}\operatorname{Im}
 \int_{\mathbb T}\ln\mu_{E+i\eta}^+(\theta)\,d\theta=\rho_1(E).
\]
Moreover, the convergence is uniform on compact energy
intervals: otherwise choose $E_j\to E_*$ and $\eta_j\rightarrow0^+$
violating uniform convergence, apply the same compactness argument at
$E_*$, and using continuity of the scalar mean and of the real rotation number, one gets a contradiction.

Finally, the original sections satisfy
\[
 m_z^+(\theta)=\frac{c(z,\theta)}{c(z,\theta+\alpha)}\mu_z^+(\theta).
\]
Set $k_+=\deg c(z,\cdot)$.  The ratio contributes $-k_+\alpha$ modulo $\mathbb Z$. For the unstable line use $w_z^-=d/b$, $k_-=\deg b$, and the first component of $Qv(t)$, giving $-2\pi(t_n-t_0)$ and hence $\rho_2$. Both degrees are constant in $z$ and both ratios $a/c,d/b$ are independent of the local sections.
\end{pf}

We are now ready to link $\rho_1$ and $\rho_2$ with the IDS.
\begin{Theorem}\label{t4}
For $\alpha\in \R\backslash\Q$, $E_0\in \R$ with $L(E_0)>0$ and $\omega(E_0)=1$, there exists a neighborhood $I$ of $E_0$, such that
\begin{enumerate}
\item for almost every $E\in \Sigma_{v,\alpha}\cap I$, $\rho'_1(E)\geq 0$ and $\rho'_2(E)\leq 0$;
\item $2\hat{\rho}(E)=\rho_1(E)-\rho_2(E)=N(E)+k'\alpha(\mod\Z)$ for some $k'$ depending on $E_0$.\footnote{The sign reflects our normalization $F^*SF=J$.
For example, in the nearest-neighbor case with positive
hopping $\lambda$, one has $S=-\lambda J$, and we may take
$F=\lambda^{-1/2}\operatorname{diag}(1,-1)$.
This change of coordinates reverses orientation, so the
rotation number of the normalized center cocycle is the
negative of the usual Schr\"odinger rotation number.
Thus $N=1-2\rho_{\mathrm{std}}$ becomes
$N=1+2\hat\rho$ with compatible lifts.}
\end{enumerate}
\end{Theorem}
\begin{pf}
Note that by \eqref{thou}, we have
$$
\Re w^\pm(z)=\mp \widehat{L}^d(z)=\mp L^d(\widehat{A}_z).
$$
By the Thouless formula obtained in \cite{ks}, we have
$$
\widehat{L}^d(z)=\Re \left(\int_\R\ln (E'-z)dN(E')-\ln |\hat{v}_d|\right).
$$
Both Weyl matrices in \eqref{aaa13} are holomorphic on $\mathbb H$ because a singular boundary projection would give a nonreal half-line Dirichlet eigenvalue. The half-line resolvent expansion and the shift rows in \eqref{sco} give, uniformly in $\theta$ as $|z|\to\infty$,
\[\det M^+=(-1)^{d-1}\overline{\hat v_d}/z+O(z^{-2}),\qquad
\det M^-=(-1)^{d-1}z/\hat v_d+O(1).\]
Thus $\det M^\pm$ have degree zero in $\theta,$ for large $z$, and therefore everywhere on $\mathbb H$. For their periodic logarithms, equality of real parts gives $w^\pm(z)=\mp(\int_{\mathbb R}\ln(E'-z)\,dN(E')+i\kappa^\pm-\ln|\hat v_d|)$. Comparing the determinant ratio at infinity gives $\kappa^++\kappa^-\in2\pi\mathbb Z$.

Thus
$$
\lim\limits_{\Im z\rightarrow 0^+} \frac{1}{2\pi}\Im (w^+(z)-w^-(z))=N(E)(\mod\Z).
$$
By Lemma \ref{l51}, Lemma \ref{l52} and \eqref{lana6}, we have for some $k'\in\Z$
\begin{align}\label{sum2}
2\hat{\rho}(E)=\rho_1(E)-\rho_2(E)=k'\alpha+N(E)(\mod\Z). 
\end{align} 
Thus $\hat{\rho}(E)$ is a continuous non-decreasing function, hence it
is differentiable at almost every E. Note also that  by Corollary \ref{c4}, we have $\int_\T\phi(E,\theta)d\theta$ is analytic in $E$.  Thus $\rho_1(E), \rho_2(E)$ are differentiable for almost every E. By Lemma \ref{l52}, \eqref{sum2} and Thouless formula, we have for some $k''\in\Z$ and $c_+\in\R$,
\begin{equation}\label{aaa20}
\lim\limits_{\Im z\rightarrow 0^+} \frac{1}{2\pi}\Im w^+(z)=\frac{1}{2\pi}\Im \int_\T\ln \det{M_E^+(\theta)}d\theta+\rho_1(E)+k''\alpha=\hat{\rho}(E)+c_+(\mod\Z),
\end{equation}
which implies for some $k'''\in\Z$ and $c_+'\in\R$,
$$
\frac{1}{2\pi}\Im \int_\T\ln \det{M_E^+(\theta)}d\theta=- \int_\T\phi(E,\theta)d\theta+c_+'(\mod\Z).
$$
By Thouless formula and poisson formula, for almost every $E$, we have
$$
2\hat{\rho}'(E)=N'(E)=\frac{1}{\pi}\lim\limits_{\eta\downarrow 0}\frac{L(E+i\eta)-L(E)}{\eta}.
$$

Note that for almost every $E\in \Sigma_{v,\alpha}\cap I$, we have
$\widehat{L}_d(E)=0$, so by \eqref{aaa20} and Cauchy-Riemann equations, we have 
\begin{align*}
\rho'_1(E)=&\frac{d\int_\T\phi(E,\theta)d\theta}{dE}+\hat{\rho}'(E)=-\frac{d\Im\int_\T\ln \det{M_E^+(\theta)}d\theta}{2\pi dE}+\hat{\rho}'(E)\\
=&\lim\limits_{\Im z\rightarrow 0^+}\left(\frac{\partial\Re\int_\T\ln \det{M_z^+(\theta)}d\theta}{2\pi\partial \Im z}-\frac{\partial \Re w^+(z)}{2\pi\partial \Im z}\right)= \lim\limits_{\Im z\rightarrow 0^+}\frac{-\widehat{L}_{d+1}(z)}{2\pi\Im z}\geq 0,
\end{align*}
\begin{align*}
\rho'_2(E)=&\frac{d\int_\T\phi(E,\theta)d\theta}{dE}-\hat{\rho}'(E)=-\frac{d\Im \int_\T\ln \det{M_E^+(\theta)}d\theta}{2\pi dE}-\hat{\rho}'(E)\\
=&\lim\limits_{\Im z\rightarrow 0^+}\left(\frac{\partial\Re\int_\T\ln \det{M_z^+(\theta)}d\theta}{2\pi \partial \Im z}+\frac{\partial \Re w^+(z)}{2\pi\partial\Im z}\right)=\lim\limits_{\Im z\rightarrow 0^+}\frac{-\widehat{L}_{d}(z)}{2\pi\Im z}\leq 0.
\end{align*}
Here we use $\Re w^\pm(z)=\Re\int_\T\ln \det{M_z^\pm(\theta)}d\theta+\Re\int_\T \ln m_z^\pm(\theta)d\theta$.
By Lemma~\ref{l51} and the Cauchy--Riemann equations,
the functions
$\Re\int_{\mathbb T}\log\det M_z^+(\theta)\,d\theta$, 
$\Re\int_{\mathbb T}\log\det M_z^-(\theta)\,d\theta$
have equal derivatives with respect to $\Im z$
at the real energy interval.
We also use $\widehat{L}_{d+1}(z)=\Re \int_\T \ln m_z^+(\theta)d\theta$ and $\widehat{L}_{d}(z)=\Re \int_\T \ln m_z^-(\theta)d\theta$.
\end{pf}
\subsection{The infinite-range case}
Let $v(x)=\sum_{k\in\Z} \hat{v}_k e^{2\pi ikx}$ be a real analytic function, and $v_n(x)=\sum_{k=-n}^n\hat{v}_k e^{2\pi ikx}$ be its $n$-th truncation. Note that  $L^v_\e(E_0)$ has at least one non-negative turning points $0\leq \e_1(E_0)<h_1$. Let $I$ be the neighborhood of $E_0$ in  Theorem \ref{theorem-main1}.

To specify the dependence on $n$, in this case, we denote the
integrated density of states of $H_{v_n,\alpha,\theta}$ by
$N_n(E)$. By continuity of the Lyapunov exponent and the bounds $N(E-\|v_n-v\|_\infty)\leq N_n(E)\leq N(E+\|v_n-v\|_\infty)$,
$$
L^{v_n}(E)\rightarrow L^v(E), \ \ N_n(E)\rightarrow N(E), \ \ \forall E\in I,
$$
uniformly in $E$.

By Theorem \ref{theorem-main1}, $\phi_n$ depends analytically on $E$ and $\phi_n\rightarrow\phi_E$ uniformly on a small complex neighborhood of $I$. We also have ${\bf C}_n\rightarrow {\bf C}$ where ${\bf C}\in C^\omega(I\times \T,SL(2,\R))$.

We now define $\hat{\rho}(E):=\rho({\bf C}(E,\cdot))$ and 
$$
\rho_1(E):=\int_\T \phi(E,\theta)d\theta+\hat{\rho}(E),\ \  \rho_2(E):=\int_\T \phi(E,\theta)d\theta-\hat{\rho}(E).
$$
We will call $\rho_1$ and $\rho_2$ the fibered rotation pair of $\widehat{H}_{v,\alpha,\theta}$.

\begin{Corollary}\label{c100}
Let $v\in C^\omega(\T,\R)$, $\alpha\in\R\backslash\Q$, and
$E_0\in\Sigma_{v,\alpha}$ satisfy $\bar\omega(E_0)=1$ and $L(E_0)>0$.
There is a neighborhood $I$ of $E_0$ such that
\begin{enumerate}
\item $2\hat{\rho}(E)=\rho_1(E)-\rho_2(E)=N(E)+c_I(\mod\Z)$ for some $c\in\R$ depending on $I$;
\item for almost every $E\in \Sigma_{v,\alpha}\cap I$, $\rho'_1(E)\geq 0$ and $\rho'_2(E)\leq 0$.
\end{enumerate}
\end{Corollary} 
\begin{pf}
Fix $\varepsilon$ inside the first positive-slope branch at $E_0$. Continuity and quantization preserve its slope $2\pi$ and $L^{v_n}(E)>0$ on a common $I$ for large $n$. At spectral energies of an approximant in $I$, convexity gives $\omega^{v_n}\le1$, and regularity at zero would imply uniform hyperbolicity, so $\omega^{v_n}=1$. Such energies exist by norm convergence of spectra. At those energies the corresponding center is subcritical, hence $\deg\mathbf C_n=0$ throughout $I$.
By Theorem~\ref{t4}, $2\rho(\mathbf C_n)-N_n$ is locally constant near these energies and also in gaps, where the dual center is uniformly hyperbolic and $N_n$ constant. Taking differences at two energies and passing to the limit proves (1). By (1), we have
$$
\rho_1(E)=\frac{N(E)+c_I}{2}+\int_\T \phi(E,\theta)d\theta,\ \ \rho_2(E)=-\frac{N(E)+c_I}{2}+\int_\T \phi(E,\theta)d\theta
$$
thus $\rho_1(E)$ and $\rho_2(E)$ are differentiable for almost every $E\in I$.

Note that by Theorem \ref{1general},
\begin{align*}
\Re\int_\T\ln \det{M_z^+(\theta)}d\theta- \Re w^+(z)=2\pi \e_1(z),
\end{align*}
where 
$$
\Re w^+(z)=\lim\limits_{n\rightarrow\infty}\Re \left(w^+_n(z)-\ln|\hat{v}_{l(n)}|\right)=-L^v(z)=-\int_\R \ln|z-E|dN(E)
$$ and $\Re\int_\T\ln \det{M_z^+(\theta)}d\theta=\lim\limits_{n\rightarrow\infty} \left(\Re\int_\T\ln \det{M_{z,n}^+(\theta)}d\theta-\ln|\hat{v}_{l(n)}|\right)$ \footnote{\label{har} By Theorem \ref{1general}, $\Re w^+_n(z)-\ln|\hat{v}_{l(n)}|=-L^{v_n}(z)$ and $\Re\int_\T\ln \det{M_{z,n}^+(\theta)}d\theta-\ln|\hat{v}_{l(n)}|=-L_\e^{v_n}(z)+2\pi \e$ for some fixed $\e>0$ where $\e>\sup_{|z-E_0|<r}\e_1(z)$ and $\e$ is in between the first and the second turning point of $L_e^{v}(z)$. Hence they both converge uniformly, in particular, since $(\alpha,A_z(\cdot+i\e))$ is uniformly hyperbolic, $\Re\int_\T\ln \det{M_z^+(\theta)}d\theta=-L^v_\e(z)+2\pi\e$ is harmonic in a neighborhood of $E_0$.}. By the Thouless formula, for $\eta>0$,
\[
\frac{\partial}{\partial\eta}L^v(E+i\eta)
=\int_{\mathbb R}\frac{\eta}{(E-t)^2+\eta^2}\,dN(t).
\]
The right-hand side is $\pi$ times the Poisson integral of
$dN$, and hence converges to $\pi N'(E)$ for almost every
$E$. Integrating in $\eta$ and using continuity at the
real axis, we obtain, for almost every $E\in I$,
\begin{equation}\label{tp}
N'(E)
=\frac{1}{\pi}\lim_{\eta\downarrow0}
\frac{L^v(E+i\eta)-L^v(E)}{\eta}
=-\frac{1}{\pi}\lim_{\eta\downarrow0}
\frac{\Re w^+(E+i\eta)-\Re w^+(E)}{\eta}.
\end{equation}
 By Footnote~\ref{har}, the normalized limit denoted by
\[
\Re\int_{\mathbb T}\log\det M_z^+(\theta)\,d\theta
\]
is harmonic in a complex neighborhood of $E_0$. For each large $n$, the differentiated identity \eqref{aaa20} extends from a spectral neighborhood to $I$ by real analyticity of the scalar mean and Footnote~\ref{har}. 
Applying the boundary phase identity and the
Cauchy--Riemann equations first to the approximants, and
then passing to the limit using locally uniform convergence,
gives
\[
\frac{d}{dE}\int_{\mathbb T}\phi(E,\theta)\,d\theta
=
\frac{1}{2\pi}
\left.\frac{\partial}{\partial\eta}
\Re\int_{\mathbb T}\log\det M_{E+i\eta}^+(\theta)\,d\theta
\right|_{\eta=0}.
\]
Since
\[
\rho_1(E)
=\int_{\mathbb T}\phi(E,\theta)\,d\theta
+\frac{N(E)+c_I}{2},
\]
By the multiplicative Jensen formula, $L_\varepsilon^v(z)$
is constant up to $\varepsilon_1(z)$, and the Type~I
 gives slope $2\pi$ on the following branch.
Thus, for the fixed $\varepsilon$ chosen in
Footnote~\ref{har},
\[
L_\varepsilon^v(z)
=L^v(z)+2\pi\bigl(\varepsilon-\varepsilon_1(z)\bigr).
\]
We therefore obtain
\[
\Re\int_{\mathbb T}\log\det M_z^+(\theta)\,d\theta
-\Re w^+(z)
=L^v(z)-L_\varepsilon^v(z)+2\pi\varepsilon\\
=2\pi\varepsilon_1(z).
\]
By \eqref{tp} we have for
almost every $E\in\Sigma_{v,\alpha}\cap I$,
\begin{align*}
\rho_1'(E)
&=\lim_{\eta\downarrow0}\frac{1}{2\pi\eta}
\biggl[
\Re\int_{\mathbb T}
\bigl(\log\det M_{E+i\eta}^+(\theta)
-\log\det M_E^+(\theta)\bigr)\,d\theta
-\Re w^+(E+i\eta)+\Re w^+(E)
\biggr]\\
&=\lim_{\eta\downarrow0}
\frac{\varepsilon_1(E+i\eta)-\varepsilon_1(E)}{\eta}
=\lim_{\eta\downarrow0}
\frac{\varepsilon_1(E+i\eta)}{\eta}\geq0,
\end{align*}
since $\varepsilon_1(E)=0$ at the supercritical E. At $z=E-i\eta$ the same identity and $L^v(\bar z)=L^v(z)$ give $\rho_2'(E)=-\lim_{\eta\downarrow0}\varepsilon_1(E-i\eta)/\eta\le0$.
\end{pf}

 \subsection{Proof of Theorem \ref{t2}} 

The gist of the proofs of
 Theorems \ref{t2}, \ref{t3} is in relation of the IDS and the rotation
 number of the hidden (limiting) $SL(2,\R)$ cocycle $(\alpha,{\bf C}(E,\cdot)),$
 whose regularity obtained in Theorem \ref{theorem-main1} allows to appeal
 to the results in the universal subcritical regime. Denote
 $$
 \Sigma^{ar}_{v,\alpha}=\{E\in\R:\omega(E)=0\},\ \  \Sigma^+_{v,\alpha}=\{E\in\Sigma_{v,\alpha}:L(E)>0\}.
 $$
 By the almost reducibility theorem \cite{avila1,avila2}, $\Sigma^{ar}_{v,\alpha}$ is open. Moreover, if $H_{v,\alpha,x}$ has no critical regime, then $\Sigma_{v,\alpha}\subset\Sigma^{ar}_{v,\alpha}\cup \Sigma^{+}_{v,\alpha}$.
 
  It was proved by Avila \cite{avila1,avila2} that $\rho(E)|_{\Sigma^{ar}_{v,\alpha}}$ is absolutely continuous, thus $N(E)|_{\Sigma^{ar}_{v,\alpha}}=1-2\rho(E)|_{\Sigma^{ar}_{v,\alpha}}$ is absolutely continuous. 
On the other hand, $\Sigma^+_{v,\alpha}$ can be covered by
countably many intervals $I_i$ from Theorem~\ref{theorem-main1}.
On each $I_i$, ${\bf C}\in C^\omega(I_i\times\T,SL(2,\R))$ and
\[
L({\bf C}(E,\cdot+i\e))=L({\bf C}(E,\cdot)),
\qquad E\in I_i,\quad |\e|<h,
\]
for some $h>0$. Thus the center cocycles are subcritical or
uniformly hyperbolic, by Theorem~\ref{t2.1}, and
$\hat\rho(E)=\rho({\bf C}(E,\cdot))$ is absolutely continuous
on $I_i$.\footnote{The argument in
\cite[Corollary 1.8]{avila1} applies here as well.
By \cite[Corollary 1.5]{avila1}, almost reducibility gives
rotations reducibility for a full-measure set of rotation
numbers. At such parameters, analytic dependence on $E$
implies that $\hat\rho$ is pointwise Lipschitz; at uniformly
hyperbolic parameters it is locally constant. Hence the image
under $\hat\rho$ of the non-Lipschitz parameters has Lebesgue
measure zero. Together with pointwise Lipschitzness on the
remaining parameters, this shows that $\hat\rho$ maps null
sets to null sets. By Corollary~\ref{c100},
$2\hat\rho(E)=N(E)+c_{I_i}$ with compatible continuous lifts,
so $\hat\rho$ is monotone. It is therefore absolutely continuous.}
By Corollary~\ref{c100},
\[
N(E)=2\hat\rho(E)-c_{I_i},\qquad E\in I_i.
\]

Thus $N$ is absolutely continuous on each $I_i$, completing
the proof. \qed

 \subsection{Proof of Theorem \ref{t3}}
Fix $E_0\in\Sigma_{v,\alpha}\cap I_i$. By almost reducibility
and \cite[Lemmas 2.2--2.3]{ajsm}, on a neighborhood $J\ni E_0$
the center family is analytically conjugate to Schr\"odinger
cocycles at a fixed energy $e_*$, with potentials $w_E$ in one
sufficiently small analytic ball on a common strip. Fix the first conjugacy at $E_0$. The iteration in Lemma~2.3 of \cite{ajsm} converges uniformly and analytically in $E$. Analyticity gives $\|w_E-w_{E'}\|_\infty\leq K|E-E'|$ locally.
 \cite{lyzz}
give a common $1/2$-H\"older constant $H$ for the IDS $N_{w_E}$.
Finite-volume min--max at $e_*\pm\|w_E-w_{E'}\|_\infty$ gives
\[
|N_{w_E}(e_*)-N_{w_{E'}}(e_*)|
\leq H\|w_E-w_{E'}\|_\infty^{1/2}
\leq H\sqrt K\,|E-E'|^{1/2}.
\]
The conjugacies have constant degree on $J$, so their
contributions cancel in rotation-number differences.
The Schr\"odinger rotation--IDS identity and Corollary~\ref{c100}
thus give a uniform $1/2$-H\"older estimate for $N$ on $J$.
The same argument applies to the original Schr\"odinger family
near subcritical energies. Compactness of the spectrum and
constancy of $N$ on the gaps then give global
$1/2$-H\"older continuity.\qed

 \section{Reducibility to localization: proof of Theorem \ref{t1}}
 In this section, we prove localization for Type I operators from reducibility of the dual center.
 The dual-reducibility-to-localization argument
 was pioneered by Avila--You--Zhou \cite{ayz} for the almost Mathieu
 operator, generalized in \cite{jk2}, and refined to the arithmetic
 level in \cite{gy,gyzh2}.

Those completeness arguments use the symmetry of the two phases associated with
an energy, $\pm \theta$ modulo their rotation orbits, thus depend on the reflection symmetry. 
However, these completeness arguments rely  on the $SL(2,
\mathbb{R})$ structure of the dual cocycle.
This creates the same two fundamental obstructions as
before for general Type I operators:

\begin{enumerate}
    \item \textbf{The (familiar) asymmetry obstruction ($v$ is not
        even):} For non-even potentials, the dual cocycle lies in
      $Sp(2, \mathbb{C})$. $\rho_1 \neq -\rho_2$ and the symmetry of the two phases associated with
an energy, $\pm \theta$ modulo their rotation orbits, crucial for the completeness arguments of \cite{ayz, gy, gyzh2} is lost. 
We develop new analytic tools to prove localization without relying on symmetric phases.
\item \textbf{The (familiar) infinite-range obstruction
($v$ is not a polynomial):}
For general analytic potentials, the dual operator has infinite range
and admits no classical transfer-matrix dynamics. We apply reducibility to the intrinsic two-dimensional
center constructed above. Together with the scalar cohomological equation, this produces analytic Bloch waves for the infinite-range dual. Their 
Fourier coefficients are exponentially decaying eigenfunctions
of the original Schr\"odinger operator. 
    
\end{enumerate}
\subsection{Nonperturbative reducibility}
Suppose that $A\in C^\omega(\T, SL(2,\R))$ admits a holomorphic
extension to $\{|\Im \theta|<h\}$.  Recall that the cocycle $(\alpha,A)$ is
said to be almost reducible if for some $0<h'<h$ there exists a
sequence $B_n\in C_{h'}^\omega(\T,PSL(2,R))$ such that
$B^{-1}_n(\theta+\alpha)A(\theta)B_n(\theta)$ converges to constant
uniformly in $|\Im\theta|<h'$. We have

\begin{Theorem}[\cite{avila1,avila2}] \label{arc-conjecture}
Any subcritical $(\alpha,A)$ with $\alpha\in \R\backslash\Q$, $A\in C^\omega(\T,SL(2,\R))$, is almost reducible.
\end{Theorem}
For every $\tau>1$ and $\gamma>0$, we define
\begin{equation}\label{newqq}
\Theta^\tau_\gamma:=\left\{\theta\in\T:\|2\theta+k\alpha\|_{\R/\Z}\geq
  \frac{\gamma}{(|k|+1)^\tau},k\in\Z\right\}; \Theta:=\cup_{\tau>1,\gamma>0}\Theta^\tau_\gamma.
\end{equation}
The quantitative estimates of \cite{ayz,gyzh2} on the subcritical strip give the following result.
\begin{Theorem}[\cite{ayz,gyzh2}]\label{glored}
Assume $\alpha\in \R\backslash\Q$, $h>\beta(\alpha)/2\pi$, $A\in
C_h^\omega(\T,SL(2,\R))$ with $L(\alpha,A(\cdot+i\e))=0$ for $|\e|<h$,
and $\rho(\alpha,A)\in \Theta.$ Then $(\alpha,A)$ is analytically reducible to a constant rotation. 
\end{Theorem}

\begin{Corollary} \label{c5}
Assume $\alpha\in \R\backslash\Q$, $v\in C^\omega(\T,\R)$, $\omega(E)=1$ with $L(E)>\beta(\alpha)$ and
$\hat{\rho}(E)\in \Theta.$ Then there exist  linearly independent $F_1,F_2\in C^\omega(\T,\C^{\Z})$ such that
$$
\widehat{H}_{v,\alpha,\theta}F(\theta)=EF(\theta),\ \ F(\theta)=(F_1(\theta), F_2(\theta));
$$
$$
TF(\theta)=F(\theta+\alpha)e^{2\pi i\hat{\phi_E}(0)}R_{\hat{\rho}(E)}.
$$
where $\hat{\phi}_E(0)=\int_\T \phi_E(\theta)d\theta,$  and
$R_{\hat{\rho}(E)}\in SO(2,\R)$ is a rotation by angle $2\pi \hat{\rho}(E).$
\end{Corollary}
\begin{pf}
By Theorems \ref{hss} and \ref{theorem-main1}, on any strip of width $h$ with $\beta(\alpha)<2\pi h<L(E)$,  There exist 
 \(
\phi\in C^\omega(\mathbb T,\mathbb R),
{\bf C}\in C^\omega(\mathbb T,\SL(2,\mathbb R)), (U,V)\in C^\omega(\T,\C^\Z),
\)
such that 
\begin{enumerate}
 \item $\widehat{H}_{v,\alpha,\theta}(U(\theta), V(\theta))=E(U(\theta),V(\theta))$;
 \item  $(TU(\theta),TV(\theta))=(U(\theta+\alpha),V(\theta+\alpha))e^{2\pi i\phi(\theta)}{\bf C}(\theta)$ ;
\item $L(\alpha,{\bf C}(\cdot+i\e))=0,\ \ \forall |\e|<h.$
\end{enumerate}

On the other hand, by Theorem \ref{glored}, multiplying a lifted projective reducing conjugacy of degree $d$ on the left by $R_{-d\theta/2}$ gives a periodic degree-zero conjugacy $B\in C^\omega(\T,SL(2,\R))$  such that 
$$
B(\theta+\alpha){\bf C}(\theta)B(\theta)^{-1}=R_{\hat{\rho}(E)}
$$
where $\hat{\rho}(E)=\rho({\bf C)}$. We define $\psi(\theta)=\sum_{k\in\Z\backslash\{0\}}\hat{\psi}(k)e^{2\pi i k\theta}$ where
\begin{equation}\label{new1010}
\hat{\psi}(k)=\frac{\hat{\phi}(k)}{e^{2\pi ik\alpha}-1}.
\end{equation}
Since $2\pi h>\beta(\alpha)$, it is easy to check that $\psi\in C^\omega(\T,\R)$ and $\psi(\theta+\alpha)-\psi(\theta)=\phi(\theta)-\hat{\phi}(0)$.

Finally, we define
\begin{equation}\label{new20}
F(\theta)=(U(\theta),V(\theta))B^{-1}(\theta)e^{2\pi i\psi(\theta)}.
\end{equation}
One can check that  
\begin{equation}\label{new21}
TF(\theta)=F(\theta+\alpha)e^{2\pi i\hat{\phi}(0)}R_{\hat{\rho}(E)}.
\end{equation}
\end{pf}
\subsection{Nonperturbative localization} 
Let
$$
\begin{gathered}
\Sigma_{v,\alpha}^\beta=\{E\in \Sigma_{v,\alpha}: \bar\omega(E)=1,\ L(E)>\beta(\alpha)\geq0\}.
\end{gathered}
$$
Fix local frames with $\beta(\alpha)<2\pi h<\inf_E L(E)$ and continuous lifts on a countable Borel partition subordinate to the intervals of Corollary~\ref{c100}. For every $\tau>1$ and $\gamma>0$, we define 
$$
\mathcal{E}^\tau_\gamma=\Sigma^\beta_{v,\alpha}\cap\{E\in\R: \hat{\rho}(E)\in \Theta_\gamma^\tau\}.
$$
For each $E\in\mathcal E_\gamma^\tau$, choose the analytic dual solutions
$F_1,F_2$ from the Corollary \ref{c5}. We define normalized sequences
$u_E,v_E\in \ell^2(\Z)$ by
\begin{equation}\label{def1}
u_{E}(n)=\frac{\hat{f}(n)}{\|f\|_{L^2}}=\frac{\int_{\T}f(\theta)e^{-2\pi in\theta}d\theta}{\|f\|_{L^2}},\ \ v_{E}(n)=\frac{\hat{g}(n)}{\|g\|_{L^2}}=\frac{\int_{\T}g(\theta)e^{-2\pi in\theta}d\theta}{\|g\|_{L^2}},
\end{equation}
where 
\begin{align}\label{gz25}
(f_1(\theta),f_2(\theta))=(F_1(0,\theta),F_2(0,\theta)).
\end{align}
\begin{align}\label{gz26}
f(\theta)=\frac{f_{1}(\theta)-if_{2}(\theta)}{\sqrt{2}},\ \ g(\theta)=\frac{f_{1}(\theta)+if_{2}(\theta)}{\sqrt{2}}.
\end{align}

Both $f$ and $g$ are nonzero. Indeed,  if either vanished
identically, then by \eqref{new21} every coordinate of  the
corresponding combination $F_1-iF_2$ or $F_1+iF_2$ would vanish
identically, contradicting the linear independence of $F_1,F_2$.
Thus the normalizations in \eqref{def1} are well defined.

\begin{Theorem}\label{eigen}
We have that $\{u_E(n)\}$ is an eigenfunction of
$H_{v,\alpha,\rho_1(E)}$ and $\{v_E(n)\}$ is an eigenfunction of
$H_{v,\alpha,\rho_2(E)}$, both with the eigenvalue $E$.
\end{Theorem}
\begin{pf}
By Corollary \ref{c5}, \eqref{gz25} and \eqref{gz26}, one has for $j=1,2$,
\begin{equation}\label{e41}
\sum\limits_{k=-\infty}^\infty \hat{v}_{k}e^{2\pi ik\rho_j(E)}h_{j}(\theta+k\alpha)+(2\cos2\pi(\theta)-E)h_{j}(\theta)=0,
\end{equation}
where $h_1=f$ and $h_2=g$. Taking the Fourier transform of \eqref{e41}, we get
\begin{equation*}
\sum\limits_{k=-\infty}^\infty \hat{v}_{k}e^{2\pi ik(\rho_1(E)+n\alpha)}\hat{f}(n)+\hat{f}(n-1)+\hat{f}(n+1)=E\hat{f}(n),
\end{equation*}
\begin{equation*}
\sum\limits_{k=-\infty}^\infty \hat{v}_{k}e^{2\pi ik(\rho_2(E)+n\alpha)}\hat{g}(n)+\hat{g}(n-1)+\hat{g}(n+1)=E\hat{g}(n).
\end{equation*}
Thus $\{u_E(n)\}$ is an eigenfunction of $H_{v,\alpha,\rho_1(E)}$ and
$\{v_E(n)\}$ is an eigenfunction of $H_{v,\alpha,\rho_2(E)}$,
corresponding to the eigenvalue $E$. Analyticity of $f,g$ gives exponential decay.

\end{pf}

For any fixed $\tau>1$, let $\mathcal{E}= \cup_\gamma \mathcal{E}^\tau_\gamma$ and $E^j(\theta)=\{E\in \mathcal{E}:\rho_j(E)=\theta (\mod\Z)\}$.

Note that by Theorem \ref{eigen}, $E^j(\theta)$ only contains eigenvalues of $H_{v,\alpha,\theta}$, so is a set that
contains at most countably many elements. Set $T\theta=\theta+\alpha.$
We also denote  $E^j_m(\theta):=E^j(T^m\theta)$, in particular, $E^j_0(\theta)=E^j(\theta)$. 
\begin{Definition}\label{Rm}
\label{defR}
$\nu_{\theta}$ is defined as:
$$
\nu_{\theta}=\sum\limits_{k\in\Z}\left(\sum\limits_{E\in E^1_k(\theta)}|u_{E}(-k)|^2+\sum\limits_{E\in E^2_k(\theta)}|v_{E}(-k)|^2\right)
$$
where $u_E,v_E$ are from \eqref{def1}.
\end{Definition}

\begin{Lemma}\label{lem:averaged-eigenfunction-mass}
We have
\(
 \int_{\mathbb T}\nu_\theta\,d\theta
 =\bigl|N(\Sigma_{v,\alpha}^{\beta})\bigr|
\)
where $|N(\Sigma_{v,\alpha}^{\beta})|$ denotes $\int_{\Sigma_{v,\alpha}^{\beta}}dN(E)$.
\end{Lemma}

\begin{pf}
Set \(S=\Sigma_{v,\alpha}^{\beta}\). By Corollary \ref{c100},
\begin{equation}\label{eq:lemma91-derivatives}
 \rho_1'(E)\geq0\geq\rho_2'(E),\ \ 
 N'(E)=\rho_1'(E)-\rho_2'(E)=2\hat\rho'(E)
 \ \ \text{for a.e. }E\in S.
\end{equation}
 By a countable disjoint Borel decomposition of $S$, it suffices
to work in a neighborhood where $\rho_1,\rho_2$ and $\widehat\rho$
admit absolutely continuous real lifts by the local argument of Section 7.3.

Since $E_k^j(\theta)=E^j(\theta+k\alpha)$, translation invariance
of Lebesgue measure and Tonelli's theorem give \footnote{The restrictions
$\rho_j|_{\mathcal E}$ are Borel and have at most countably many fibers by Theorem~\ref{eigen} (each fiber
$E^j(\theta)$ is contained in the point spectrum of $H_{v,\alpha,\theta}$). Moreover, by simplicity of the eigenvalues, 
$|u_E(k)|^2$ and $|v_E(k)|^2$ are the corresponding diagonal
spectral-projection weights, which are Borel functions of $E$. Thus the sums below are measurable.}
\begin{align*}
\int_{\mathbb T}\nu_\theta\,d\theta
&=
\sum_{k\in\mathbb Z}\int_{\mathbb T}
\left(
\sum_{E\in E^1(\theta+k\alpha)}|u_E(-k)|^2
+
\sum_{E\in E^2(\theta+k\alpha)}|v_E(-k)|^2
\right)d\theta\\
&=
\sum_{k\in\mathbb Z}\int_{\mathbb T}
\left(
\sum_{E\in E^1(\theta)}|u_E(-k)|^2
+
\sum_{E\in E^2(\theta)}|v_E(-k)|^2
\right)d\theta\\
&=
\int_{\mathbb T}
\left(
\sum_{E\in E^1(\theta)}\sum_{k\in\mathbb Z}|u_E(-k)|^2
+
\sum_{E\in E^2(\theta)}\sum_{k\in\mathbb Z}|v_E(-k)|^2
\right)d\theta\\
&=
\int_{\mathbb T}
\bigl(\#E^1(\theta)+\#E^2(\theta)\bigr)\,d\theta.
\end{align*}

In the last equality we used the normalization in \eqref{def1}. Within one branch, an energy belonging to both $E_k^j(\theta)$ and $E_\ell^j(\theta)$ forces $(k-\ell)\alpha\in\Z$, hence $k=\ell$. Across branches, an overlap would force $\rho_1(E)-\rho_2(E)=2\hat{\rho}(E)\in\Z+\alpha\Z$, excluded on $\mathcal{E}$. After shifting to phase $\theta$, orthogonality of different eigenfunctions gives $\nu_\theta\leq \mu_{\theta,\delta_0}^{pp}(S)$.

 The one-dimensional change-of-variables formula, counting
preimages with multiplicity, gives
\[
\int_{\mathbb T}\#E^j(\theta)\,d\theta
=
\int_{\mathcal E}|\rho_j'(E)|\,dE,
\qquad j=1,2.
\]

It follows from \eqref{eq:lemma91-derivatives} that
\begin{equation}\label{eq:lemma91-good-mass}
 \int_{\mathbb T}\nu_\theta\,d\theta
 =\int_{\mathcal E}\bigl(|\rho_1'(E)|+|\rho_2'(E)|\bigr)dE
 =\int_{\mathcal E}N'(E)\,dE.
\end{equation}
Set $\Theta^\tau=\bigcup_{\gamma>0}\Theta_\gamma^\tau$.
Since $\tau>1$, 
$|\mathbb T\backslash\Theta^\tau|=0$. Put
$B=\mathbb T\backslash\Theta^\tau$. Applying the  the same change-of-variables formula to
$\hat\rho$ restricted to $S\backslash\mathcal E$ gives
\[
 \int_{S\backslash\mathcal E}|\hat\rho'(E)|\,dE
 =\int_B\#\{E\in S\backslash\mathcal E:
       \hat\rho(E)=t\}\,dt=0.
\]
By \eqref{eq:lemma91-derivatives} and absolute continuity of the
density-of-states measure on $S$, we obtain
\[
 \int_{S\setminus\mathcal E}dN(E)
 =\int_{S\setminus\mathcal E}N'(E)\,dE=0.
\]
Combining this with \eqref{eq:lemma91-good-mass} proves
\[
 \int_{\mathbb T}\nu_\theta\,d\theta
 =\int_S N'(E)\,dE=|N(S)|.
\]
\end{pf}

\noindent {\bf Proof of Theorem \ref{t1}:}
Note that
$$
 \nu_{\theta}\leq |\mu^{pp}_{\theta,\delta_0}(\Sigma_{v,\alpha}^\beta)|\leq |\mu_{\theta,\delta_0}(\Sigma_{v,\alpha}^\beta)|,
$$
where $\mu_{\theta,\delta_0}$ is the spectral measure of $H_{v,\alpha,\theta}$ defined by
$$
\langle\delta_0,\chi_{B}(H_{v,\alpha,\theta})\delta_0\rangle=\int_{\R}\chi_{B} d\mu_{\theta,\delta_0}.
$$
Moreover, $\int_\T |\mu_{\theta,\delta_0}(\Sigma^\beta_{v,\alpha})|d\theta=|N(\Sigma_{v,\alpha}^\beta)|$.
It follows that $\nu_\theta=\mu_{\theta,\delta_0}(\Sigma^\beta_{v,\alpha})$ for a.e. $\theta$. Thus the spectral projection of $\delta_0$ on this set lies in the closed span of the exponentially decaying eigenfunctions. Shift covariance gives the same for $\delta_1$ for a.e. $\theta$. Taking intersection of these phases completes the proof (since $\{\delta_0,\delta_1\}$ is a cyclic set). \qed
\appendix
\renewcommand{\thesection}{\Alph{section}}
\section{Genericity of Type I for \texorpdfstring{$GL(1,\C)$}{GL(1,C)}
  cocycles}\label{app}
The Type~I condition is open in each $C^\omega_h$ with $h>\sup_{E\in\Sigma_{v,\alpha}}\e_1(E)$ \cite{gjy1}, and a natural
conjecture is that \emph{Type I is generic} i.e. that
\emph{Type~I energies are (open and) dense in the spectrum} for generic (i.e. open and
dense) analytic one-frequency Schrödinger operators. One piece of supporting evidence is that simplicity of every turning point that occurs is generic on each fixed closed annulus
 for analytic $GL(1,\mathbb C)$ cocycles---equivalently, analytic scalar functions on an
 annulus, nonvanishing on the unit circle---after Laurent-polynomial approximation, the condition that two roots have
the same modulus is contained in a proper real-algebraic subset of the fixed-degree
polynomial coefficient space and is therefore nongeneric.  Equivalently, the complementary simplicity condition is dense in
the corresponding analytic normed spaces. Certainly, the conjecture is a lot more challenging in
the non-commutative setting, where it is equivalent to
simplicity of the smallest nonnegative dual Lyapunov exponent, a problem of the
sort known to be quite difficult (e.g. \cite{gm}). However, this 1D
case may be viewed as a toy model for the conjectural density of
Type~I energies: failure of simplicity is expressed by real-analytic conditions
on the zero locations, and can be removed by arbitrarily small
analytic perturbations. We  provide detail of an elementary argument, showing directly that multiple zeros and coincident zero moduli can be removed by arbitrarily
small analytic perturbations.

Let
\(
A=\{z\in\C:r<|z|<R\}
\)
be an annulus with $0<r<R<\infty$. Denote by $C^\omega(\overline A)$ the space of functions
holomorphic on a neighborhood of $\overline A$, equipped with
the norm
$
\|f\|_{\overline A}=\sup_{z\in\overline A}|f(z)|.
$

Let \(\mathcal G\subset C^\omega(\overline A)\) be the set of functions \(f\) such that:
\begin{enumerate}
\item \(f\) has no zeros on the boundary circles \(\{|z|=r\}\cup\{|z|=R\}\),
\item all zeros of $f$ are simple and any two distinct zeros \(z_1,z_2\in A\) of \(f\) satisfy
\(
|z_1|\neq |z_2|.
\)
\end{enumerate}

\begin{Theorem}\label{thm:GL1-generic}
The set \(\mathcal G\) is open and dense in \(C^\omega(\overline A)\) with respect to the norm
\(\|\cdot\|_{\overline A}\). In particular, \(\mathcal G\) is generic.
\end{Theorem}

\begin{pf}
\(\mathcal G\) is open by an easy compactness argument. 
We now prove density. Let $f\in C^\omega(\overline A)$ and  $\varepsilon>0$ be given.
Since $f$ is holomorphic on an open neighborhood of $\overline A$, its Laurent series
$
f(z)=\sum_{n=-\infty}^{\infty} a_n z^n
$
converges uniformly on $\overline A$. Therefore for some $N,M<\infty,$
$
L(z):=\sum_{n=-N}^{M} a_n z^n = z^{-N}P(z),
$
where $P$ is a polynomial of degree at most $d:=N+M$, satisfies
$
\|f-L\|_{\overline A}<\frac{\varepsilon}{2}.
$

 We perturb the roots directly. If $P$ is a constant, take $P^*=P$ unless $P=0$, in which case take $0<|P^*|<\varepsilon r^N/2$. For nonconstant $P$,  write
\[
P(z)=a\prod_{j=1}^{d_0}(z-\alpha_j),\ \ a\neq0,\,d_0>0.
\]
 Choose $\beta_j$ arbitrarily close to $\alpha_j$ so that
their moduli are pairwise distinct and avoid $r$ and $R$. This can be
done successively: at each step only finitely many circles are excluded.

Set
\(
P^*(z)=a\prod_{j=1}^{d_0}(z-\beta_j).
\)
The coefficients of this product depend continuously on the roots.
Consequently $P^*\to P$ uniformly on $\overline A$ as
$\beta_j\to\alpha_j$. Since $z^{-N}$ is bounded on $\overline A$,
the perturbations may be chosen so that
\(
\|z^{-N}(P-P^*)\|_{\overline A}<\frac{\varepsilon}{2}.
\)
 In all cases put $g(z)=z^{-N}P^*(z)$. Its zeros in $A$ are simple,
have pairwise distinct moduli, and avoid $\partial A$. Thus $g\in\mathcal G$ and $\|f-g\|<
\varepsilon.$

\end{pf}

\section*{Acknowledgements}
L. Ge was partially supported by NSFC grant (12371185) and the Fundamental Research Funds for the Central Universities (the start-up fund), Peking University. SJ's work was supported by NSF DMS-2052899,
DMS-2155211, and Simons 896624. She is also grateful to School of
Mathematics at Georgia Institute of Technology and UC Irvine where  parts of this
work were done. We would like to thank A. Avila, J. You and Q, Zhou
for useful discussions. SJ is grateful to I. Spitkovsky for a tutorial
on factorization maps. We are also grateful to X. Li and Q. Zhou for
their careful
reading of the previous version that has prompted a significant improvement.
 
\end{document}